\input amstex\documentstyle{amsppt} 
\pagewidth{12.5cm}\pageheight{19cm}\magnification\magstep1
\topmatter
\title $\bold Z/m$-graded Lie algebras and perverse sheaves, I\endtitle
\author George Lusztig and Zhiwei Yun\endauthor
\address{Department of Mathematics, M.I.T., Cambridge, MA 02139,
Department of Mathematics, Stanford University, Stanford, CA 94305}\endaddress
\thanks{G.L. supported by NSF grant DMS-1566618; Z.Y. supported by NSF 
grant DMS-1302071 and the Packard Foundation.}\endthanks
\abstract We give a block decomposition of the equivariant derived category arising from a cyclically graded Lie algebra. This generalizes  certain aspects of the generalized Springer correspondence to the graded setting. \endabstract
\endtopmatter   
\document

       \define\tfl{\ti{\fl}}   \define\tfq{\ti{\fq}}

     \define\bco{\bar{\co}}

\define\mpb{\medpagebreak}

 \define\hfp{\hat{\fp}} \define\hfl{\hat{\fl}} \define\hfu{\hat{\fu}}

    \define\hP{\hat P}   \define\hz{\hat z}

\define\frl{\forall}

\define\si{\sim}
\define\wt{\widetilde}
\define\sqc{\sqcup}

\define\ovsc{\overset\cir\to}
\define\qua{\quad}

\define\hL{\hat L}

\define\bX{\bar X}

\define\lb{\linebreak}

\define\op{\oplus}
   
\define\part{\partial}
\define\emp{\emptyset}
\define\imp{\implies}
\define\ra{\rangle}
\define\n{\notin}
\define\iy{\infty}
\define\m{\mapsto}
\define\do{\dots}
\define\la{\langle}

\define\lra{\leftrightarrow}

\define\sm{\smallmatrix}
\define\esm{\endsmallmatrix}
\define\sub{\subset}    

\define\T{\times}
\define\ti{\tilde}
\define\nl{\newline}
\redefine\i{^{-1}}
\define\fra{\frac}
\define\un{\underline}
\define\ov{\overline}
\define\ot{\otimes}
\define\bbq{\bar{\QQ}_l}

\define\ad{\text{\rm ad}}
\define\Ad{\text{\rm Ad}}
\define\Hom{\text{\rm Hom}}

\define\Aut{\text{\rm Aut}}
\define\Ind{\text{\rm Ind}}    \define\tInd{\wt{\Ind}}

\define\ind{\text{\rm ind}}
\define\Res{\text{\rm Res}}    \define\tRes{\wt{\Res}}
\define\res{\text{\rm res}}

\define\supp{\text{\rm supp}}

\define\a{\alpha}
\redefine\b{\beta}
\redefine\c{\chi}
\define\g{\gamma}
\redefine\d{\delta}
\define\e{\epsilon}
\define\et{\eta}
\define\io{\iota}
\redefine\o{\omega}
\define\p{\pi}
\define\ph{\phi}
\define\ps{\psi}
\define\r{\rho}
\define\s{\sigma}
\redefine\t{\tau}
\define\th{\theta}
\define\k{\kappa}
\redefine\l{\lambda}
\define\z{\zeta}
\define\x{\xi}

\define\vt{\vartheta}

\redefine\D{\Delta}
\define\Om{\Omega}

\define\Ph{\Phi}
\define\Ps{\Psi}

\define\kk{\bold k}

\define\CC{\bold C}
\define\DD{\bold D}

\define\FF{\bold F}

\define\NN{\bold N}

\define\QQ{\bold Q}

\define\ZZ{\bold Z}

\define\ca{\Cal A}
\define\cb{\Cal B}
\define\cc{\Cal C}
\define\cd{\Cal D}

\define\cf{\Cal F}

\define\ch{\Cal H}
\define\ci{\Cal I}

\define\ck{\Cal K}
\define\cl{\Cal L}

\define\co{\Cal O}
\define\cp{\Cal P}
\define\cq{\Cal Q}
\define\car{\Cal R}

\define\ct{\Cal T}

\define\cz{\Cal Z}
\define\cx{\Cal X}

\define\fg{\frak g}
\define\fh{\frak h}

\define\fk{\frak k}
\define\fl{\frak l}
\define\fm{\frak m}

\define\fp{\frak p}
\define\fq{\frak q}

\define\fs{\frak s}

\define\fu{\frak u}

\define\fz{\frak z}

\define\fG{\frak G}

\define\fL{\frak L}
\define\fM{\frak M}

\define\fP{\frak P}

\define\fR{\frak R}
\define\fS{\frak S}
\define\fT{\frak T}

\define\tz{\ti z}

\define\tA{\ti A}

\define\tC{\ti C}

\define\tE{\ti E}

\define\tL{\ti L}

\define\tQ{\ti Q}

\define\sha{\sharp}

\define\bul{\bullet}

\define\cir{\circ}

\define\KL{KL}
\define\KO{Ko}
\define\INTERS{L1}
\define\CSI{L2}
\define\CSII{L3}
\define\GRA{L4}
\define\ANTIORB{L5}
\define\RIR{RR}
\define\ST{St}
\define\VI{Vi}

\define\da{\delta}

\head Contents\endhead
Introduction.

1. Recollections on $\bold Z$-graded Lie algebras.

2. $\bold Z/m$-gradings and $\epsilon$-spirals.

3. Admissible systems.

4. Spiral induction.

5. Study of a pair of spirals.

6. Spiral restriction.

7. The categories ${\Cal Q}(\frak g_\delta)$, ${\Cal Q}'(\frak g_\delta)$.

8. Monomial and quasi-monomial objects.

9. Examples.

\head Introduction\endhead
\subhead 0.1\endsubhead
Let $\kk$ be an algebraically closed field of characteristic $p\ge0$. We fix 
an integer $m>0$ such that $m<p$ whenever $p>0$ and we write $\ZZ/m$ instead 
of $\ZZ/m\ZZ$. For $n\in\ZZ$, let $\un n$ denote the image of $n$ in  $\ZZ/m$.
 
We also fix $G$, a semisimple simply connected algebraic group 
over $\kk$ and a $\ZZ/m$-grading $\fg=\op_{i\in\ZZ/m}\fg_i$ (see 0.11) for 
the Lie algebra $\fg$ of $G$; we shall assume that either $p=0$ or that $p$ is
so large relative to $G$, that we can operate with $\fg$ as if $p$ was $0$.

For any integer $d$ invertible in $\kk$ let $\mu_d=\{z\in\kk^*;z^d=1\}$. The
$\ZZ/m$-grading on $\fg$ is the same as an action of $\mu_m$ on $G$ or a 
homomorphism $\ti\vt:\mu_m@>>>\Aut(G)$. ($\ti\vt$ induces a homomorphism
$\ti\th:\mu_m@>>>\Aut(\fg)$ and for $i\in\ZZ/m$ we have
$\fg_i=\{x\in\fg;\ti\th(z)x=z^ix\qua\frl z\in\mu_m\}$.) Let 
$G_{\un0}=\{g\in G;g\ti\vt(z)=\ti\vt(z)g\qua\frl z\in\mu_m\}$, a connected 
reductive subgroup of $G$ with Lie algebra $\fg_{\un0}$. 
For any $i\in\ZZ/m$, the $\Ad$-action of $G_{\un0}$ on $\fg$ leaves stable 
$\fg_i$ and its closed subset $\fg_i^{nil}:=\fg_i\cap\fg^{nil}$. (Here 
$\fg^{nil}$ is the variety of nilpotent elements in $\fg$.)

We are interested in studying the equivariant derived categories (see 0.11)
$\cd_{G_{\un0}}(\fg_i)$, $\cd_{G_{\un0}}(\fg^{nil}_i)$.
More specifically we would like to classify
$G_{\un0}$-equivariant simple perverse sheaves with
support in $\fg^{nil}_i$ and (in the case where $p>0$) their 
Fourier-Deligne transform. The simple perverse sheaves in 
$\cd_{G_{\un0}}(\fg^{nil}_i)$ are in bijection with pairs $(\co,\cl)$
where $\co$ is a nilpotent $G_{\un0}$-orbit in $\fg_i$ and $\cl$ is
(the isomorphism class of) an irreducible $G_{\un0}$-equivariant local system
on $\co$. (The pair $(\co,\cl)$ gives rise to the simple perverse sheaf $P$
with support equal to the closure of $\co$ and with $P|_\co=\cl[\dim\co]$.)
We denote the set of such $(\co,\cl)$ by $\ci(\fg_i)$. This is a finite
set, since the $G_{\un0}$-action on $\fg_i$ has only finitely many orbits.
Alternatively, if we choose $e\in\co$, then the local system $\cl$ corresponds
to an irreducible representation of $\p_0(G_{\un0}(e))$ (see 0.11) where 
$G_{\un0}(e)$ is the stabilizer of $e$ under $G_{\un0}$.

There are many $\ZZ/m$-graded Lie algebras which appear in nature.

\subhead 0.2\endsubhead
In this subsection we assume that $m=2$ and $\kk=\CC$. Then the 
$\ZZ/2$-grading $\fg=\fk\op\fp$ (with $\fk=\fg_{\un0}$, $\fp=\fg_{\un1}$)  
has been extensively studied in connection with the theory of symmetric
spaces and the representation theory of real semisimple groups. In particular,
the nilpotent $G_{\un0}$-orbits on $\fp$ are known to
be in bijection with the nilpotent orbits in the Lie algebra of a real form
of $G$ determined by the $\ZZ/2$-grading (Kostant and Sekiguchi).

\subhead 0.3\endsubhead
Another class of examples comes from cyclic quivers.  In this subsection
we assume that $m\ge2$. 
We consider the simplest
such example where $V$ is a $\kk$-vector space equipped with a $\ZZ/m$-grading
$V=\op_{i\in\ZZ/m}V_i$ (see 0.11) 
and $G=SL(V)$ with the $\ZZ/m$-grading given by
$$\fg_i=\{T\in\fg=\fs\fl(V);T(V_j)\sub V_{j+i}\qua\frl j\in\ZZ/m\}.$$
In this case we have $G_{\un0}=S(\prod_{i\in\ZZ/m}GL(V_i))$, the intersection 
of $SL(V)$ with the Levi subgroup $\prod_{i}GL(V_i)$ of a parabolic subgroup
of $GL(V)$. 
The subspace $\fg_{\un1}$
is $\op_{i\in\ZZ/m}\Hom(V_i,V_{i+1})$. We may consider a quiver $Q$
with $m$ vertices indexed by $\ZZ/m$ and an arrow $i\m i+1$ for each 
$i\in\ZZ/m$,
$$\CD
V_2@<<<V_1\\
@VVV   @AAA\\
\do @.   V_0\\
@VVV   @AAA\\
V_{m-2}@>>>V_{m-1}\\
\endCD$$
Then $\fg_{\un1}$ 
is the space of representations of $Q$ where we put $V_i$ at the
vertex $i$. 

More generally, if $G$ is a classical group, then the $G_{\un0}$-action on
$\fg_{\un1}$ can be interpreted in terms of a cyclic quiver with some extra
structure (see 9.5 for the case where $G$ is a symplectic group).  

\subhead 0.4\endsubhead  
In this subsection we forget the $\ZZ/m$-grading. Instead of the action of
$G_{\un0}$ on $\fg_i$ and $\fg_i^{nil}$ we consider the adjoint 
action of $G$ on $\fg$ and on $\fg^{nil}$. 
Let $\ci(\fg)$ be the set of pairs $(\co,\cl)$ where $\co$ is a $G$-orbit on
$\fg^{nil}$ and $\cl$ is an irreducible $G$-equivariant local system on $\co$
(up to isomorphism). 
From the results on the generalized Springer theory in \cite{\INTERS} we have
a canonical decomposition
$$\ci(\fg)=\sqc_{(L,C)}\ci(\fg)_{(L,C)}\tag a$$
where $(L,C)$ runs over the $G$-conjugacy classes of data $L,C$ with $L$ a 
Levi subgroup of a parabolic subgroup of $G$ and $C$ an $L$-equivariant
cuspidal perverse sheaf on the nilpotent cone of the Lie algebra of $L$.
(Actually, the results of \cite{\INTERS} are stated for unipotent elements in
$G$ instead of nilpotent elements in $\fg$.) We call (a) the {\it block 
decomposition} of $\ci(\fg)$.  

Let $P(\fg^{nil})$ be the subcategory of $\cd(\fg^{nil})$ consisting of 
complexes whose perverse cohomology sheaves are $G$-equivariant. Using (a) 
and \cite{\CSII, (7.3.1)}, we see that we have a direct sum 
decomposition
$$P(\fg^{nil})=\op_{(L,C)}P(\fg^{nil})_{(L,C)},\tag b$$
where $(L,C)$ is as in (a).
We call (b) the {\it block decomposition} of $P(\fg^{nil})$. 
In \cite{\RIR} it is shown that the following variant of (b) holds: we have
a direct sum decomposition 
$$\cd_G(\fg^{nil})=\op_{(L,C)}\cd_G(\fg^{nil})_{(L,C)}\tag c$$
where $(L,C)$ is as in (a). We call (c) the {\it block decomposition} of 
$\cd_G(\fg^{nil})$. 

In this paper we find a $\ZZ/m$-graded analogue of this (ungraded) block 
decomposition.

\subhead 0.5\endsubhead
We fix $\z$, a primitive $m$-th root of $1$ in $\kk$ and we set 
$\vt=\ti\vt(\z):G@>>>G$, $\th=\ti\th(\z):\fg@>>>\fg$.
Then for $i\in\ZZ/m$ we have $\fg_i=\{x\in\fg;\th(x)=\z^ix\}$.

Let $\et\in\ZZ-\{0\}$. We consider systems 
$(M,\fm_*,\tC)$ where $M=\{g\in G;\Ad(\t)\vt g=g\}$ for some semisimple
element of finite order $\t\in G_{\un0}$, $\fm_*=\{\fm_N\}_{N\in\ZZ}$ is a 
$\ZZ$-grading of the Lie algebra $\fm$ of $M$ (see 0.11) such that 
$\fm_N\sub\fg_{\un N}$ for all $N$, $M_0$ is the closed connected
subgroup of $M$ with Lie algebra $\fm_0$ and $\tC$ is an $M_0$-equivariant 
cuspidal perverse sheaf on $\fm_\et$. We will review the notion of 
$M_0$-equivariant cuspidal perverse sheaf (already defined in \cite{\GRA}) on 
$\fm_\et$ in 1.2. 
Such a system $(M,\fm_*,\tC)$ is said to be {\it admissible}
if a certain technical condition involving the group of components of the
centre of $M$ is satisfied (see 3.1). 

Let $\un\fT_\et$ be the set of admissible 
systems up to $G_{\un0}$-conjugacy. The following result is proved in 7.9.

\proclaim{Theorem 0.6} There is a canonical direct sum decomposition of
$\cd_{G_{\un0}}(\fg_{\un\et}^{nil})$ into full subcategories
$$\cd_{G_{\un0}}(\fg_{\un\et}^{nil})=\op_{(M,\fm_*,\tC)\in\un\fT_\et}
\cd_{G_{\un0}}(\fg_{\un\et}^{nil})_{(M,\fm_*,\tC)}$$
indexed by $\un\fT_\et$.
\endproclaim   

In particular, any simple perverse sheaf in $\cd_{G_{\un0}}(\fg_{\un\et}^{nil})$
belongs to a well defined block 
$\cd_{G_{\un0}}(\fg_{\un\et}^{nil})_{(M,\fm_*,\tC)}$.
This gives a map
$$\Psi:\ci(\fg_{\un\et})@>>>\un\fT_\et.$$
In fact, we will first establish the map $\Psi$ in 3.5 and then prove the theorem in 7.9, 
using a key calculation in proposition 6.4. 

Note that in the case where $m=1$, the theorem can be deduced from 0.4(a).
On the other hand, for large $m$, a $\ZZ/m$-grading on $\fg$ is the same as a
$\ZZ$-grading, so that in this case the theorem can be deduced from the
results of \cite{\GRA}. Thus, the result about block decomposition in this
paper generalizes results in \cite{\INTERS} and \cite{\GRA}.

\subhead 0.7\endsubhead
As an explicit example, let us consider the case where $G=SL_n(\kk),\et=1$. In the
ungraded case, blocks are in bijection with pairs $(d,\c)$ where $d$ is a
divisor of $n$ and $\c:\mu_d@>>>\bbq^*$ is a primitive character. (See
\cite{\INTERS}.) To $d$ we attach the subgroup $M=S(GL_d^{n/d})$ (a Levi
subgroup of a parabolic subgroup) and $\c$ gives a cuspidal perverse sheaf
$C_\c$ with support equal to the nilpotent cone of the Lie algebra of $M$.
Now in the $\ZZ/m$-graded case, we have $G=SL(V)$, $V=\op_{i\in \ZZ/m}V_i$ as in 0.3, and
we identify $\fg_{\un1}$ 
with $\op_i\Hom(V_i,V_{i+1})$. In this case, the set 
of blocks $\un\fT_1$ has a similar explicit description. We have a natural
bijection 
$$\un\fT_1\lra\{(d,f,\c)\}/\si.\tag a$$
Here the right hand side is the set of equivalence classes of triples 
$(d,f,\c)$ where $(d,\c)$ is as in the ungraded case and 
$f:\{1,2,\do,n/d\}@>>>\ZZ/m$ is a map such that
$$\sha\{(b,y)\in\ZZ\T\ZZ;1\le b\le n/d,0\le y\le d-1,f(b)+\un y=i\}=\dim V_i
\tag b$$ 
for all $i\in\ZZ/m$. Two triples $(d,f,\c)$ and $(d',f',\c')$ are equivalent
if and only if $d=d',\c=\c'$ and $f'$ is obtained from $f$ by composition with
a permutation of $\{1,2,\do,n/d\}$.

\subhead 0.8\endsubhead
In the ungraded case, the objects in the block $\cd_G(\fg^{nil})_{(L,C)}$ are
obtained from $C$ via parabolic induction (and decomposition) through any
parabolic subgroup $P$ of $G$ containing $L$ as a Levi subgroup. In the
$\ZZ/m$-graded case, each block
$\cd_{G_{\un0}}(\fg_{\un\et}^{nil})_{(M,\fm_*,\tC)}$
is constructed from a certain induction procedure which we call {\it spiral
induction}, see Section 4. 
Instead of a parabolic subgroup of $G$ compatible with the
$\ZZ/m$-grading on $\fg$, we introduce the notion of a spiral $\fp_*$ 
which is a sequence of subspaces $\fp_N\sub\fg_{\un N}$, one for each $N\in\ZZ$,
satisfying certain conditions (see Section 2).  
It turns out that spirals are the correct
analogues of parabolic subalgebras in the $\ZZ/m$-graded case. In fact there are
two kinds of spiral inductions, one giving objects in 
$\cd_{G_{\un0}}(\fg_{\un\et}^{nil})$ and the other giving (assuming that $p>0$)
Fourier-Deligne transforms of objects
in $\cd_{G_{\un0}}(\fg_{-\un\et}^{nil})$. The latter may be viewed as an analogue
of character sheaves in the $\ZZ/m$-graded setting.

\subhead 0.9\endsubhead
We now discuss the contents of the various sections.
Many arguments in this paper rely on results from \cite{\GRA} concerning
$\ZZ$-graded Lie algebras; in Section 1 we review some results from 
\cite{\GRA} that we need. In Section 2 we introduce the $\e$-spirals attached
to a $\ZZ/m$-graded Lie algebra and their splittings; the analogous concepts
in the $\ZZ$-graded cases are the parabolic subalgebras compatible with the
$\ZZ$-grading and their Levi subalgebras compatible with the $\ZZ$-grading.
We also attach a canonical spiral to any element of $\fg_{\un\et}^{nil}$ which
plays a crucial role in the arguments of this paper.
In Section 3 we introduce the admissible systems, which eventually will be
used to index the blocks in $\cd_{G_{\un0}}(\fg_{\un\et}^{nil})$. 
In Section 4 we introduce the operation of spiral induction which is our
main tool in the study of $\cd_{G_{\un0}}(\fg_{\un\et}^{nil})$. 
In Sections 5 and 6 we calculate explicitly the $Hom$ space between two
spiral inductions, generalizing to the $\ZZ/m$-graded case a result from
\cite{\GRA}. This is used in Section 7 to prove Theorem 0.6.
In Section 8 we introduce monomial and quasi-monomial complexes on
$\fg_{\un\et}^{nil}$; we show that the monomial complexes (resp. quasi-monomial)
complexes generate the appropriate Grothendieck group over $\QQ(v)$
(resp. over $\ZZ[v,v\i]$) where $v$ is an indeterminate; this again
generalizes to the $\ZZ/m$-graded case a result from \cite{\GRA}. This result
is of the same type as that which says that the plus part of a quantized
enveloping algebra defined in terms of perverse sheaves is generated over
$\QQ(v)$ by monomials in the $E_i$ and over $\ZZ[v,v\i]$ by the products of
divided powers of the $E_i$ (which could be called quasi-monomials).
In Section 9 we discuss the examples where $G=SL(V)$ or $G=Sp(V)$; 
in these cases we describe the spirals and in the case of $G=SL(V)$ we describe
the blocks.

\subhead 0.10\endsubhead
It is known that, in the ungraded case, 
each block of $\cd_G(\fg^{nil})$ can be related to
the category of representations of a certain Weyl group; if $m$ is large, so 
that the $\ZZ/m$ grading of $\fg$ is a $\ZZ$-grading and 
$\fg_{\un\et}^{nil}=\fg_{\un\et}$, then each block of 
$\cd_{G_{\un0}}(\fg_{\un\et}^{nil})$ is related to the
category of representations of a certain graded affine Hecke algebra with
possibly unequal parameters. 
In fact, without assumptions on $m$, each block of
$\cd_{G_{\un0}}(\fg_{\un\et}^{nil})$ is related to a certain 
graded double affine Hecke algebra (corresponding to an affine Weyl group 
attached to the block) with possibly unequal parameters; this will be considered
in a sequel to this paper. We also plan to describe explicitly
the blocks in the case where $G$ is a classical group and relate them to 
cyclic quivers with extra structure.
The case of the symplectic group is partially carried out in 9.5-9.7.

\subhead 0.11. Notation \endsubhead 
All algebraic varieties are assumed to be over $\kk$; all 
algebraic groups are assumed to be affine. 
For any algebraic variety $X$ we denote by $\cd(X)$ the bounded derived category of 
$\bbq$-complexes on $X$. 
Let $D:\cd(X)@>>>\cd(X)$ be Verdier duality. 

For $K\in\cd(X)$ we denote by $\ch^nK$ the $n$-th cohomology sheaf of $K$ and 
by $\ch^n_xK$ the stalk of $\ch^nK$ at $x\in X$.

If $X'$ is a locally closed smooth irreducible subvariety of $X$ with closure
$\bX'$ and $\cl$ is an irreducible local system on $X'$ we denote by 
$\cl^\sha\in\cd(X)$ the intersection cohomolgy complex of $\bar X'$ with
coefficients in $\cl$, extended by $0$ on $X-\bX'$.

If $X$  has a given action of an algebraic group $H$ we denote by $\cd_H(X)$ 
the corresponding equivariant derived category.

If $H$ is an 
 algebraic group we denote by $H^0$ the identity component of $H$, by $\cz_H$
the centre of $H$. We set $\p_0(H)=H/H^0$.
Now assume that $H$ is connected. We denote by $\fL H$ 
the Lie algebra of $H$ and by $U_H$ the unipotent radical of $H$. Let $\fh=\fL H$.
If $\fh'$ is a Lie subalgebra of $\fh$ we write 
$e^{\fh'}\sub H$ for the closed connected subgroup of $H$ such that 
$\fL(e^{\fh'})=\fh'$, assuming that such a subgroup exists.

We shall often denote a collection $\{V_N;N\in\ZZ\}$ of vector spaces indexed 
by $N\in\ZZ$ by the symbol $V_*$.

If $V$ is a $\kk$-vector space, a $\ZZ$-grading on $V$ is a collection
of subspaces $V_*=\{V_k;k\in\ZZ\}$ such that $V=\op_{k\in\ZZ}V_k$; 
a $\ZZ/m$-grading 
on $V$ is a collection of subspaces $\{V_i;i\in\ZZ/m\}$ such that 
$V=\op_{i\in\ZZ/m}V_i$; a $\QQ$-grading on $V$ is a collection of subspaces 
$\{{}_\k V;\k\in\QQ\}$ such that $V=\op_{\k\in\QQ}({}_\k V)$.

A $\ZZ$-grading for the Lie algebra $\fh$ is a $\ZZ$-grading
$\fh_*=\{\fh_k;k\in\ZZ\}$ of $\fh$ as a vector space satisfying 
$[\fh_k,\fh_{k'}]\sub\fh_{k+k'}$ for all $k,k'\in\ZZ$; 
a $\ZZ/m$-grading for $\fh$ 
is a $\ZZ/m$-grading $\{\fh_i;i\in\ZZ/m\}$ of $\fh$ as a vector 
space satisfying $[\fh_i,\fh_{i'}]\sub\fh_{i+i'}$ for all $i,i'\in\ZZ/m$; 
a $\QQ$-grading for $\fh$ is a $\QQ$-grading
$\{{}_\k\fh;\k\in\QQ\}$ of $\fh$ as a vector space satisfying 
$[{}_\k\fh,{}_{\k'}\fh]\sub{}_{\k+\k'}\fh$ for all $\k,\k'\in\QQ$. 

Let $Y_H$ be the set of homomorphisms of algebraic groups $\kk^*@>>>H$. For
$\l\in Y_H$ and $c\in\ZZ$, we define $c\l\in Y_H$ by $(c\l)(t)=\l(t^c)$ for 
$t\in\kk^*$. We define an equivalence relation on $Y_H\T\ZZ_{>0}$ by
$(\l,r)\si(\l',r')$ whenever there exist $c,c'$ in $\ZZ_{>0}$ such that
$c\l=c'\l'$, $cr=c'r'$; the set of equivalence classes for this relation is 
denoted by $Y_{H,\QQ}$. Let $\l/r=(1/r)\l$ be the equivalence class of $(\l,r)$. Now 
$\l\m\l/1$ identifies $Y_H$ with a subset of $Y_{H,\QQ}$. For 
$\k\in\QQ,\mu\in Y_{H,\QQ}$ we define $\k\mu\in Y_{H,\QQ}$ by
$\k\mu=(k\l)/(k'r)$ where $k\in\ZZ,k'\in\ZZ_{>0},r\in\ZZ_{>0},\l\in Y_H$ are
such that $\k=k/k'$, $\mu=\l/r$; this is independent of the choices. In 
particular, we have $r\mu\in Y_H$ for some $r\in\ZZ_{>0}$.

Let $\l\in Y_H$. For $k\in\ZZ$ we set
$${}^\l_k\fh=\{x\in\fh;\Ad(\l(t))x=t^kx\qua\frl t\in\kk^*\}.$$
Note that $\{{}^\l_k\fh,k\in\ZZ\}$ is a $\ZZ$-grading of $\fh$.

Now let $\mu\in Y_{H,\QQ}$. For $\k\in\QQ$ we set
${}^\mu_\k\fh={}^{r\mu}_{r\k}\fh$
where $r\in\ZZ_{>0}$ is chosen so that
$r\mu\in Y_H$, $r\k\in\ZZ$. This is well defined (independent of the choice of
$r$). Note that $\{{}^\mu_\k\fh,\k\in\QQ\}$ is a $\QQ$-grading of $\fh$.

\subhead 0.12\endsubhead
Let $H$ be a connected algebraic group acting on an algebraic variety $X$ and
let $A,B$ be two $H$-equivariant semisimple complexes on $X$; let $j\in\ZZ$.
We define a finite dimensional $\bbq$-vector space
$\DD_j(X,H;A,B)$ as in \cite{\GRA, 1.7}. For the purpose of this paper, we can take the following formula as the definition of $\DD_j(X,H;A,B)$

(a)  $\DD_j(X,H;A,B)=\Hom_{\cd_H(X)}(A,D(B)[-j])^{*}$.
\nl
Let $d_j(X;A,B)=\dim\DD_j(X,H;A,B)$,
$\{A,B\}=\sum_{j\in\ZZ}d_j(X;A,B)v^{-j}\in\NN((v))$
where $v$ is an indeterminate.

If $A,B$ are $H$-equivariant simple perverse sheaves on $X$, then

$\{A,B\}\in1+v\NN[[v]]$ if $B\cong D(A)$

$\{A,B\}\in v\NN[[v]]$ if $B\not\cong D(A)$.
\nl
(See \cite{\GRA, 1.8(d)}.) 

For an algebraic variety $X$ we denote by $\r_X$ the map $X@>>>(\text{point})$.

Let $v$ be an indeterminate and let $\ca=\ZZ[v,v\i]$.
Let $\bar{}:\QQ(v)@>>>\QQ(v)$ be the field involution such that $\bar v=v\i$.
This restricts to a ring involution $\bar{}:\ca@>>>\ca$.

\define\doet{\dot\et}

For any $\et\in\ZZ-\{0\}$ we define $\doet=\et/|\et|\in\{1,-1\}$ where $|\et|$ is 
the absolute value of $\et$.

\head 1. Recollections on $\ZZ$-graded Lie algebras\endhead
In this section we recall notation and results from \cite{\GRA} that will be used in this paper.

\subhead 1.1\endsubhead
In this section we fix a connected reductive group $H$; let $\fh=\fL H$.

Let $J^H$ be the variety consisting of all triples $(e,h,f)\in\fh^3$ such that 
$[h,e]=2e,[h,f]=-2f,[e,f]=h$ (then $e,f$ are necessarily in $\fh^{nil}$). If 
$\ph=(e,h,f)\in J^H$, there is a unique homomorphism of algebraic groups 
$\ti\ph:SL_2(\kk)@>>>H$ such that the differential of $\ti\ph$ carries 
$\left(\sm 1&0\\0&0\esm\right),\left(\sm 1&0\\0&-1\esm\right), 
\left(\sm 0&0\\1&0\esm\right)$ to $e,h,f$ respectively; we then define 
$\io_\ph\in Y_H$ by $\io_\ph(t)=\ti\ph\left(\sm t&0\\0&t\i\esm\right)$.

\subhead 1.2\endsubhead
In the remainder of this section we assume that a $\ZZ$-grading $\fh_*$ for 
$\fh$ is given. Then there exists $\l\in Y_H$ and $r\in\ZZ_{>0}$ with 
$\fh_k={}^\l_{rk}\fh$ 
for all $k\in\ZZ$. (It follows that ${}^\l_\k\fh=0$ for all $\k\in\QQ-r\ZZ$.)

(In this paper we will often refer to results in \cite{\GRA}, even though,
strictly speaking, in \cite{\GRA} a stronger assumption on the $\ZZ$-grading 
of $\fh$ is made, namely that $r$ above can be taken to be $1$. Note that the 
results of \cite{\GRA} hold with the same proof when the stronger assumption 
is replaced by the present assumption.)

We have $\fh_k\sub\fh^{nil}$ for any $k\in\ZZ-\{0\}$. Note that $\fh_0$ is a 
Lie subalgebra of $\fh$ and that $H_0:=e^{\fh_0}\sub H$ is well defined and it
acts by the $\Ad$-action on each $\fh_k$. If $k\ne0$, this action has only 
finitely many orbits, see \cite{\GRA, 3.5}; we denote by $\ovsc\fh_k$ the 
unique open $H_0$-orbit in $\fh_k$.

Let $\et\in\ZZ-\{0\}$. 

(a) {\it We say that the $\ZZ$-grading $\fh_*$ of $\fh$ is {\it $\et$-rigid} 
if there exists $\io\in Y_H$ such that (i),(ii) below hold.

(i) ${}^\io_k\fh=\fh_{\et k/2}$ for any $k\in\ZZ$ such that $\et k/2\in\ZZ$ and
${}^\io_k\fh=0$ for any $k\in\ZZ$ such that $\et k/2\n\ZZ$;

(ii) $\io=\io_\ph$ for some $\ph=(e,h,f)\in J^H$ such that $e\in\ovsc\fh_\et$,
$h\in\fh_0$, $f\in\fh_{-\et}$.}
\nl
It follows that $2k'\in\et\ZZ$ whenever $\fh_{k'}\ne0$.
Note that $\io$ is unique if it exists, since, by (ii), $\io(\kk^*)$ is 
contained in the derived group of $H$.

We show:

(b) {\it In the setup of (a), let $\ph'=(e',h',f')\in J^H$ be such that
$e'\in\ovsc\fh_\et$, $h'\in\fh_0$, $f'\in\fh_{-\et}$. Let $\io'=\io_{\ph'}$.
Then $\io'=\io$.}
\nl
Let $\ph$ be as in (ii). Using \cite{\GRA, 3.3}, we can find
$g_0\in H_0$ such that $\Ad(g_0)$ carries $\ph$ to $\ph'$. It follows that
$\Ad(g_0)\io(t)=\io'(t)$ for any $t\in\kk^*$. For $k\in\ZZ$ such that $\et k/2\in\ZZ$ we have
$${}^{\io'}_k\fh=\Ad(g_0)({}^\io_k\fh)=\Ad(g_0)\fh_k=\fh_k;$$
for $k\in\ZZ$ such that $\et k/2\n\ZZ$ we have
$${}^{\io'}_k\fh=\Ad(g_0)({}^\io_k\fh)=0.$$
$${}^{\io'}_{2k\et}\fh=\Ad(g_0)({}^{\io}_{2k\et}\fh)=\Ad(g_0)\fh_k=\fh_k.$$
Thus $\io'$ satisfies the defining properties of $\io$ in (a). By uniqueness
we have $\io'=\io$ as required.

\mpb

Let $\ci(\fh_\et)$ be the set of all pairs $(\co,\cl)$ where $\co$ is an
$H_0$-orbit in $\fh_\et$ and $\cl$ is an $H_0$-equivariant irreducible local
system on $\fh_\et$ (up to isomorphism).

Let $\cq(\fh_\et)$ be the category of $\bbq$-complexes on $\fh_\et$ which are 
direct sums of shifts of simple $H_0$-equivariant perverse sheaves on 
$\fh_\et$. There are up to isomorphism only finitely many such simple perverse
sheaves; they form a set in bijection with $\ci(\fh_\et)$.

An $H_0$-equivariant perverse sheaf $A$ on $\fh_\et$ is said to
be {\it cuspidal} if there exists a nilpotent $H$-orbit $\cc$ in $\fh$ and an 
irreducible $H$-equivariant cuspidal local system $\cf$ on $\cc$ such that 
$\ovsc\fh_\et\sub\cc$ and $A|_{\ovsc\fh_\et}=\cf|_{\ovsc\fh_\et}[\dim\fh_\et]$. 
If such $(\cc,\cf)$ exists, it is unique, see \cite{\GRA, 4.2(c)}. 
Note that if $A$ is cuspidal then it is necessarily a simple perverse sheaf.

(c) {\it If there exists a cuspidal $H_0$-equivariant perverse sheaf $A$ on $\fh_\et$,
 the grading $\fh_*$ of $\fh$ is necessarily $\et$-rigid;
moreover, we have $A|_{\fh_\et-\ovsc\fh_\et}=0$.} 
\nl
(See \cite{\GRA, 4.4(a), 4.4(b)}.) 

In the setup of (c), the element $\io\in Y_H$ provided by (a) is known to satisfy 

(d) ${}^\io_k\fh=0$ unless $k\in2\ZZ$; 
\nl
we deduce that:

(e) {\it If $k'\in\ZZ$ and $\fh_{k'}\ne0$ then $k'/\et\in\ZZ$.} 

\subhead 1.3. Parabolic induction\endsubhead
In the setup of 1.2 assume that $P$ is a parabolic subgroup of $H$ with 
$\fp:=\fL P$ satisfying $\fp=\op_{k\in\ZZ}\fp_k$ where $\fp_k=\fp\cap\fh_k$. 
We set $U=U_P,L=P/U$, $\fu=\fL U,\fl=\fL L=\fp/\fu$. We have 
$\fu=\op_{k\in\ZZ}\fu_k$ where $\fu_k=\fu\cap\fh_k$. Setting 
$\fl_k=\fp_k/\fu_k$, we have $\fl=\op_{k\in\ZZ}\fl_k$; this gives a 
$\ZZ$-grading of the Lie algebra $\fl$.

Now $\fp_0$ is a parabolic subalgebra of the reductive Lie algebra $\fh_0$; we
have $\fp_0=\fL P_0$ where $P_0$ is a parabolic subgroup of the connected 
reductive group $H_0$. Let $L_0$ be the image of $P_0$ under the obvious 
homomorphism $P@>>>L$. Then $L_0=e^{\fl_0}\sub L$. Now $P_0$ acts by the 
$\Ad$-action on each $\fp_k$. Let 
$\p:\fp_\et@>>>\fl_\et$ be the obvious projection. We have a diagram
$$\fl_\et@<a<<H_0\T\fp_\et@>b>>E@>c>>\fh_\et$$
where 
$$E=\{(hP_0,z)\in H_{0}/P_0\T\fh_\et;\Ad(h\i)z\in\fp_\et\},$$ 
$$a(h,z)=\p(\Ad(h\i)z), b(h,z)=(hP_0,z), c(gP_0,z)=z.$$
Now $a$ is smooth with connected fibres, $b$ is a principal $P_0$-bundle and 
$c$ is proper. If $A\in\cq(\fl_\et)$, then $a^*A$ is a $P_0$-equivariant 
semisimple complex on $H_0\T\fp_\et$ hence there is a well-defined semisimple
complex 
$A_1$ on $E$ 
such that $b^*A_1=a^*A$. Since $c$ is proper, the
complex
$$\ind_{\fp_\et}^{\fh_\et}(A):=c_!A_1$$
belongs to $\cq(\fh_\et)$. For $B\in\cd(\fh_\et)$ we can form
$$\res_{\fp_\et}^{\fh_\et}(B):=\p_!(B|_{\fp_\et})\in\cd(\fl_\et).$$
Thus we have functors
$\res_{\fp_\et}^{\fh_\et}:\cd(\fh_\et)@>>>\cd(\fl_\et)$,
$\ind_{\fp_\et}^{\fh_\et}:\cq(\fl_\et)@>>>\cq(\fh_\et)$.

When $\tfl$ is a Levi subalgebra of $\fp$ such that $\tfl=\op_{k\in\ZZ}\tfl_k$ 
with $\tfl_k=\tfl\cap\fh_k$, we will sometime consider 
$\ind_{\fp_\et}^{\fh_\et}(A)$ with $A\in\cq(\tfl_\et)$ by identifying 
$\tfl_\et=\fl_\et$ in an obvious way and $A$ with an object in $\cq(\fl_\et)$.  

\subhead 1.4\endsubhead
In the setup of 1.3 let $S'_P$ be the set of Levi subgroups of $P$ and let 
$S_P$ be the set of all $M\in S'_P$ such that, setting $\fL M=\fm$, 
$\fm_k=\fm\cap\fh_k$, we have $\fm=\op_{k\in\ZZ}\fm_k$, or equivalently such 
that $\Ad(\l(t))\fm=\fm$ for all $t\in\kk^*$. We have $S_P\ne\emp$; indeed, 
we can find $M\in S'_P$ such that $\l(k^*)\sub M$; then $M\in S_P$. Since $U$ 
acts simply transitively by conjugation on $S'_P$, it follows that 

(a) {\it The unipotent group $\{u\in U;u\l(t)=\l(t)u\qua\frl t\in \kk^*\}$ 
acts simply transitively by conjugation on $S_P$.}

\subhead 1.5. Blocks for $\cq(h_{\et})$\endsubhead
Let $\fM_\et(H)$ be the set of all systems 
$$(M,M_0,\fm,\fm_*,\tC)$$ 
where $M$ is 
a Levi subgroup of some parabolic subgroup of $H$, $\fm=\fL M$, $\fm_*$ is a 
$\ZZ$-grading of $\fm$ such that $\fm_k=\fm\cap\fh_k$ for all $k$, 
$M_0=e^{\fm_0}\sub M$ and $\tC$ is a cuspidal $M_0$-equivariant 
perverse sheaf on $\fm_\et$ (up to isomorphism). Note that $H_0$ acts by 
conjugation on $\fM_\et(H)$. Let $\un\fM_\et(H)$ be the set of orbits for this 
action. 

In the setup of 1.2 assume that $A$ is a simple $H_0$-equivariant perverse 
sheaf on $\fh_\et$. By \cite{\GRA, 7.5},

(a) {\it there exists $P,L,L_0,\fp,\fl$ as in 1.3 and a cuspidal 
$L_0$-equivariant perverse sheaf $C$ on $\fl_\et$ such that some shift of $A$ 
is a direct summand of $\ind_{\fp_\et}^{\fh_\et}(C)$.}
\nl
Assume now that $P',L',L'_0,\fp',\fl'$ is another quintuple like 
$P,L,L_0,\fp,\fl$ and that $C'$ is a cuspidal $L'_0$-equivariant 
perverse sheaf on $\fl'_\et$ such that some shift of $A$ is a direct summand of
$\ind_{\fp'_\et}^{\fh_\et}(C')$.

Let $M\in S_P,M'\in S_{P'}$, let $\fL M=\fm=\op_k\fm_k$ be as in 1.4 and let 
$\fL M'=\fm'=\op_k\fm'_k$ where $\fm'_k=\fm'\cap\fh_k$. Let 
$M_0=e^{\fm_0}\sub M$, $M'_0=e^{\fm'_0}\sub M'$. We can identify 
$M,M_0,\fm,\fm_k$ with $L,L_0,\fl,\fl_k$ via $P@>>>L$ and we can identify 
$M',M'_0,\fm',\fm'_k$ with $L',L'_0,\fl',\fl'_k$ via $P'@>>>L'$. Then $C$ 
(resp. $C'$) becomes a cuspidal $M_0$-equivariant (resp. 
$M'_0$-equivariant) perverse sheaf $\tC$ (resp. $\tC'$) on $\fm_\et$ (resp. 
$\fm'_\et$). 

Using the last sentence of \cite{\GRA, 15.3}, we see that there exists 
$h\in H_0$ such that $\Ad(h)$ carries $M,M_0,\fm,\fm_k$ to 
$M',M'_0,\fm',\fm'_k$ and $\tC$ to $\tC'$. Thus,

(b) {\it $A\m(M,M_0,\fm,\fm_k,\tC)$ is a well defined map from the set of 
(isomorphism classes) of simple $H_0$-equivariant perverse sheaves on $\fh_\et$
to the set $\un\fM_\et(H)$.}

\subhead 1.6\endsubhead
Let $(M,M_0,\fm,\fm_k,\tC)\in\fM_\et(H)$. We show:

(a) {\it there exists a parabolic subgroup $P$ of $H$ such that $M$ is a Levi
subgroup of $P$ and such that, setting $\fp=\fL P$, $\fp_k=\fp\cap\fh_k$, we 
have $\fp=\op_{k\in\ZZ}\fp_k$.}
\nl
Let $\cz=\cz^0_M$. Then $\fz=\fL\cz$ is the centre of $\fm$. Since 
$\fm_0$ is a Levi subalgebra of a parabolic subalgebra of $\fm$, we have 
$\fz\sub\fm_0$ hence $\cz\sub M_0$. We can find $\l_1\in Y_{\cz}$ such that the
centralizer of $\l_1(\kk^*)$ in $H$ is equal to the centralizer of $\cz$ in $H$
which equals $M$. Let $\l\in Y_H,r$ be as in 1.2. Then 
$\l(\kk^*)\sub\cz_{H_0}$. Now $\l_1(\kk^*)\sub\cz$ hence 
$\l_1(\kk^*)\sub H_0$. 
It follows that $\l_1(t)\l(t')=\l(t')\l_1(t)$ for any $t,t'$ in $\kk^*$. Thus 
we have $\fh=\op_{k\in\ZZ,k'\in\ZZ}({}^\l_{kr}\fh\cap{}^{\l_1}_{k'}\fh)$. 
Since the centralizer of $\l_1(\kk^*)$ in $\fh$ equals $\fm$, we have
$\fm=\op_{k\in\ZZ}({}^\l_{kr}\fh\cap{}^{\l_1}_0\fh)$. We set
$$\fp=\op_{k\in\ZZ,k'\in\ZZ_{\ge0}}({}^\l_{kr}\fh\cap{}^{\l_1}_{k'}\fh).$$
Clearly, $\fp$ is a parabolic subalgebra of $\fh$ with Levi subalgebra $\fm$ 
and such that, setting $\fp_k=\fp\cap\fh_k$, we have $\fp=\op_{k\in\ZZ}\fp_k$.
This proves (a).

\subhead 1.7\endsubhead
To any $(M,M_0,\fm,\fm_*,\tC)\in\fM_\et(H)$ we associate a simple perverse
sheaf $A$ in $\cq(\fh_\et)$ as follows. Let $\co$ be the $H_0$-orbit in 
$\fh_\et$ which contains $\ovsc\fm_\et$. Let $\cl'$ be the irreducible 
$M_0$-equivariant local system on $\ovsc\fm_\et$ such that 
$\tC|_{\ovsc\fm_\et}=\cl'[\dim\fm_\et]$. By \cite{\GRA, 11.2}, there is a well 
defined irreducible $H_0$-equivariant local system $\cl$ on $\co$ such that 
$\cl|_{\ovsc\fm_\et}=\cl'$. By definition, $A$ is the simple perverse sheaf on 
$\fh_\et$ such that $\supp A$ is contained in the closure of $\co$ and 
$A|_\co=\cl[\dim\co]$. 

\subhead 1.8\endsubhead
Assume that the $\ZZ$-grading $\fh_*$ of $\fh$ is $\et$-rigid. A  
perverse sheaf $A$ in $\cq(\fh_\et)$ is said to be {\it $\et$-semicuspidal} if 
$\supp A=\fh_\et$ and $A$ is attached to some 
$(M,M_0,\fm,\fm_*,\tC)\in\fM_\et(H)$ as in 1.7 (in particular, $A$ is a simple
perverse sheaf). In this case we have $\ovsc\fm_\et\sub\ovsc\fh_\et$; moreover, 

(a) {\it $H_0$ acts transitively on the set of systems 
$(M,M_0,\fp,\fp_*,\fm,\fm_*,\tC)$ such that 
$(M,M_0,\fm,\fm_*,\tC)\in\fM_\et(H)$, $A$ is attached to 
$(M,M_0,\fm,\fm_*,\tC)$ as in 1.7, $\fp$ is a parabolic subalgebra of $\fh$ 
with Levi subalgebra $\fm$ and $\fp=\op_{k\in\ZZ}\fp_k$ where 
$\fp_k=\fp\cap\fh_k$.}
\nl
(See \cite{\GRA, 11.9}.)

If $(M,M_0,\fp,\fp_*,\fm,\fm_*,\tC)$ is as in (a), then 
$$\ind_{\fp_\et}^{\fh_\et}(\tC)\cong\op_jA[-2s_j][\dim\fm_\et-\dim\fh_\et]
\tag b$$
where $s_j\in\NN$ are defined as follows. Choose $\ph=(e,h,f)\in J^H$ as in 
1.2(ii); let $H_\ph=\{g\in H;\Ad(g)(e)=e,\Ad(g)(h)=h,\Ad(g)(f)=f\}$, let 
$\cb$ be the variety of Borel subgroups of $H_\ph^0$; then $s_j$ are defined 
by $\r_{\cb!}\bbq=\op_j\bbq[-2s_j]$. (See \cite{\GRA, 11.13}.)

\subhead 1.9\endsubhead
Let $\cx$ be the set of all systems $(M,M_0,\fp,\fp_*,\fm,\fm_*,\tA)$ 
where $\fp$ is a parabolic subalgebra of $\fh$ with Levi subalgebra $\fm$,
$\fp=\op_{k\in\ZZ}\fp_k$, $\fm=\op_{k\in\ZZ}\fm_k$ where $\fp_k=\fp\cap\fh_k$,
$\fm_k=\fm\cap\fh_k$, $M=e^{\fm},M_0=e^{\fm_0}$ and $\tA$ is a simple perverse
sheaf in $\cq(\fm_\et)$ (up to isomorphism) which is $\et$-semicuspidal. We have
the following result, see \cite{\GRA, 13.3}.

(a) {\it Let $A_1\in\cq(\fh_\et)$. There exists
$C_1,C_2,\do,C_t,C_{t+1},\do,C_{t+t'}$ in $\cq(\fh_\et)$ such that 
$$A_1\op C_1\op C_2\op\do\op C_t=C_{t+1}\op\do\op C_{t+t'}$$ 
and each $C_j$ is of the form $\ind_{\fp_\et}^{\fh_\et}(\tA)[a_j]$
for some $(M,M_0,\fp,\fp_*,\fm,\fm_*,\tA)\in\cx$ (depending on $j$) and some 
$a_j\in\ZZ$.}
\nl

{\it Erratum to \cite{\GRA}.} In the definition of a good object in the second
paragraph of \cite{\GRA, 13.2}, one should insert the words "shifts of" after 
"direct sum of" (twice).

\subhead 1.10\endsubhead
Let $s\in\ZZ-\{0\}$. We show:

(a) {\it the subspace $\fh^{(1)}:=\op_{k\in s\ZZ}\fh_k$ of $\fh$ is the Lie algebra of
a well defined connected reductive subgroup $H^{(1)}$ of $H$.}
\nl
We can assume that $s>0$. We shall define $e\in Z_{\ge0}$ as follows: if $p=0$ we have $e=0$; if $p>0$ we define
$e$ by $s=s'p^e$ where $s'\in\ZZ_{>0}$ is not divisible by $p$. We shall argue by induction on $e$. (When $p=0$ we
only have to consider the case $e=0$.)
Assume first that $e=0$.

Let $\bar H$ be the adjoint group of $H$ and let $\bar{\fh}$ be its Lie algebra. Then $\bar{\fh}$ inherits a 
$\ZZ$-grading $\bar{\fh}=\op_k\bar{\fh}_k$ from $\fh$. If we assume known that 
$\bar{\fh}^{(1)}:=\op_{k\in s\ZZ}\bar{\fh}_k$ is the Lie algebra of
a well defined connected reductive subgroup $\bar H^{(1)}$ of $\bar H$ then we can take $H^{(1)}$ to be the
identity component of the inverse image of $\bar H^{(1)}$ under the obvious map $H@>>>\bar H$.
Thus we can assume that $H$ is adjoint. 
Let $\l\in Y_H$ be such that ${}^\l_k\fh=\fh_k$ for all $k$.
Let $\z'$ be a primitive $s$-th root of $1$ in $\kk$. (Note that if $p>0$, $s=s'$ is not divisible by $p$.)
We define $\o:H@>>>H$ by $\o(g)=\Ad(\l(\z'))(g)$; this is an automorphism of $H$.
The automorphism $\o':\fh@>>>\fh$ induced by $\o$ sends $x\in\fh_k$ (where $k\in\ZZ$) to $\z'{}^kx$. 
Hence $\o^s=1$ and $\fh^{(1)}$ is equal to $\{x\in\fh;\o(x)=x\}$.
Let $H^{(1)}$ be the identity component of $\{g\in H;\o(g)=g\}$. This is a connected reductive 
group with Lie algebra $\fh^{(1)}$. Thus (a) is proved in the case $e=0$.
We now assume that $e\ge1$ hence $p>0$.
We can find an element $x_0\in\fh$ such that $[x_0,x]=kx$ for any $k\in\ZZ$ and any $x\in\fg_k$. (We can take $x_0$ in the
image of the tangent map of $\l:\kk^*@>>>H$.)
Let $\ti\fh=\{x\in\fh;[x_0,x]=0\}$. We have $\ti\fh=\op_{k\in p\ZZ}\fh_k$.
Let $\ti H$ be the identity component of $\{g\in H;\Ad(g)x_0=x_0\}$. Since $x_0\in\fh$ is semisimple, it follows
that $\ti H$ is reductive with Lie algebra $\ti\fh$. We define a $\ZZ$-grading $\ti\fh=\op_{k'\in\ZZ}\ti\fh_{k'}$ by
$\ti\fh_{k'}=\fh_{pk'}$. By the induction hypothesis applied to $\ti H,\ti\fh$ we see that
there is a well defined connected reductive subgroup $\ti H^{(1)}$ of $\ti H$ whose Lie algebra is
$\op_{k'\in(s/p)\ZZ}\ti\fh_{k'}=\op_{k'\in(s/p)\ZZ}\fh_{pk'}=\op_{k\in s\ZZ}\fh_k=\fh^{(1)}$.
We can take $H^{(1)}=\ti H^{(1)}$. This completes the inductive proof.

\head 2. $\ZZ/m$-gradings and $\e$-spirals\endhead
In this section we introduce the key notion of this paper, namely a spiral. 
Spirals are analogues in the $\ZZ/m$-graded setting of parabolic subalgebras in the 
ungraded or $\ZZ$-graded setting. We also attach a canonical spiral to each nilpotent element in $\fg_{\d}$.

\subhead 2.1\endsubhead
In the rest of this paper, $m\ge1$, $G$, $\fg=\op_{i\in\ZZ/m}\fg_i$ are as in 
0.1 and $\z,\vt,\th$ are as in 0.5. Recall that for $i\in\ZZ/m$ we have
$\fg_i=\{x\in\fg;\th(x)=\z^ix\}$ and that $\vt:G@>>>G$ is the (semisimple) 
automorphism of $G$ which induces $\th:\fg@>>>\fg$; note 
that $\th(\Ad(g)x)=\Ad(\vt(g))\th(x)$ for all $x\in\fg$, $g\in G$.

We shall fix $\da\in\ZZ/m$.

For any semisimple automorphism $\g:G@>>>G$, we set $G^\g=\{g\in G;\g(g)=g\}$.
By a theorem of Steinberg \cite{\ST},

(a) {\it $G^\g$ is a connected reductive subgroup of $G$.}
\nl
Now $\fg_{\un0}$ is a Lie subalgebra of $\fg$. Recall that $G_{\un0}=G^\vt$
and that the $\Ad$-action of $G_{\un0}$ on $\fg$ leaves stable $\fg_i$ and its
closed subset $\fg_i^{nil}:=\fg_i\cap\fg^{nil}$ for any $i\in\ZZ/m$.

Let $\fG$ be the set of subgroups of $G$ of the form $G^{\Ad(\t)\vt}$ for some 
semisimple element of finite order $\t\in G_{\un0}$; by (a), any group in 
$\fG$ is a connected reductive subgroup of $G$. For example, we have 
$G_{\un0}\in\fG$; hence we have $G_{\un0}=e^{\fg_{\un0}}$.

\subhead 2.2\endsubhead
Let $\la,\ra:\fg\T\fg@>>>\kk$ be a Killing form;  
it is nondegenerate and it 
satisfies $\la\fg_i,\fg_j\ra=0$ whenever $i+j\ne\un0$ in $\ZZ/m$. Hence for 
any $i\in\ZZ/m$, $\la,\ra:\fg_i\T\fg_{-i}@>>>\kk$ is nondegenerate.

\subhead 2.3. The Morozov-Jacobson theorem in the $\ZZ/m$-graded setting\endsubhead
We set $J=J^G$, see 1.1. For $x\in\fg^{nil}$ let $J(x)=\{(e,h,f)\in J;e=x\}$,
$G(x)=\{g\in G;\Ad(g)x=x\}$ and let $U=U_{G(x)^0}$. Recall the following 
result of Morozov-Jacobson and Kostant, see \cite{\KO}.

(a) {\it We have $J(x)\ne\emp$. The $U$-action on $J(x)$ given by 
$$u:(e,h,f)\m u(e,h,f):=(e,\Ad(u)h,\Ad(u)f)$$
is simply transitive.}
\nl
Assume now that $x\in\fg^{nil}_\da$. We set
$$J_\da(x)=\{(e,h,f)\in J(x);e=x,h\in\fg_{\un0},f\in\fg_{-\da}\}.$$
We show:

(b) {\it We have $J_\da(x)\ne\emp$. The $(U\cap G_{\un0})$-action on $J_\da(x)$ 
(restriction of the $U$-action in (a)) is simply transitive.}
\nl
If $(e,h,f)\in J(x)$, then $(\z^{-\da}e,h,\z^\da f)\in J_\da(\z^{-\da}x)$ and 
$$(\z^{-\da}\th(e),\th(h),\z^\da\th(f))\in J(\z^{-\da}\th(x))=J(x)$$
(we use that $\th(e)=\z^\da e$). Hence 
$(e,h,f)\m(\z^{-\da}\th(e),\th(h),\z^\da\th(f))$ is a 
morphism $\th':J(x)@>>>J(x)$. Next we note that $g\m\vt(g)$ defines a 
homomorphism $G(x)@>>>G(x)$. (If $\Ad(g)x=x$ then 
$\th(x)=\th(\Ad(g)x)=\Ad(\vt(g))\th(x)$. 
Since $\th(x)=\z^\da x$, we see that $\z^\da x=\Ad(\vt(g))\z^\da x$ hence 
$x=\Ad(\vt(g))x$ and $\vt(g)\in G(x)$.) This restricts to a homomorphism 
$\th'':U@>>>U$ with fixed point set $U^{\th''}$. For $u\in U$, 
$(e,h,f)\in J(x)$ we have $\th'(u(e,h,f))=\th''(u)\th'(e,h,f)$. By (a), $J(x)$
is an affine space. Since $\th'{}^m=1$ and $m$ is invertible in $\kk$, the 
fixed point set $J(x)^{\th'}$ is nonempty. Since the $U$-action on $J(x)$ is 
simply transitive, it follows that this restricts 
to a simply transitive action of $U^{\th''}$ on $J(x)^{\th'}$. We have 
$J(x)^{\th'}=J_\da(x)$ and $U^{\th''}=U\cap G_{\un0}$. We see that (b) holds.

\subhead 2.4\endsubhead
Let $\l\in Y_{G_{\un0}}$ (resp. $\mu\in Y_{G_{\un0},\QQ}$). Since $\l$ (resp.
$\mu$) can be viewed as an element of $Y_G$ (resp. $Y_{G,\QQ}$), the 
decomposition $\fg=\op_{k\in\ZZ}({}_k^\l\fg)$ (resp.
$\fg=\op_{\k\in\QQ}({}_\k^\mu\fg))$ is defined as in 1.1. For $i\in\ZZ/m$ and 
for $k\in\ZZ$ (resp. $\k\in\QQ$) we set ${}_k^\l\fg_i={}_k^\l\fg\cap\fg_i$ 
(resp. ${}_\k^\mu\fg_i={}_\k^\mu\fg\cap\fg_i$); we then have
$\fg_i=\op_{k\in\ZZ}({}_k^\l\fg_i)$ (resp.
$\fg_i=\op_{\k\in\QQ}({}_\k^\mu\fg_i))$ for any $i\in\ZZ/m$ (we now use that 
$\l\in Y_{G_{\un0}}$ (resp. $\mu\in Y_{G_{\un0},\QQ}$)).

Let $s\in\ZZ-\{0\}$. We show:

(a) {\it the subspace $\fg^{(1)}:=\op_{k\in s\ZZ}({}^\l_k\fg_{\un{k/s}})$ of $\fg$ is the Lie algebra of
a well defined connected reductive subgroup $G^{(1)}$ of $G$.}
\nl
We apply 1.10(a) to $H=G,\fh=\fg$ with the $\ZZ$-grading $\fg=\op_k({}^\l_k\fg)$. We see that there
is a well defined reductive connected subgroup $H^{(1)}$ of $G$ whose Lie algebra is 
$\fh^{(1)}=\op_{k\in s\ZZ}({}^\l_k\fg)$.
Note that $H^{(1)}$ contains $\l(\kk^*)$ and is $\vt$-stable.
We choose $\z'\in\kk^*$ such that $\z'{}^s=\z$.
 We define $\o:H^{(1)}@>>>H^{(1)}$ by $\o(h)=\Ad(\l(\z'))\i\vt(h)$; this is an automorphism of $H^{(1)}$.
The automorphism $\o':\fh^{(1)}@>>>\fh^{(1)}$ induced by $\o$ sends $x\in{}^\l_k\fg_i$ (where $k\in s\ZZ$, $i\in\ZZ/m$)
to $\z'{}^{-k}\z^ix=\z^{i-\un{k/s}}x$. 
Hence $\o'{}^m=1$ and $\fg^{(1)}$ is equal to $\{x\in\fh^{(1)};\o'(x)=x\}$. Let 
$G^{(1)}$ be the identity component of $\{h\in H^{(1)};\o(h)=h\}$. Then $G^{(1)}$
is a connected reductive subgroup of $H^{(1)}$ with Lie algebra $\fg^{(1)}$. This proves (a).

Now ${}^\l_0\fg_{\un0}$ is a Levi subalgebra of a parabolic subalgebra of $\fg_{\un0}$. Hence
$e^{{}^\l_0\fg_{\un0}}$ is a well defined subgroup of $G_{\un0}$ (a Levi subgroup of a parabolic subgroup of $G_{\un0}$).
We have

(b) $e^{{}^\l_0\fg_{\un0}}\sub G^{(1)}$.

\subhead 2.5. The definition of $\e$-spirals\endsubhead
In the rest of this section we fix $\e\in\{1,-1\}$.  
For any $\mu\in Y_{G_{\un0},\QQ}$ and any $N\in\ZZ$ we set
$${}^\e\fp^\mu_N=\op_{\k\in\QQ;\k\ge N\e}({}_\k^\mu\fg_{\un N}).\tag a$$
If $r\in\ZZ_{>0}$ is such that $\l:=r\mu\in Y_{G_{\un0}}$ then we have
$${}^\e\fp^\mu_N=\op_{k\in\ZZ;k\ge rN\e}({}_k^\l\fg_{\un N}).$$ 
A collection $\{\fp_N;N\in\ZZ\}$ (or $\fp_*$) of subspaces of $\fg$ is said to
be an {\it $\e$-spiral} if there exists $\mu\in Y_{G_{\un0},\QQ}$ such that
$\fp_N={}^\e\fp^\mu_N$ for any $N\in\ZZ$. We then set (for $N\in\ZZ$)
$$\fu_N=\{x\in\fg_{\un N};\la x,{}^\e\fp^\mu_{-N}\ra=0\}=
\op_{\k\in\QQ;\k>N\e}({}_\k^\mu\fg_{\un N}).$$
We say that $\fu_*=\{\fu_N;N\in\ZZ\}$ is the nilradical of $\fp_*$.

The following properties of $\fp_*,\fu_*$ are immediate:  

$\do\sub\fp_N\sub\fp_{N-\e m}\sub\fp_{N-2\e m}\sub\do$ for any $N$;

$\fp_N\sub\fg_{\un N}$ for any $N$; $\fp_N=0$ if $N\e\gg0$; 
$\fp_N=\fg_{\un N}$ if $N\e\ll0$;

$[\fp_N,\fp_{N'}]\sub\fp_{N+N'}$ for any $N,N'$ in $\ZZ$;

$\do\sub\fu_N\sub\fu_{N-\e m}\sub\fu_{N-2\e m}\sub\do$ for any $N$;

$\fu_N\sub\fp_N$ for any $N$; $\fu_N=\fg_{\un N}$ if $N\e\ll0$;

$[\fu_N,\fp_{N'}]\sub\fu_{N+N'}$ for any $N,N'$ in $\ZZ$.
\nl
For $N\in\ZZ$ we set $\fl_N=\fp_N/\fu_N$ and $\fl=\op_{N\in\ZZ}\fl_N$. We have 
$\fl_N=0$ if $N\gg0$ or if $N\ll0$ hence $\dim\fl<\iy$; moreover, 
$[,]:\fp_N\T\fp_{N'}@>>>\fp_{N+N'}$ induces an operation 
$\fl_N\T\fl_{N'}@>>>\fl_{N+N'}$ which defines a Lie algebra structure on $\fl$.

Note that $\fp_0$ is a parabolic subagebra of the reductive Lie algebra 
$\fg_{\un0}$ and $\fu_0=\{x\in\fg_{\un0};\la x,\fp_0\ra=0\}$ is the 
nilradical of $\fp_0$. We set $P_0=e^{\fp_0}\sub G_{\un0}$, 
$U_0=e^{\fu_0}\sub G_{\un0}$. Then $P_0$ is a parabolic subgroup of $G_{\un0}$
and $U_0=U_{P_0}$, so that $L_0:=P_0/U_0$ is a connected reductive group. We 
note that:

(b) {\it the $\Ad$-action of $P_0$ on $\fg$ leaves stable $\fp_N$ and $\fu_N$ 
for any $N$.}
\nl
From (b) we see that for any $N$ there is an induced action of $P_0$ on 
$\fl_N=\fp_N/\fu_N$. We show:

(c) {\it the restriction of this action to $U_0$ is trivial.}
\nl
It is enough to show that the $\ad$-action of $\fu_0$ on $\fp_N/\fu_N$ is 
zero. This follows from the inclusion $[\fu_0,\fp_N]\sub\fu_N$ which has been 
noted earlier.

From (b),(c) we see that for any $N$ there is an induced action of 
$L_0=P_0/U_0$ on $\fl_N=\fp_N/\fu_N$. We show:

(d) {\it if $x\in\fp_N$, $N\e>0$, then $x\in\fg^{nil}_{\un N}$.}
\nl
It is enough to show that for any $x'\in\fg$ we have $\ad(x)^n(x')=0$ for 
$n\gg0$. We can assume that $x'\in\fg_i$ for some $i\in\ZZ/m$. If $N'\in\ZZ$ 
satisfies $\un N'=i$ and $N'\e\ll0$, then $\fp_{N'}=\fg_i$; thus we have 
$x'\in\fp_{N'}$ for some $N'$. We have $\ad(x)x'=[x,x']\in\fp_{N+N'}$, 
$\ad(x)^2(x')\in\fp_{2n+N'}$ and, more generally, $\ad(x)^n(x')\in\fp_{nN+N'}$
for $n\ge1$. If $n\gg0$ we have $nN\e+N'\e\gg0$ hence $\fp_{nN+N'}=0$; thus,
$\ad(x)^n(x')=0$. This proves (d).

\mpb

An element $\mu\in Y_{G_{\un0},\QQ}$ is said to be {\it $p$-regular} if $r\mu\in Y_{G_{\un0}}$ for some
$r\in\ZZ_{>0}$ such that $r\n p\ZZ$. (This condition holds automatically if $p=0$.)
An $\e$-spiral $\fp_*$ is said to be $p$-reqular if $\fp_*={}^\e\fp^\mu_*$ for some $p$-regular
$\mu\in Y_{G_{\un0},\QQ}$.

\subhead 2.6. Splittings of $\e$-spirals \endsubhead
For $\mu\in Y_{G_{\un0},\QQ}$ and $N\in\ZZ$ we set
$${}^\e\tfl^\mu_N=\op_{\k\in\QQ;\k=N\e}({}_\k^\mu\fg_{\un N})
={}_{N\e}^\mu\fg_{\un N}.$$
If $r\in\ZZ_{>0}$ is such that $\l:=r\mu\in Y_{G_{\un0}}$, then we have
$${}^\e\tfl^\mu_N={}_{rN\e}^\l\fg_{\un N}.$$ 

A {\it splitting} of an $\e$-spiral $\fp_*$ is a collection 
$\{\tfl_N;N\in\ZZ\}$ (or $\tfl_*$) of subspaces of $\fg$ such that for some 
$\mu\in Y_{G_{\un0},\QQ}$ we have $\fp_*={}^\e\fp^\mu_*$ and 
$\tfl_N={}^\e\tfl^\mu_N$ for any $N\in\ZZ$. 
Let $\fu_*$ be the nilradical of $\fp_*$. From the definitions we see that 
$\fp_N=\fu_N\op\tfl_N$ for any $N$, $[\tfl_N,\tfl_{N'}]\sub\tfl_{N+N'}$ for 
any $N,N'$ and the sum $\tfl:=\sum_{N\in\ZZ}\tfl_N$ (in $\fg$) is direct. Now 
$\tfl$ is a Lie subalgebra of $\fg$ which is $\ZZ$-graded by the subspaces 
$\tfl_N$. Note that the isomorphisms $\tfl_N@>\si>>\fl_N$ (restrictions of the 
obvious maps $\fp_N@>>>\fl_N$) give rise after taking $\op_N$ to an 
isomorphism $\tfl@>\si>>\fl$ which is compatible with the Lie algebra 
structures and the $\ZZ$-gradings.

For $\mu$ as above we can find $\l\in Y_{G_{\un0}}$ and $r\in\ZZ_{>0}$ such 
that $r\mu=\l$. 
Applying 2.4(a) with $s=r\e$ we see that

(a) {\it There is a well defined connected reductive subgroup $\tL$ of $G$ whose Lie
algebra is $\tfl$. In particular, $\tfl$ and $\fl$ are reductive Lie algebras.}
\nl

Let $\tL_0=e^{\tfl_0}$. From 2.4(b) we have

(b) $\tL_0\sub\tL$.

We show:

(c) {\it Assume that we have $\tfl_*={}^\e\tfl^\mu_*$, $\fp_*={}^\e\fp^\mu_*$ where 
$\mu$ is $p$-regular, that is, $\mu=r\l$ with $\l\in Y_{G_{\un0}}$ and $r\in\ZZ_{>0}$ such that $r\n p\ZZ$. 
Then there exists $\z'$, a root of $1$ in $\kk^*$ such that
$\tfl=\{x\in\fg;\Ad(\l(\z')\i)\th(x))=x\}$,
$\tL=G^{\Ad(\l(\z')\i)\vt}=e^{\tfl}\sub G$; note that $\tL\in\fG$.}
\nl
Let $\z'$ be a primitive root of $1$ of order $rm$ in $\kk^*$ 
such that $\z'{}^r=\z$. 
We have $\fg=\op_{k\in\ZZ,i\in\ZZ/m}({}^\l_k\fg_i)$, 
$\tfl_N={}_{Nr\e}^\l\fg_{\un N}$ for all $N\in\ZZ$.
For $k,N'\in\ZZ$ and $x\in{}_k^\l\fg_{\un N'}$ we have 
$$\Ad(\l(\z')\i)(\th(x))=\z'{}^{-k}\z^{N'}x=\z'{}^{rN'\e-k}x.$$ 
The condition that $\z'{}^{rN'\e-k}=1$ is that $rN'\e-k\in rm\ZZ$ or that 
$k\in r\ZZ$ and $\un{N'}=\un{k/(r\e)}$. We see that 
$$\{x\in\fg;\Ad(\l(\z')\i)(\th(x))=x\}
=\op_{k\in r\ZZ,i\in\ZZ/m;\un{k/(r\e)}=i}({}_k^\l\fg_i)
=\op_{N\in\ZZ}({}_{rN\e}^\l\fg_{\un N})=\tfl,$$
and (c) follows.

We return to the general case.

We have $\l(\kk^*)\sub\tL_0$; moreover 
$\Ad(\l(t))$ acts as identity on $\tfl_0={}_0^\l\fg_0=\fL\tL_0$; thus, 
$\l(\kk^*)\sub\cz_{\tL_0}$. Since $\kk^*$ is connected, we deduce

(d)) {\it $\l(\kk^*)\sub\cz^0_{\tL_0}$.}
\nl
Note that:

(e) {\it for $t\in\kk^*,N\in\ZZ$, $\Ad(\l(t))$ acts on $\fl_N$ as $t^{rN\e}$ 
times identity.}
\nl
We show:

(f) {\it If $\tfl_*$ is a splitting of an $\e$-spiral $\fp_*$ then $\tfl_*$ is
a splitting of an $(-\e)$-spiral.}
\nl
Let $\mu\in Y_{G_{\un0},\QQ}$ be such that $\tfl_*={}^\e\tfl^\mu_*$,
$\fp_*={}^\e\fp^\mu_*$. Let $\mu'=(-1)\mu\in Y_{G_{\un0},\QQ}$. Then 
$\tfl_*={}^{-\e}\tfl^{\mu'}_*$ is a splitting of the $(-\e)$-spiral 
${}^{-\e}\fp^{\mu'}_*$.

\subhead 2.7\endsubhead
Let $\fS$ be the set of splittings of an $\e$-spiral $\fp_*$. Clearly, 
$\fS\ne\emp$. Let $U_0$ be as in 2.5. Now $U_0$ acts on $\fS$ by
$u:\tfl_*\m\{\Ad(u)\tfl_N;N\in\ZZ\}$. (We use that $\Ad(u)\fp_N=\fp_N$ 
for any $N$.) We show:

(a) {\it this $U_0$-action on $\fS$ is simply transitive.}
\nl
Let $\fu_*$ be the nilradical of $\fp_*$. Let $\tfl_*\in\fS$, $\tfl'_*\in\fS$.
Since $\tfl_0,\tfl'_0$ are Levi subalgebras of $\fp_0$, there is a unique 
$u\in U_0$ such that $\Ad(u)\tfl_0=\tfl'_0$. It remains to show that this $u$ 
satisfies $\Ad(u)\tfl_N=\tfl'_N$ for any $N$. Let $\tfl=\op_N\tfl_N$, 
$\tfl'=\op_N\tfl'_N$ (these are Lie subalgebras of $\fg$) and let 
$\tL=e^{\tfl}\sub G$, $\tL'=e^{\tfl'}\sub G$. Let $\mu,\mu'$ in 
$Y_{G_{\un0},\QQ}$ be such that $\fp_*={}^\e\fp^\mu_*={}^\e\fp^{\mu'}_*$,
$\tfl_*={}^\e\tfl^\mu_*$, $\tfl'_*={}^\e\tfl^{\mu'}_*$. We can find 
$r\in\ZZ_{>0}$ such that $\l:=r\mu\in Y_{G_{\un0}}$, 
$\l':=r\mu'\in Y_{G_{\un0}}$.
Let $\tL_0$ be as in 2.6 and let $\tL'_0$ be the analogous subgroup of $\tL'$.
We now fix $N\in\ZZ$. The $\Ad$-action of $\tL_0$ (resp. $\tL'_0$) on $\fg$ 
leaves stable $\tfl_N,\fu_N$ (resp. $\tfl'_N,\fu_N$). Let $\tL''_0=u\tL_0u\i$, 
$\tfl''_N=\Ad(u)\tfl_N$; then the $\Ad$-action of $\tL''_0$ on $\fg$ leaves 
stable $\tfl''_N,\fu_N$. Since $\Ad(u)\tfl_0=\tfl'_0$, we have 
$u\tL_0u\i=\tL'_0$ hence $\tL'_0=\tL''_0$. Let $T$ be a maximal torus of 
$\tL'_0=\tL''_0$. Now the $\Ad$-action of $T$ on $\fg$ leaves stable 
$\tfl'_N,\tfl''_N,\fu_N,\fp_N$. Let $\cx=\Hom(T,\kk^*)$. For any $\a\in\cx$ 
let 
$$\fp_{N,\a}=\{x\in\fp_N;\Ad(\t)x=\a(\t)x\qua\frl\t\in T\},
\qua\fu_{N,\a}=\fu_N\cap\fp_{N,\a},$$
$$\tfl'_{N,\a}=\tfl'_N\cap\fp_{N,\a},\qua\tfl''_{N,\a}=\tfl''_N\cap\fp_{N,\a}.
$$
We have $\tfl'_N=\op_{\a\in\cx}\tfl'_{N,\a}$, 
$\tfl''_N=\op_{\a\in\cx}\tfl''_{N,\a}$, $\fu_N=\op_{\a\in\cx}\fu_{N,\a}$. Let 
$\car'=\{\a\in\cx;\tfl'_{N,\a}\ne0\}$, 
$\car''=\{\a\in\cx;\tfl''_{N,\a}\ne0\}$, 
$\ti\car=\{\a\in\cx;\fu_{N,\a}\ne0\}$. Since $\tfl'_N,\tfl''_N$ are $T$-stable 
complements of $\fu_N$ in $\fp_N$, the $T$-modules $\tfl'_N,\tfl''_N$ are 
isomorphic, hence $\car'=\car''$. Since $\l'(\kk^*)\sub\cz^0_{\tL'_0}$ 
(see 2.6(d)), we have $\l'(\kk^*)\sub T$; hence for any $\a\in\cx$ we can define 
$\a\bul\l'\in\ZZ$ by $\a(\l'(t))=t^{\a\bul\l'}$ for all $t\in\kk^*$. 

Assume that $\a\in\ti\car$. Then for any $t\in\kk^*$, $\Ad(\l'(t))$ acts on 
$\fu_{N,\a}$ as multiplication by $t^{\a\bul\l'}$ hence 
$\fu_{N,\a}\sub{}_{\a\bul\l'}^{\l'}\fg_{\un N}$; thus 
${}_{\a\bul\l'}^{\l'}\fg_{\un N}$ has a nonzero intersection with $\fu_N$, so 
that $\a\bul\l'>rN\e$. We see that $\ti\car\sub\{\a\in\cx;\a\bul\l'>rN\e\}$.
Assume now that $\a\in\car'$. Then for any $t\in\kk^*$, $\Ad(\l'(t))$ 
acts on $\fl'_{N,\a}$ as multiplication by $t^{\a\bul\l'}$ hence 
$\fl'_{N,\a}\sub{}_{\a\bul\l'}^{\l'}\fg_{\un N}$; thus, 
${}_{\a\bul\l'}^{\l'}\fg_{\un N}$ has a nonzero intersection with $\tfl'_N$, 
so that $\a\bul\l'=rN\e$. We see that 
$\car'\sub\{\a\in\cx;\a\bul\l'=rN\e\}$. It 
follows that $\car'\cap\ti\car=\emp$ so that $\fp_{N,\a}=\tfl'_{N,\a}$ for 
$\a\in\car'$.  
Since $\car'=\car''$, we have also $\car''\cap\ti\car=\emp$, so that 
$\fp_{N,\a}=\tfl''_{N,\a}$ for $\a\in\car''=\car'$. Thus, for 
$\a\in\car'=\car''$ we have $\tfl'_{N,\a}=\tfl''_{N,\a}$ hence 
$\tfl'_N=\tfl''_N$ and $\tfl'_N=\Ad(u)\tfl_N$. This proves (a). 

\mpb

For any splitting $\tfl_*$ of $\fp_*$ we denote by $\tL(\tfl_*)$ the connected 
reductive subgroup $\tL$ of $G$ associated to $\tfl_*$ in 2.6. The family of 
groups $(\tL(\tfl_*))$ indexed by the various splittings $\tfl_*$ of $\fp_*$ 
has the 
property that any two groups in the family are canonically isomorphic to each 
other; the isomorphism is provided by conjugation by a well defined $u\in U_0$
(this follows from (a)). It follows that the groups in the family can be 
identified with a single connected reductive group $L$ which is canonically 
isomorphic to each group in the family. Note that $L$ is canonically attached 
to the $\e$-spiral $\fp_*$ and that $\fL L=\fl$ canonically. Note also that 
$L_0$ in 2.5 is naturally a closed subgroup of $L$.

\subhead 2.8. Subspirals coming from parabolics of $\fl_{*}$ \endsubhead
Let $\fp_*$ be an $\e$-spiral. We define $\fu_*$, $\fl_*$, $\fl$ in terms of 
$\fp_*$ as in 2.5. Let $\fq$ be a parabolic subalgebra of $\fl$ compatible 
with the $\ZZ$-grading of $\fl$ that is, such that $\fq=\op_{N\in\ZZ}\fq_N$ 
where $\fq_N=\fq\cap\fl_N$. For any $N\in\ZZ$ let $\hfp_N$ be the inverse 
image of $\fq_N$ under the obvious map $\fp_N@>>>\fl_N$. We show:

(a) {\it $\hfp_*$ is an $\e$-spiral. Moreover, if $\fp_*$ is $p$-reqular then
$\hfp_*$ is $p$-regular.}
\nl
We can find $\mu\in Y_{G_{\un0},\QQ}$ such that $\fp_*={}^\e\fp^\mu_*$; let 
$\tfl_*={}^\e\tfl^\mu_*$. Let $\tL$ be as in 2.6.
Let $\tfq$ be the Lie subalgebra of $\tfl$ 
corresponding to $\fq$ under the obvious isomorphism $\tfl@>\si>>\fl$ and let 
$\tfq_N=\tfq\cap\tfl_N$ so that $\tfq=\op_{N\in\ZZ}\tfq_N$. We then have 
$\hfp_N=\fu_N\op\tfq_N$ for all $N$. Let $r\in\ZZ_{>0}$ be such that
$\l:=r\mu\in Y_{G_{\un0}}$; if $\fp_*$ is $p$-regular we assume in addition that
$r\n p\ZZ$.

From 2.6(e) we see that for $t\in\kk^*$, 
$\Ad(\l(t))$ leaves stable each $\tfq_N$ hence it leaves stable $\tfq$. It 
follows that $\kk^*$ acts via $t\m\Ad(\l(t))$ 
on the variety of Levi 
subalgebras of $\tfq$; since this variety is isomorphic to an affine space, 
there exists a Levi subalgebra $\fm$ of $\tfq$ such that 
$\Ad(\l(t))\fm=\fm$ for all $t\in\kk^*$. Let $R$ be the closed connected 
subgroup of $\tL$ (a torus) such that $\fL R$ is the centre of $\fm$. Since 
$\tfq$ is a parabolic subalgebra of $\tfl$ with Levi subalgebra $\fm$, we can 
find $\l'\in Y_R$ such that, setting for any $N'\in\ZZ$:
$${}_{N'}^{\l'}\tfl=\{x\in\tfl;\Ad(\l'(t))x=t^{N'}x\qua\frl t\in\kk^*\},$$ 
we have $\tfq=\op_{N'\in\ZZ_{\ge0}}({}_{N'}^{\l'}\tfl)$, $\fm={}_0^{\l'}\tfl$. 
We have $\fm=\op_N\fm_N$ where $\fm_N=\fm\cap\tfl_N$ and $\fm_0$ is a Levi 
subalgebra of a parabolic subalgebra of $\fm$. Hence a Cartan subalgebra of 
$\fm\cap\tfl_0$ is also a Cartan subalgebra of $\fm$, so that it contains the 
centre of $\fm$. Thus the centre of $\fm$ is contained in $\tfl_0$, so that 
$R\sub\tL_0$. Since for any $t,t'\in\kk^*$, $\l(t)$ is contained in 
$\cz_{\tL_0}$ and $\l'(t')\in\tL_0$, we have $\l(t)\l'(t')=\l'(t')\l(t)$. We 
can view $\l'$ as an element of $Y_{G_{\un0}}$ hence ${}_k^{\l'}\fg_i$ is 
defined for $k\in\ZZ,i\in\ZZ/m$ and we have
$\fg_i=\op_{k\in\ZZ}({}_k^{\l'}\fg_i)$ for any $i\in\ZZ/m$. We can find
$a\in\ZZ_{>0}$ such that ${}_k^{\l'}\fg_i=0$ for any $i\in\ZZ/m$ and any 
$k\in\ZZ-[-a,a]$. Let $b$ be an integer such that $b>2a$, $b\n p\ZZ$. We define 
$\l''\in Y_{G_{\un0}}$ by $\l''(t)=\l(t^b)\l'(t)=\l'(t)\l(t^b)$ for all 
$t\in\kk^*$. By definition, for $k\in\ZZ,i\in\ZZ/m$ we have
$$\align&{}_k^{\l''}\fg_i=\{x\in\fg_i;\Ad(\l(t^b)\l'(t))x=t^kx\qua
\frl t\in\kk^*\}\\&=\op_{k',k_2;k'\in b\ZZ,k_2\in\ZZ,k'+k_2=k}
({}_{k'/b}^\l\fg_i\cap{}_{k_2}^{\l'}\fg_i).\endalign$$ 
Note that ${}_k^{\l''}\fg_i=0$ unless $k=bk_1+k_2$ for some 
$k_1\in\ZZ\cap[-a,a]$, $k_2\in\ZZ$; in this case, $k_1,k_2$ are uniquely 
determined by $k$ since $b>2a$. Thus, we have
$${}_k^{\l''}\fg_i={}_{k_1}^\l\fg_i\cap{}_{k_2}^{\l'}\fg_i
\text{ if $k=bk_1+k_2$ with $k_1,k_2$ in }\ZZ,$$
$${}_k^{\l''}\fg_i=0,\text{ otherwise }.$$
Let $\mu'=\fra{1}{br}\l''\in Y_{G_{\un0},\QQ}$ and let 
$\fp'_*={}^\e\fp^{\mu'}_*$. For $N\in\ZZ$ we have
$$\fp'_N=\op_{k_1,k_2\in\ZZ;bk_1+k_2\ge Nbr\e,|k_2|\le a}
({}_{k_1}^\l\fg_{\un N}\cap{}_{k_2}^{\l'}\fg_{\un N}).$$
The only integer multiple of $b$ in $[-a,a]$ is $0$; hence the condition that 
$k_2\ge b(rN\e-k_1)$ (with $k_2\in[-a,a]$) is equivalent to the condition 
that either $0>b(rN\e-k_1)$, $k_2\in[-a,a]$ or that $0=b(rN\e-k_1)$, 
$k_2\in[0,a]$. Thus, $\fp'_N=X\op X'$ where 
$$X=\op_{k_1,k_2\in\ZZ;k_1>rN\e}
({}_{k_1}^\l\fg_{\un N}\cap{}_{k_2}^{\l'}\fg_{\un N})=
\op_{k_1\in\ZZ;k_1>rN\e}({}_{k_1}^\l\fg_{\un N})=\fu_N,$$
$$X'=\op_{k_1,k_2\in\ZZ;k_1=rN\e,k_2\ge0}
({}_{k_1}^\l\fg_{\un N}\cap{}_{k_2}^{\l'}\fg_{\un N})
=\tfl_N\cap(\op_{k_2\in\ZZ_{\ge0}}({}_{k_2}^{\l'}\fg_{\un N}))
=\tfl_N\cap\tfq=\tfq_N.$$
Thus, we have $\fp'_N=\fu_N\op\tfq_N=\hfp_N$. This proves (a).

\mpb

From the computation in the previous proof we can extract the following: 

(b) {\it the splitting ${}^\e\tfl^{\mu'}_*$ of the $\e$-spiral 
$\hfp_*={}^\e\fp^{\mu'}_*$ is equal to $\fm_*$}.

\subhead 2.9. The spiral attached to an element $x\in\fg^{nil}_\da$\endsubhead
{\it In the remainder of this paper we fix $\et\in\ZZ-\{0\}$ such that $\un\et=\da$.}
\nl
In this subsection we assume that $\e=\doet$, see 0.12. Let $x\in\fg_\da^{nil}$. 
We associate to $x$ an $\e$-spiral as follows. By 2.3(b), we can find 
$\ph=(e,h,f)\in J_\da(x)$ such that $e=x$. Let $\io=\io_\ph\in Y_G$ be as in 
1.1. Since the differential of $\io$ is the linear map $\kk@>>>\fg$, 
$z\m zh\in\fg_0$, we have $\io(\kk^*)\sub G_{\un0}$ so that $\io$ can be 
viewed as an element of $Y_{G_{\un0}}$. 
Then 

$\fp^\ph_*:={}^\e\fp^{(|\et|/2)\io}_*$ is an $\e$-spiral with splitting
$\tfl^\ph_*:={}^\e\tfl^{(|\et|/2)\io}_*$.
Note that for $N\in\ZZ$ we have
$$\fp^\ph_N=\op_{k\in\ZZ;k\ge2N\e}({}_{k/|\et|}^\io\fg_{\un N}),\qua
\tfl^\ph_N={}_{2N/\et}^\io\fg_{\un N}\text{ if }2N/\et\in\ZZ,\qua
\tfl^\ph_N=0\text{ if }2N/\et\n\ZZ.$$

We show that:

(a) {\it the $\e$-spiral $\fp^\ph_*$ is $p$-regular; it depends only on $x$, not on $\ph$.}
\nl
The $p$-regularity follows from the fact that $2\n p\ZZ$. We now prove the
second statement of (a).
By 2.3(b), another choice for $\ph$ must be of the form $u\ph$ where 
$u\in U_{G(x)^0}\cap G_{\un0}$. Let $\io'=\io_{u\ph}$. For $t\in\kk^*$ we have
$\io'(t)=u\io(t)u\i$ hence ${}_k^{\io'}\fg_i=\Ad(u)({}_k^\io\fg_i)$ for any
$k\in\ZZ,i\in\ZZ/m$. It follows that for $N\in\ZZ$ we have
$\fp^{u\ph}_N=\Ad(u)\fp^\ph_N$. To show that $\fp^{u\ph}_N=\fp^\ph_N$, it is 
enough to show that $\Ad(u)\fp^\ph_N=\fp^\ph_N$. It is enough to show:

$\Ad(u)({}_k^\io\fg)\sub\op_{k';k'\ge k}({}_{k'}^\io\fg)$ 
for any $u\in G(x),k\in\ZZ$.
\nl
Let $P$ be the parabolic subgroup of $G$ such that
$\fL P=\op_{k\in\ZZ;k\ge0}({}_k^\io\fg)$. Clearly,
$\Ad(g)({}_k^\io\fg)\sub\op_{k';k'\ge k}({}_{k'}^\io\fg)$ for any 
$g\in P,k\in\ZZ$. Hence it is enough to note the known inclusion 
$G(x)\sub P$. This proves (a).

\mpb

In view of (a) we will write $\fp_*^x$ instead of $\fp_*^\ph$, where $\ph$ is 
any element in $J_\da(x)$; let $\fu_*^x$ be the nilradical of $\fp_*^x$.
Now the splitting $\tfl^\ph_*$ depends in general on $\ph$. We set
$\tfl^\ph=\op_{N\in\ZZ}\tfl^\ph_N$; this is a $\ZZ$-graded Lie subalgebra of
$\fg$. Let $\tL^\ph=e^{\tfl^\ph}\sub G$; we have $\tL^\ph\in\fG$. Let 
$\tL^\ph_0=e^{\tfl^\ph_0}\sub\tL^\ph$. We show:

(b) {\it we have $x\in\tfl^\ph_\et$; more precisely, $x$ belongs to 
$\ovsc\tfl^\ph_\et$ (the open $\tL^\ph_0$-orbit on $\tfl^\ph_\et$).}
\nl
The first statement is the same as $x\in{}_2^\io\fg_\da$; this follows from the 
equality $[h,x]=2x$. The second statement can be deduced from 
\cite{\GRA, 4.2(a)}.

\mpb

We set $\tL^\ph_0(x)=\tL^\ph_0\cap G(x)$, $G_{\un0}(x)=G_{\un0}\cap G(x)$. We 
show:

(c) {\it the inclusion $\tL^\ph_0(x)@>>>G_{\un0}(x)$ induces an isomorphism on
the groups of components.}
\nl
Let $P_0$ be the parabolic subgroup of $G_{\un0}$ such that
$\fL P_0=\fp_0^x=\op_{k\in\ZZ;k\ge0}({}_k^\io\fg_{\un0})$ and let 
$U_0=U_{P_0}$. We 
set $P_0(x)=P_0\cap G(x)$, $U_0(x)=U_0\cap G(x)$. Then $\tL^\ph_0$ is a Levi 
subgroup of $P_0$ so that $P_0=\tL^\ph_0U_0$ (semidirect product) and 
$P_0(x)=\tL^\ph_0(x)U_0(x)$ (semidirect product). Since $U_0(x)$ is a connected
unipotent group we see that the inclusion $\tL^\ph_0(x)@>>>P_0(x)$ induces an 
isomorphism on the groups of components. It remains to show that 
$P_0(x)=G_{\un0}(x)$. As we have noted in the proof of (a), we have 
$G(x)\sub P$ hence $G_{\un0}(x)\sub P\cap G_{\un0}$; since $P\cap G_{\un0}$ 
and $P_0$ have the same Lie algebra, namely $\fp_0^x$, they must have the same
identity component; since $P_0$ is parabolic in $G_{\un0}$, we must have 
$P\cap G_{\un0}=P_0$, so that $G_{\un0}(x)\sub P_0$ and therefore 
$G_{\un0}(x)\sub P_0(x)$. Since the reverse inclusion is obvious, we see that 
$P_0(x)=G_{\un0}(x)$ and (c) is proved.

We show:

(d) {\it if $g\in G_{\un0}$ is such that $\Ad(g)\i(x)\in\fp^x_\et$, then 
$g\in P_0$.}
\nl
The assumption of (d) implies that $g\in P$. (We use \cite{\GRA, 5.7} applied 
to the trivial $\ZZ$-grading of $\fg$ that is, such that in \cite{\GRA, 3.1} 
we have $\fg_N=0$ for $N\ne0$.) Thus, we have $g\in P\cap G_{\un0}$. As in the
proof of (c) we have $P\cap G_{\un0}=P_0$ and (d) follows.

We show:

(e) {\it the $P_0$-orbit of $x$ in $\fp^x_\et$ is open dense in $\fp^x_\et$.}
\nl
We argue as in \cite{\GRA, 5.9}. It is enough to show that 
$\dim(P_0)-\dim(P_0\cap G(x))=\dim\fp_\et$ or equivalently that
$$\dim\fp^x_0-\dim\ker(\ad(x):\fp^x_0@>>>\fg_\da)=\dim\fp^x_\et.$$
Since $x\in\fp^x_\et$ (see (b)) and $[\fp^x_0,\fp^x_\et]\sub\fp^x_\et$, we have
$\ad(x)(\fp^x_0)\sub\fp^x_\et$ so that it is enough to show that
$$\dim\ker(\ad(x):\fp^x_0@>>>\fp^x_\et)=\dim\fp^x_0-\dim\fp^x_\et,$$
or equivalently, that $\ad(x):\fp^x_0@>>>\fp^x_\et$ is surjective. By the 
representation theory of $\fs\fl_2$, the linear map
$$\ad(x):\op_{k\in\ZZ;k\ge0}({}_k^\io\fg)@>>>\op_{k\in\ZZ;k\ge2}({}_k^\io\fg)$$
is surjective. This restricts for any $i\in\ZZ/m$ to a (necessarily 
surjective) map
$$\ad(x):\op_{k\in\ZZ;k\ge0}({}_k^\io\fg_i)@>>>\op_{k\in\ZZ;k\ge2}
({}_k^\io\fg_{i+\da}).$$
Taking $i=0$ we see that $\ad(x):\fp^x_0@>>>\fp^x_\et$ is surjective. This 
proves (e).

\mpb

The assignment $x\m\fp^x_*$ is a $\ZZ/m$-analogue of an assignment in the
case of $\ZZ$-graded Lie algebras given in \cite{\GRA, \S5} which is in turn 
modelled on a construction in \cite{\KL, 7.1}.  

\head 3. Admissible systems\endhead
In this section we introduce the set $\un\fT_\et$ of $G_{\un0}$-conjugacy 
classes of admissible systems, which will be used to index the blocks in 
$\cd_{G_{\un0}}(\fg_{\da}^{nil})$. 
We also define a map that assigns a pair $(\co,\cl)$ (where $\co$ is a $G_{\un0}$-orbit 
in $\fg_{\da}^{nil}$ and $\cl$ is an irreducible $G_{\un0}$-equivariant local system on it) an element in $\un\fT_\et$.

\subhead 3.1. Definition of admissible systems \endsubhead
We preserve the setup of 2.1. 

Let $\fT'_\et$ be the set consisting of all systems 
$(M,M_0,\fm,\fm_*,\tC)$ where $M\in\fG$, $\fm=\fL M$, $\fm_*$ is a 
$\ZZ$-grading of $\fm$, $M_0=e^{\fm_0}\sub M$, $\tC$ is a simple cuspidal 
$M_0$-equivariant perverse sheaf on $\fm_\et$ (up to isomorphism).

Until the end of 3.4 we fix $\dot\x=(M,M_0,\fm,\fm_*,\tC)\in\fT'_\et$. Let 
$\io\in Y_M$ be associated to $\tC$ as in 1.2(c),(a) (with $M,\tC$ instead of 
$H,A$), so that ${}^\io_k\fm=\fm_{\et k/2}$ for any $k\in\ZZ$ such that $\et k/2\in\ZZ$ and
${}^\io_k\fm=0$ for any $k\in\ZZ$ such that $\et k/2\n\ZZ$.
Then we have $\fm_{k'}={}^\io_{2k'/\et}\fm$ for $k'\in\ZZ$ such that $2k'/\et\in\ZZ$ and
$\fm_{k'}=0$ for $k'\in\ZZ$ such that $2k'/\et\n\ZZ$.
Note that $\io(\kk^*)$ is contained in $\cz_{M_0}^0$.

The system $\dot\x$ is said to be {\it admissible} if conditions (a),(b) below
are satisfied:

(a) we have $\fm_N\sub\fg_{\un N}$ for any $N\in\ZZ$;

(b) there exists an element $\t$ of finite order in the torus
$\io(\kk^*)\cz_M^0$ of $M_0$ such that $M=G^{\Ad(\t)\vt}$.
\nl
We now consider the following condition on $\dot\x$ which may or may not hold.

(c) {\it $\fm_*$ is a splitting of some $p$-regular $1$-spiral or, equivalently (see 
2.6(f)), of some $p$-regular $(-1)$-spiral.}
\nl
The following result will be proved in 3.2-3.4.

(d) {\it $\dot\x$ is admissible if and only if $\dot\x$ satisfies (c).}
\nl
We now make some comments on the significance of condition (b). Assume that 
condition (a) is satisfied and that $\t$ is any semisimple element of finite 
order of $G_{\un0}$ such that $M=G^{\Ad(\t)\vt}$. We show that we have 
automatically 
$$\t\in\io(\kk^*)\cz_M.\tag e$$
Note that $\vt(\t)=\t$ since $\t\in G_{\un0}$ hence $\t\in G^{\Ad(\t)\vt}=M$. 

Let $N\in\ZZ$ be such that $2N/\et\in\ZZ$.
Since $\fm_N\sub\fg_{\un N}$, $\th$ acts on $\fm_N$ as $\z^N$; since 
$\Ad(\t)\th$ acts as $1$ on $\fm$ we see that $\Ad(\t)$ acts on $\fm_N$ as 
$\z^{-N}$. On the other hand for $t\in\kk^*$, $\Ad(\io(t))$ acts on $\fm_N$ as
$t^{2N/\et}$. Hence if $t_0\in\kk^*$ satisfies $t_0^{2\et}=\z\i$, then we have 
$\Ad(\io(t_0))\Ad(\t\i)=t_0^{2N/\et}\z^N=\z^{-N}\z^N=1$ on $\fm_N$. It follows that
$\Ad(\io(t_0))\Ad(\t\i)=1$ on $\fm$. Since $\io(t_0)\t\i\in M$, we deduce that 
$\io(t_0)\t\i\in\cz_M$ hence $\t\in\io(\kk^*)\cz_M$, as asserted.

We see that condition (b) is a strengthening of (e) in which $\t$ is required 
to lie not only in $\io(\kk^*)\cz_M$ but in its identity component.

\subhead 3.2\endsubhead
We show:

(a) {\it for any element $\t_0$ of finite order in a torus $T$ there 
exists $\l_0\in Y_T$ such that $\t_0\in\l_0(\kk^*)$.}
\nl
We can find $c\in\ZZ_{>0}$ such that $c\n p\ZZ$ and $t_0^c=1$. Let
$\mu_c=\{z\in\kk^*;z^c=1\}$. For some $a\in\NN$ we can identify 
$T=(\kk^*)^a$ and $\t_0$ with $(z_1,\do,z_a)\in(\mu_c)^a\sub T$. Now $\mu_c$ 
is cyclic with generator $z_0$. Thus we have $z_1=z_0^{k_1},\do,z_a=z_0^{k_a}$
where $k_1,\do,k_a$ are integers. We define $\l_0\in Y_T$ by 
$t\m(t^{k_1},\do,t^{k_a})$. Then $\t_0=\l_0(z_0)$, as desired.

\mpb

We remark that the proof of (a) we can assume that 

(b) $k_1\in\ZZ_{>0},k_1\n p\ZZ$.
\nl
Indeed if $p=0$, then $k_1\n p\ZZ$ is automatic. Assume now that $p>0$. We write $k_1=k'_1p^e$ where 
$k'_1\in\ZZ-p\ZZ$, $e\in\ZZ_{\ge0}$. Define $z'_0\in\mu_c$ by $z'_0=z_0^{p^e}$. This is again a generator
of $\mu_c$.  (Recall that $c\n p\ZZ$.) We have $z_1=(z'_0)^{k'_1}$, $z_j=(z'_0)^{k'_j}$ where $k'_j\in\ZZ_{>0}$
for $j=2,3,\do,a$. Thus we can replace $z_0,k_1,\do,k_s$ by $z'_0,k'_1,\do,k'_s$ where $k'_1\in\ZZ_{>0},k'_1\n p\ZZ$.
This proves (b).

\mpb

We now assume that $\t$ as in 3.1(b) is given. We show:

(c) {\it there exist $f\in\ZZ_{>0}$ and $\l'\in Y_{\cz_M^0}$ such that $f\n p\ZZ$ and such that,
 if
$\l\in Y_{\io(\kk^*)\cz_M^0}$ is defined by $\l(t)=\io(t^f)\l'(t)$ for all 
$t$, then $\t\in\l(\kk^*)$.}
\nl
If $\io$ is identically $1$ then (c) follows from (a) applied to $T=\cz_M^0$ 
(we can take $f=1$). Assume now that $\io$ is not identically $1$. Then 
$\io:\kk^*@>>>M$ has finite kernel. Let $T=\kk^*\T\cz_M^0$; we define 
$d:T@>>>\io(k^*)\cz_M^0$ by $d(t,g)=\io(t)g$. By definition, $\io(k^*)$ is 
contained in the derived subgroup of $M$ hence it has finite intersection with 
$\cz_M^0$. It follows that $d$ has finite kernel. It is also surjective, hence
we can find $\ti\t\in T$ of finite order such that $d(\ti\t)=\t$. Using (a), we
can find $\l_0\in Y_T$  such that $\ti\t\in\l_0(\kk^*)$; moreover, by (b),
we can assume that, setting $\l_0(t)=(\l_1(t),\l'(t))$ with 
$\l_1\in Y_{\kk^*}$, $\l'\in Y_{\cz_M^0}$, we have $\l_1(t)=t^f$ for all $t$ 
where $f\in\ZZ_{>0}$, $f\n p\ZZ$. Let $\l=d\l_0:\kk^*@>>>\io(k^*)\cz_M^0$. We have
$\l(t)=\io(\l_1(t))\l'(t)=\io(t^f)\l'(t)$ for $t\in\kk^*$. Since $d(\ti\t)=\t$
and $\ti\t\in\l_0(\kk^*)$, we have $\t\in\l(\kk^*)$. This proves (c).

\subhead 3.3\endsubhead
We now assume that $\t$ as in 3.1(b) is given; let $\l,\l',f$ be as in 3.2(c).
We assume also that 3.1(a) holds. We can find $c\in\kk^*$ of finite order 
such that $\l(c)=\t$. (If $\t\ne1$ then $\l$ is not identically $1$ so it 
has finite kernel and any $c\in\l\i(\t)$ has finite order; if $\t=1$ we can
take $c=1$.)

Since $\l(\kk^*)\sub M_0$ and $M_0\sub G_{\un0}$ (as a consequence of our
assumption 3.1(a)), we can view $\l$ as an element of $Y_{G_{\un0}}$ hence
${}^\l_k\fg_i$ is defined for any $k\in\ZZ,i\in\ZZ/m$. Since 
$\l(\kk^*)\sub M$, we can view $\l$ as an element of $Y_M$ hence ${}^\l_k\fm$ 
is defined for any $k\in\ZZ$.

For $t\in\kk^*$, $k\in\ZZ$ such that $2k/\et\in\ZZ$ and $x\in\fm_k$ we have
$\Ad(\l(t))x=\Ad(\io(t^f))\Ad(\l'(t))x=\Ad(\io(t^f))x=t^{2kf/\et}x$ (we use 
that $\l'(t)\in\cz_M^0$). 
Thus $\fm_k\sub{}^\l_{2kf/\et}\fm$. Recall also that $\fm_k\ne0$ implies $k/\et\in\ZZ$, see 1.2(e).
Since the
subspaces $\fm_k$ form a direct sum decomposition of $\fm$ and the subspaces 
${}^\l_j\fm$ form a direct sum decomposition of $\fm$, it follows that
$$\align&
\fm_k={}^\l_{2kf/\et}\fm\text{ for any }k\in\et\ZZ \text{ and }\\&
{}^\l_j\fm=0\text{ unless $j=2kf/\et$ for some $k\in\et\ZZ$}.\tag a\endalign$$
For $k\in\ZZ,i\in\ZZ/m$ and $x\in{}^\l_k\fg_i$  we have
$$\Ad(\t)\th(x)=\Ad(\l(c))\th(x)=\z^i\Ad(\l(c))x=\z^ic^kx.$$
Since $\fm=\{x\in\fg;\Ad(\t)(\th(x))=x\}$, we see that 
$$\fm=\op_{j\in\ZZ,i\in\ZZ/m;\z^ic^j=1}({}^\l_j\fg_i).\tag b$$
If ${}^\l_j\fg_i$ is nonzero and contained in $\fm$ then ${}^\l_j\fm$ is 
nonzero hence by (a) we have $j=2fk/\et$ for some $k\in\ZZ$ and $\fm_k$ is 
a nonzero subspace of $\fg_i$; thus, by 3.1(a), we have $i=\un{k}$ and $2k/\et\in\ZZ$. 
Thus we can rewrite (b) as follows:
$$\fm=\op_{k\in\et\ZZ;\z^kc^{2fk/\et}=1}({}^\l_{2fk/\et}\fg_{\un{k}})$$
that is,
$$\fm=\op_{k\in\et\ZZ; (\z^\et c^{2f})^{k/\et}=1}({}^\l_{2fk/\et}\fg_{\un{k}}).\tag c$$
Assume now that $\fm_\et\ne0$. Using (a) we have $\fm_\et={}^\l_{2f}\fm\ne0$. By 
3.1(a) we have $\fm_\et\sub\fg_\da$. It follows that $\fm$ has nonzero 
intersection with ${}^\l_{2f}\fg_\da$. Now $\Ad(\t)\th$ acts on 
${}^\l_{2f}\fg_\da$ as multiplication by $\z^\et c^{2f}$ and it acts on $\fm$ 
as the identity. It follows that $\z^\et c^{2f}=1$. Thus (c) can be rewritten
as
$$\fm=\op_{k\in\et\ZZ}({}^\l_{2fk/\et}\fg_{\un{k}}).\tag d$$
Next we assume that $\fm_\et=0$. By the definition of $\io$ (see 3.1) this 
implies that $\io$ is identically $1$ hence $\fm=\fm_0$. From (a) we see that 
$\fm={}^\l_0\fm$, hence in (c) all summands corresponding to $k\ne0$ are 
zero. Thus (d) remains true in this case. We see also that
$$\fm_k={}^{|\et|\l}_{2fk\e}\fg_{\un k}$$
for all $k\in\ZZ$. Setting $\mu=|\et|\l/(2f)$ we see that
$\fm_*$ is a splitting of the $p$-regular $\e$-spiral 
${}^\e\fp^{\fra{1}{2f}|\et|\l}_*$. We see that if $\dot\x$ is admissible then it 
satisfies 3.1(c).

\subhead 3.4\endsubhead
Assume now that $\dot\x$ satisfies 3.1(c). 
Thus $\fm_*$ is a splitting of an $\e$-spiral $\fp_*={}^\e\fp^\mu_*$ where
$\mu$ is $p$-regular. Applying the conjugacy result 2.7(a) to the
two splittings $\fm_*,{}^\e\tfl^\mu_*$ we see that there exists
a $p$-regular $\mu'$ such that $\fp_*={}^\e\fp^{\mu'}_*$,
$\fm_*={}^\e\tfl^{\mu'}_*$. Thus we can find 
$\l\in Y_{G_{\un0}}$, $r\in\ZZ_{>0}$ such that $r\n p\ZZ$ and
$$\fm_N={}_{rN\e}^\l\fg_{\un N}$$
for any $N\in\ZZ$. In particular, 3.1(a) holds. We now show that 3.1(b) holds.
From 2.6(c) we see that $M=G^{\Ad(\l(\z')\i)\vt}$ for some root of unity 
$\z'\in\kk^*$. Let $\t=\l(\z')\i$. It remains to show that 
$\l(\z')\i\in\io(\kk^*)\cz_M^0$. More generally, we show that 
$\l(t)\in\io(\kk^*)\cz_M^0$ for any $t\in\kk^*$. Now $\l$ can be viewed as an 
element of $Y_M$ hence ${}^\l_k\fm$ is well defined for any $k\in\ZZ$ and we 
have ${}^\l_{rN\e}\fm=\fm_N$ for any $N\in\ZZ$. Recall that for $N\in\ZZ$ we have
$\fm_N={}^\io_{2N/\et}\fm$ if $N/\et\in\ZZ$ and $\fm_N=0$ if $N/\et\n\ZZ$.
We see that for any $N\in\et\ZZ$ and any $t\in\kk^*$, $\Ad(\l(t))$ 
acts on $\fm_N$ as $t^{rN\e}$ while $\Ad(\io(t^{|\et|}))$ acts
on $\fm_N$ as $t^{2N\e}$. Hence $\Ad(\l(t)^2\io(t)^{-r|\et|})$ acts on $\fm_N$ as 
$1$. Since $\fm$ is the sum of the subspaces $\fm_N$, we see that
$\Ad(\l(t)^2\io(t)^{-r|\et|})$ acts on $\fm$ as $1$. It follows that 
$\l(t)^2\io(t)^{-r|\et|}\in\cz_M$. Since $t\m\l(t)^2\io(t)^{-r|\et|}$ is a homomorphism
of the connected group $\kk^*$ into $\cz_M$, its image must be contained in 
$\cz_M^0$. Thus, for any $t\in\kk^*$ we have $\l(t)^2\io(t)^{-r|\et|}\in\cz^0_M$ 
hence $\l(t^2)\in\io(\kk^*)\cz^0_M$. Since any $t'\in\kk^*$ is a square, it
follows that $\l(t')\in\io(\kk^*)\cz^0_M$ for any $t'\in\kk^*$. We see that,
if $\dot\x$ satisfies 3.1(c) then $\dot\x$ is admissible. 
This completes the proof of 3.1(d).

\subhead 3.5. The map $\Ps:\ci(\fg_\da)@>>>\un\fT_\et$ \endsubhead
Let $\ci(\fg_\da)$ be the set of pairs $(\co,\cl)$ where $\co$ is a 
$G_{\un0}$-orbit on $\fg_\da^{nil}$ and $\cl$ is an irreducible 
$G_{\un0}$-equivariant local system on $\co$ defined up to isomorphism. Since 
$G_{\un0}$ acts on $\fg_\da^{nil}$ with finitely many orbits, see \cite{\VI},
the set $\ci(\fg_\da)$ is finite.

Let $\fT_\et$ be the set of all $(M,M_0,\fm,\fm_*,\tC)\in\fT'_\et$ which are 
admissible (see 3.1) or equivalently (see 3.1(d)) are such that $\fm_*$ is a 
splitting of some $p$-regular $\e$-spiral. The group $G_{\un0}$ acts in an obvious way by 
conjugation on $\fT_\et$; we denote by $\un\fT_\et$ the set of orbits, which is a finite 
set.
We will define a map $\Ps:\ci(\fg_\da)@>>>\un\fT_\et$. Let
$(\co,\cl)\in\ci(\fg_\da)$. Choose $x\in\co$ and $\ph\in J_\da(x)$; define 
$\fu^x_*$, $\tfl^\ph_*$, $\tfl^\ph,\tL^\ph,\tL^\ph_0$ as in 2.9. Recall that 
$\tL^\ph\in\fG$. We have $x\in\ovsc\tfl^\ph_\et$ (see 2.9(b)). By 2.9(c),
$\cl_1:=\cl|_{\ovsc\tfl^\ph_\et}$ is an irreducible $\tL^\ph_0$-equivariant 
local system on $\ovsc\tfl^\ph_\et$. Let $A$ be the simple 
$\tL^\ph_0$-equivariant perverse sheaf on $\tfl^\ph_\et$ whose restriction to 
$\ovsc\tfl^\ph_\et$ is $\cl_1[\dim\tfl^\ph_\et]$. The map 1.5(b) associates to 
$A$ an element $(M,M_0,\fm,\fm_*,\tC)$ of $\fM_\et(\tL^\ph)$ well defined up to
conjugation by $\tL^\ph_0$. Using 1.6(a) we can find a parabolic subalgebra 
$\fq$ of $\tfl^\ph$ compatible with the $\ZZ$-grading of $\tfl^\ph$ and such 
that $\fm$ is a Levi subalgebra of $\fq$. Setting $\fp'_N=\fu_N^\ph+\fq_N$ for
any $N\in\ZZ$, we see from 2.8(a) that $\fp'_*$ is a $p$-regular $\e$-spiral and from 
2.8(b) that $\fm_*$ is a splitting of $\fp'_*$. We see that 
$(M,M_0,\fm,\fm_*,\tC)\in\fT_\et$.

We now show that the $G_{\un0}$-orbit of $(M,M_0,\fm,\fm_*,\tC)$ is 
independent of the choices made. First, if $x,\ph$ are already chosen, then 
the $\tL^\ph_0$-orbit of \lb $(M,M_0,\fm,\fm_*,\tC)$ is well defined hence the 
$G_{\un0}$-orbit of $(M,M_0,\fm,\fm_*,\tC)$ is well defined (since 
$\tL^\ph_0\sub G_{\un0}$). The independence of the choice of $\ph$ (when $x$ is
given) follows from the homogeneity of $J_\da(x)$ under the group 
$U\cap G_{\un0}$ in 2.3(b). Finally, the independence of the choice of $x$ 
follows from the homogeneity of $\co$ under the group $G_{\un0}$. Thus, 
$$(\co,\cl)\m(G_{\un0}-\text{orbit of }(M,M_0,\fm,\fm_*,\tC))$$ 
is a well defined map $\Ps:\ci(\fg_\da)@>>>\un\fT_\et$.

\subhead 3.6\endsubhead 
Let $\dot\x=(M,M_0,\fm,\fm_*,\tC)\in\fT_\et$. Let $\co_{\dot\x}$ be the unique 
$G_{\un0}$-orbit in $\fg_{\da}^{nil}$ that contains $\ovsc\fm_\et$. Let 
$\dot\x'=(M',M'_0,\fm',\fm'_*,\tC')\in\fT_\et$. We show:

(a) {\it if $\co_{\dot\x}=\co_{\dot\x'}$, then there exists $g\in G_{\un0}$ 
such that $\Ad(g)$ carries $(M,M_0,\fm,\fm_*)$ to $(M',M'_0,\fm',\fm'_*)$.}
\nl
By \cite{\GRA, 3.3}, we can find $\ph=(e,h,f)\in J^M$, 
$\ph'=(e',h',f')\in J^{M'}$ such that 

(b) $e\in\ovsc\fm_\et$, $h\in\fm_0$, $f\in\fm_{-\et}$, $e'\in\ovsc\fm'_\et$, 
$h'\in\fm'_0$, $f'\in\fm'_{-\et}$.
\nl
We set $\io=\io_\ph\in Y_M$, $\io'=\io_{\ph'}\in Y_{M'}$. By 1.2(a),(c),(e), we have
$$\fm_k={}^\io_{2k/\et}\fm,\qua \fm'_k={}^{\io'}_{2k/\et}\fm'\text{ if }k\in\et\ZZ,
\fm_k=\fm'_k=0\text{ if }k\in\ZZ-\et\ZZ.\tag c$$
By assumption, we have $e'=\Ad(g_1)e$ for some $g_1\in G_{\un0}$. Replacing the system 
$(M,M_0,\fm,\fm_*,\tC,\ph)$ by its image under $\Ad(g_1)$, we see that we can 
assume that $e=e'$. Using 3.1(a) for $\dot\x$ and $\dot\x'$, we can view 
$\ph,\ph'$ as 
elements of $J^G_\da$ with the same first component. By 2.3(b), we can find 
$g_2\in G_{\un0}$ such that $\Ad(g_2)$ carries $\ph$ to $\ph'$. Replacing \lb
$(M,M_0,\fm,\fm_*,\tC,\ph)$ by its image under $\Ad(g_1)$, we see that we can 
assume that $\ph=\ph'$ as elements of $J^G$. It follows that $\io=\io'$ as 
elements of $Y_G$. 

Let 
$$G_\ph=\{g\in G;\Ad(g)(e)=e,\Ad(g)(h)=h,\Ad(g)(f)=f\}.$$
Since $e,h,f$ are contained in $\fm$ we have $\cz_M\sub G_\ph$. Similarly, 
since $e,h,f$ are contained in $\fm'$, we have $\cz_{M'}\sub G_\ph$. We have 
also $\cz_M^0\sub G_{\un0}$ (since the centre of $\fm$ is contained in
$\fm_0\sub\fg_{\un0}$); similarly we have $\cz_{M'}^0\sub G_{\un0}$. Thus,
$\cz_M^0$ and $\cz_{M'}^0$ are tori in $(G_\ph\cap G_{\un0})^0$. We show that 
$\cz_M^0$ is a maximal torus of $(G_\ph\cap G_{\un0})^0$. Indeed, assume that 
$S$ is a torus of $(G_\ph\cap G_{\un0})^0$ that contains $\cz_M^0$. Since 
$S\sub G_\ph$, for any $s\in S$ we have $\Ad(s)h=h$ hence $s\io(t)=\io(t)s$ 
that is, $\Ad(\io(t))s=s$ for $t\in\kk^*$. Since $S$ contains $\cz_M^0$, for 
any $s\in S,z\in\cz_M^0$ we have $\Ad(z)s=s$. Since $S\sub G_{\un0}$ we have 
$\vt(s)=s$ for any $s\in S$. We see that $\Ad(\io(t))\Ad(z)\vt(s)=s$ for any 
$t\in\kk^*$, $z\in\cz_M^0$, $s\in S$. We can find $\t\in\io(\kk^*)\cz_M^0$ such
that $M=G^{\Ad(\t)\vt}$. We have seen that $\Ad(\t)\vt(s)=s$ for $s\in S$. 
Thus $S\sub M$. Since $S\sub G_\ph$, we have 
$$S\sub M_\ph:=\{g\in M;\Ad(g)(e)=e,\Ad(g)(h)=h,\Ad(g)(f)=f\},$$
hence $S\sub M_\ph^0$. Since $e$ is a distinguished nilpotent element of 
$\fm$, we have $M_\ph^0=\cz_M^0$. Thus we have $S\sub\cz_M^0$. By assumption,
we have $\cz_M^0\sub S$, hence $\cz_M^0=S$. Thus $\cz_M^0$ is indeed a maximal
torus of $(G_\ph\cap G_{\un0})^0$, as claimed. Similarly we see that
$\cz_{M'}^0$ is a maximal torus of $(G_\ph\cap G_{\un0})^0$. Since any two
maximal tori of $(G_\ph\cap G_{\un0})^0$ are conjugate, we can find $g_3$ in
$(G_\ph\cap G_{\un0})^0$ such that $\Ad(g_3)$ carries $\cz_M^0$ to 
$\cz_{M'}^0$. (It also carries $\ph$ to $\ph$.) 

Replacing $(M,M_0,\fm,\fm_*,\tC,\ph)$ by its image under $\Ad(g_3)$, we see 
that we can assume that $\cz_M^0=\cz_{M'}^0$ and $\ph=\ph'$. 

Assume now that $e=0$ so that $e'=0$. By the definition of $\io=\io'$ we see that
$\io=\io'$ is identically $1$ hence $\fm=\fm_0$, $\fm'=\fm'_0$ and $G_\ph=G$.
Since $e=0$ is distinguished in $\fm$ it follows that $M$ is a torus.
Hence $M=\cz_M^0$. Similarly $M'=\cz_{M'}^0$. Since 
$\cz_M^0=\cz_{M'}^0$ it follows that $M=M'$. We see that (a) holds in this case.

In the remainder of the proof we assume that $e\ne0$ hence $e'\ne0$.
Recall that 
$M=G^{\Ad(\io(t))\Ad(z)\vt}$, $M'=G^{\Ad(\io(t'))\Ad(z')\vt}$,
for some $t,t'$ in $\kk^*$ and some $z,z'$ in $\cz_M^0=\cz_{M'}^0$. Since 
$e\in\fm_\et$, we have $\Ad(\io(t))\Ad(z)\th(e)=e$; since $\Ad(z)$ acts as $1$ 
on $\fm$, we deduce that $t^2\z^\et e=e$ and since $e\ne0$, we see that 
$t^2=\z^{-\et}$. Similarly, since $e\in\fm'_\et$ we have 
$\Ad(\io(t'))\Ad(z')\th(e)=e$ and $t'{}^2=\z^{-\et}$. 

We show that for any $k\in\ZZ$ we have $\fm_k\sub\fm'$. By 1.2(e) we can assume that 
$k\in\et\ZZ$.
Let $x\in\fm_k$. We must show that $\Ad(\io(t'))\Ad(z')\th(x)=x$. Since
$\Ad(z')$ acts by $1$ on $\fm$,  
it is enough to show that $\z^kt'{}^{2k/\et}x=x$
or that $(\z^\et t'{}^2)^{k/\et}x=x$. This follows from $t'{}^2=\z^{-\et}$.

Thus we have $\fm_k\sub\fm'$. Since this holds for any $k\in\ZZ$, we deduce 
that $\fm\sub\fm'$. Interchanging the roles of $\fm,\fm'$ we see that 
$\fm'\sub\fm$ hence $\fm=\fm'$. This implies that $M=M'$. Since $\io=\io'$, we
see from (c) that $\fm_*=\fm'_*$. From $\fm_0=\fm'_0$ we deduce that 
$M_0=M'_0$. This completes the proof of (a).

\mpb

The following result can be extracted from the proof of (a).

(d) {\it If $\fm_\et=0$ (so that $e=0$) then $\fm=\fm_0$ is a Cartan subalgebra of 
$\fg_{\un0}$.}

\subhead 3.7\endsubhead
Let $(M,M_0,\fm,\fm_*,\tC)\in\fT_\et$. Let $x\in\ovsc\fm_\et$. We choose 
$\ph=(e,h,f)\in J^M$ such that $e=x,h\in\fm_0,f\in\fm_{-\et}$ (see 
\cite{\GRA, 3.3}). We can view $x$ as an element of $\fg_\da^{nil}$ and $\ph$ 
as an element of $J_\da(x)$. We define $\tfl_*=\tfl^\ph_*$ as in 2.9.
Recall that for $N\in\ZZ$ we have 
$$\tfl_N={}_{2N/\et}^\io\fg_{\un N}\text{ if }2N/\et\in\ZZ,\qua
\tfl_N=0\text{ if }2N/\et\n\ZZ,$$
where $\io=\io_\ph\in Y_G$. Let $\tfl=\op_N\tfl_N\sub\fg$ and let 
$\tL=e^{\tfl}\sub G$. We show:

(a) {\it $\fm$ is a Levi subalgebra of a parabolic subalgebra of $\tfl$ which 
is compatible with the $\ZZ$-grading of $\tfl$.}
\nl
We shall prove (a) without the statement of compatibility with the 
$\ZZ$-grading; then the full statement of (a) would follow from 1.6(a).

Assume first that $x=0$. Then $h=0$ hence $\io$ is the constant map with
image $1$. It follows that $\tfl=\tfl_0=\fg_{\un0}$ and $\fm=\fm_0$; moreover,
by 3.6(d), $\fm$ is a Cartan subalgebra of $\fg_{\un0}$. Hence in
this case (a) is immediate. In the rest of the proof we assume that $x\ne0$.

Since $\ovsc\fm_\et$ carries a cuspidal local system, for any $N\in\ZZ$ such that
$2N/\et\in\ZZ$ we have $\fm_N={}^\io_{2N/\et}\fm$. Since $\fm_N\sub\fg_{\un N}$, we have
$\fm_N\sub{}^\io_{2N/\et}\fg_{\un N}$ hence $\fm_N\sub\tfl_N$. Taking sum over all 
$N\in\ZZ$ such that $2N/\et\in\ZZ$, we get $\fm\sub\tfl$. We can find $t_0\in\kk^*$, 
$z\in\cz_M^0$, both of 
finite order, such that $\fm=\{y\in\fg;\Ad(\io(t_0))\Ad(z)\th(y)=y\}$. 
Note that $\tfl_*={}^{\doet}\tfl^{(|\et|/2)\io}_*$

By 2.6(c),
we can find $\z'\in\kk^*$ such that 
$\tfl=\{y\in\fg;\Ad(\io(\z')\i)\th(y))=y\}$. Since $\fm\sub\tfl$, we have
$$\fm=\{y\in\tfl;\Ad(\io(t_0))\Ad(z)\th(y)=y\}=
\{y\in\tfl;\Ad(\io(t_0))\Ad(z)\Ad(\io(\z'))y=y\}.\tag b$$
(Note that 2.6(c) is applicable since $\tfl_*={}^{\doet}\tfl^{(|\et|/2)\io}_*$.)

Since $x\in\fm_\et\sub{}^\io_2\fg$, we have $\Ad(\io(t))x=t^2x$ for any $t$.
Taking $t=t_0\i$ or $t=\z'$ we see that $t_0^{-2}x=\Ad(\io(t_0))\i x$ and 
$\z'{}^2x=\Ad(\io(\z'))x$. Since $x\in\fm$ and $x\in\tfl$ we have
$\Ad(\io(t_0))\i x=\th(x)$ and $\Ad(\io(\z'))x=\th(x)$. It follows that
$t_0^{-2}x=\z'{}^2x$ so that (since $x\ne0$) we have $t_0^{-2}=\z'{}^2$.

If $N\in\ZZ,2N/\et\in\ZZ$ and $y\in\tfl_N$, we have $\Ad(\io(t_0\z'))y=(t_0\z')^{2N/\et}y=y$. Since 
$\tfl=\op_N\tfl_N$ we have $\Ad(\io(t_0\z'))y=y$ for all $y\in\tfl$. Hence (b)
implies

(c) $\fm=\{y\in\tfl;\Ad(z)y=y\}$.
\nl 
It remains to show that (c) implies (a).
Since $z$ is of finite order and $z\in\cz_M^0$, we can find 
$\l\in Y_{\cz_M^0}$ such that $z=\l(t_1)$ for some $t_1\in\kk^*$. (See 3.2(a).)

Let $\fm'=\{y\in\tfl;\Ad(\l(t))y=y\qua\frl t\in\kk^*\}$. Note that $\fm'$ 
is a Levi subalgebra of a parabolic subalgebra $\fq$ of $\tfl$. Since
$\l(\kk^*)\sub\cz_M^0$ we see that $\Ad(\l(t))$ acts as $1$ on $\fm$ for any 
$t$ hence $\fm\sub\fm'$. Now $\Ad(\l(t_1))$ acts as $1$ on $\fm'$. Since
$\fm=\{y\in\tfl;\Ad(\l(t_1))y=y\}$ it follows that $\fm'=\fm$ hence 
$\fm'=\fm$. Thus (a) holds.

\subhead 3.8. Primitive pairs \endsubhead
Let $(M,M_0,\fm,\fm_*,\tC)\in\fT_\et$. Let $x\in\ovsc\fm_\et$. We can view $x$ 
as an element of $\fg_\da^{nil}$. We set $M_0(x)=M_0\cap G(x)$, 
$G_{\un0}(x)=G_{\un0}\cap G(x)$. We show:

(a) {\it The inclusion $M_0(x)@>>>G_{\un0}(x)$ induces an isomorphism on the 
groups of components.} 
\nl
Let $\ph\in J^M$, $\tfl,\tfl_*,\tL$ be as in 3.7.  
Let $\tL_0=e^{\tfl_0}\sub\tL$. We have $x\in\ovsc\tfl_\et$ (see 
\cite{\GRA, 4.2(a)}). Let $\tL_0(x)=\tL_0\cap G(x)$. To prove (a) it is enough
to prove (i) and (ii) below.

(i) {\it The inclusion $M_0(x)@>>>\tL_0(x)$ induces an isomorphism on the 
groups of components.} 

(ii) {\it The inclusion $\tL_0(x)@>>>G_{\un0}(x)$ induces an isomorphism on 
the groups of components.} 
\nl
Now (i) follows from \cite{\GRA, 11.2} (we use 3.7(a)) and (ii) is a special case
of 2.9(c). This proves (a).

\mpb

Let $\co$ be the $G_{\un0}$-orbit of $x$ in $\fg_\da^{nil}$. Let $\cl'$ be the 
irreducible $M_0$-equivariant local system on $\ovsc\fm_\et$ such that 
$\tC|_{\ovsc\fm_\et}=\cl'[\dim\fm_\et]$. Let $\cl$ be the irreducible 
$G_{\un0}$-equivariant local system on $\co$ which corresponds to $\cl'$ under
(a). We say that $(\co,\cl)\in\ci(\fg_\da)$ is the {\it primitive pair} 
corresponding to $(M,M_0,\fm,\fm_*,\tC)\in\fT_\et$; it is clearly independent 
of the choice of $x,\ph$ (we use \cite{\GRA, 3.3}).

Let $\cl''$ be the irreducible $\tL_0$-equivariant local system on 
$\ovsc\tfl_\et$ which corresponds to $\cl'$ under (i). Let 
$\cl''{}^\sha\in\cd(\tfl_\et)$ be as in 0.11.
From 1.8(b) we see that

(b) {\it $\ind_{\fq_\et}^{\tfl_\et}(\tC)$ is a nonzero direct sum of shifts of 
$\cl''{}^\sha$.}
\nl
Consider the map $(M,M_0,\fm,\fm_*,\tC)\m(\co,\cl)$ (as above) from $\fT_\et$ to
$\ci(\fg_\da)$; the image of this map is denoted by $\ci^{prim}(\fg_\da)$. From 
3.6(a) and (a) we see that:

(c) {\it this induces a bijection $\o:\un\fT_\et@>\si>>\ci^{prim}(\fg_\da)$.}

\mpb

Using the definitions and 1.8(b), we see that: 

(d) {\it for $\x\in\un\fT_\et$ we have $\Ps(\o(\x))=\x$, where $\Ps:\ci(\fg_\da)@>>>\un\fT_\et$ is as in 3.5.}
\nl
Combining (c) and (d), we have

(e) {\it the restriction of  $\Ps$ to $\ci^{prim}(\fg_\da)$ gives the inverse of $\o$.}
\nl
From (d) we get

(f) {\it the map $\Ps:\ci(\fg_\da)@>>>\un\fT_\et$ is surjective.}
\nl
Another proof of (f) is given in 7.3.

\subhead 3.9\endsubhead  
Now let $\et_1\in\ZZ-\{0\}$ be such that $\un\et_1=\da$. We define a bijection
$\fT'_\et@>\si>>\fT'_{\et_1}$ by 
$(M,M_0,\fm,\fm_*,\tC)\m(M,M_0,\fm,\fm_{(*)},\tC)$
where $\fm_{(*)}$ is the new $\ZZ$-grading on $\fm_*$ whose $k$-component $\fm_{(k)}$
is equal to $\fm_{k\et/\et_1}$ for any $k\in\et_1\ZZ$ and is equal to $0$ for any
$k\in\ZZ-\et_1\ZZ$. (This is well defined since $\fm_{k'}=0$ for any $k'\in\ZZ-\et\ZZ$, see 
1.2(e).) Here we regard $\tC$ as a simple perverse sheaf on $\fm_\et=\fm_{(\et_1)}$.
This restricts to a bijection $\fT_\et@>\si>>\fT_{\et_1}$, which induces a bijection
$\un\fT_\et@>\si>>\un\fT_{\et_1}$. This allows us to identify canonically all the sets
$\fT_{\et_1}$ (for various $\et_1\in\ZZ-\{0\}$ such that $\un\et_1=\da$) with a single set
$\fT_\da$ and also 
all the sets $\un\fT_{\et_1}$ (for various $\et_1\in\ZZ-\{0\}$ such that $\un\et_1=\da$) with a
single set $\un\fT_\da$. Here $\fT_\da$, $\un\fT_\da$ are defined purely in terms of $\da$ 
(independently of the choice of $\et$).

\head 4. Spiral induction\endhead
In this section we introduce the key tool in studying the block decomposition for $\cd_{G_{\un0}}(\fg_{\da}^{nil})$, 
namely the spiral induction. This is the 
analogue in the $\ZZ/m$-graded setting of the parabolic induction in the ungraded or $\ZZ$-graded setting. 

\subhead 4.1. Definition of spiral induction \endsubhead
In addition to $\et\in\ZZ-\{0\}$ which has been fixed in 2.9, in this section we fix $\e\in\{1,-1\}$. 
We denote by $\fP^\e$ the set of all data of the form
$$(\fp_*,L,P_0,\fl,\fl_*,\fu_*)\tag a$$ 
where $\fp_*$ is an $\e$-spiral and $L,P_0,\fl,\fl_*,\fu_*$ are associated to 
$\fp_*$ as in 2.5. Let 
$$(\fp_*,L,P_0,\fl,\fl_*,\fu_*)\in\fP^\e.$$ 
Let $\p:\fp_\et@>>>\fl_\et$ be the obvious projection. We have a diagram 
$$\fl_\et@<a<<G_{\un0}\T\fp_\et@>b>>E'@>c>>\fg_\da\tag b$$
where $E'=\{(gP_0,z)\in G_{\un0}/P_0\T\fg_\da;\Ad(g\i)z\in\fp_\et\}$,
$a(g,z)=\p(\Ad(g\i)z)$, $b(g,z)=(gP_0,z)$, $c(gP_0,z)=z$. Here $a$ is smooth with 
connected fibres, $b$ is a principal $P_0$-bundle and $c$ is proper. Now 
$\cq(\fl_\et)$ is defined as in 1.2, with $H,\fh$ replaced by $L,\fl$. If 
$A\in\cq(\fl_\et)$, then $a^*A$ is a $P_0$-equivariant semisimple complex on 
$G_{\un0}\T\fp_\et$, hence there is a well-defined semisimple complex $A_1$ on 
$E'$ such that $b^*A_1=a^*A$. We can form the complex
$${}^\e\Ind_{\fp_\et}^{\fg_\da}(A)=c_!A_1.$$
Since $c$ is proper, this is a semisimple, $G_{\un0}$-equivariant complex on 
$\fg_\da$. 

If $\tfl_*$ is a splitting of $\fp_*$, we will sometimes 
consider
${}^\e\Ind_{\fp_\et}^{\fg_\da}(A)$ with $A\in\cq(\tfl_\et)$ by identifying
$\tfl_\et$ with $\fl_\et$ 
in an obvious way and $A$ with an object in $\cq(\fl_\et)$.

For any $A\in\cq(\fl_\et)$ we have
$$D({}^\e\Ind_{\fp_\et}^{\fg_\da}(A))
={}^\e\Ind_{\fp_\et}^{\fg_\da}(D(A))[2e]\tag c$$
where $e$ is the dimension of a fibre of $a$ minus the dimension of a fibre of
$b$, or equivalently
$$e=\dim\fg_{\un0}+\dim\fp_\et-\dim\fu_0-(\dim\fp_\et-\dim\fu_\et)
-(\dim\fp_0-\dim\fu_0)=\dim\fu_0+\dim\fu_\et.$$
Hence, if for $A\in\cq(\fl_\et)$ we set
$${}^\e\tInd_{\fp_\et}^{\fg_\da}(A)=
{}^\e\Ind_{\fp_\et}^{\fg_\da}(A)[\dim\fu_0+\dim\fu_\et],$$
then
$$D({}^\e\tInd_{\fp_\et}^{\fg_\da}(A))
={}^\e\tInd_{\fp_\et}^{\fg_\da}(D(A)).\tag d$$

\subhead 4.2. Transitivity \endsubhead
We state a transitivity property of induction. In addition to the datum 4.1(a)
we consider a parabolic subalgebra $\fq$ of $\fl$ such that 
$\fq=\op_{N\in\ZZ}\fq_N$ where $\fq_N=\fq\cap\fl_N$. For any $N\in\ZZ$ let 
$\hfp_N$ be the inverse image of $\fq_N$ under the obvious map 
$\fp_N@>>>\fl_N$. Then $\hfp_*$ is an $\e$-spiral, see 2.8(a). Let
$$(\hfp_*,\hL,\hP_0,\hfl,\hfl_*,\hfu_*)\in\fP^\e$$ 
be the datum analogous to 4.1(a) defined by $\hfp_*$. Now $\cq(\hfl_\et)$ is
defined as in 1.2, with $H,\fh$ replaced by $\hL,\hfl$. If $A\in\cq(\hfl_\et)$,
then $\ind_{\fq_\et}^{\fl_\et}(A)\in\cq(\fl_\et)$ is defined as in 1.3 
and we have canonically
$${}^\e\Ind_{\hfp_\et}^{\fg_\da}(A)
={}^\e\Ind_{\fp_\et}^{\fg_\da}(\ind_{\fq_\et}^{\fl_\et}(A)).
\tag a$$
The proof is similar to that of \cite{\CSI, 4.2}; it is omitted.

\subhead 4.3\endsubhead
In the setup of 4.1, assume that $A\in\cq(\fl_\et)$ is a cuspidal 
perverse sheaf (see 1.2). We have $A=\cl^\sha[\dim\fl_\et]$ where $\cl$ is an
irreducible local system on $\ovsc\fl_\et$ and $\cl^\sha\in\cd(\fl_\et)$ is as
in 0.11. In this case we can give an alternative 
description of ${}^\e\Ind_{\fp_\et}^{\fg_\da}(\cl^\sha)$. Let $P_0$, 
$\p:\fp_\et@>>>\fl_\et$ be as in 4.1. Let 
$$\dot\fg_\da=\{(gP_0,z)\in G_{\un0}/P_0\T\fg_\da;
\Ad(g\i)z\in\p\i(\ovsc\fl_\et))\},$$
an open smooth irreducible subvariety of $E'$ in 4.1. Let $\dot\cl$ be the 
local system on $\dot\fg_\da$ 
defined by $b'{}^*\dot\cl=a'{}^*\cl$ where
$$\ovsc\fl_\et@<a'<<G_{\un0}\T(\p\i(\ovsc\fl_\et))@>b'>>\dot\fg_\da,$$
$a'(g,zx)=\p(\Ad(g\i)z),b'(g,z)=(gP_0,z)$. Let $\dot\cl^\sha$ be the 
intersection cohomology complex of $E'$ with coefficients in $\dot\cl$. From 
the definitions we have $a^*\cl^\sha=b^*\dot\cl^\sha$ ($a,b$ as in 4.1). We 
define $c':\dot\fg_\da@>>>\fg_\da$ by $c'(g,z)=z$. We show:
$${}^\e\Ind_{\fp_\et}^{\fg_\da}(\cl^\sha)=c'_!\dot\cl.\tag a$$ 
Using the definitions we see that it is enough to show that the restriction of 
$\dot\cl^\sha$ to $E'-\dot\fg_\da$ is zero. This can be deduced from 1.2(c).

\subhead 4.4\endsubhead
Let $\cq^\e_\et(\fg_\da)$ be the subcategory of $\cd(\fg_\da)$ consisting of 
complexes which are direct sums of shifts of simple $G_{\un0}$-equivariant 
perverse sheaves $B$ on $\fg_\da$ with the following property: there exists a 
datum $(\fp_*,L,P_0,\fl,\fl_*,\fu_*)$ as in 4.1(a) and a simple cuspidal 
perverse sheaf $A$ in $\cq(\fl_\et)$ such that some shift of $B$ 
is a direct summand of ${}^\e\Ind_{\fp_\et}^{\fg_\da}(A)$. We show:

(a) {\it If $(\fp_*,L,P_0,\fl,\fl_*,\fu_*)\in\fP^\e$ and $A'$ is a simple 
(not necessarily cuspidal) perverse sheaf in $\cq(\fl_\et)$, then 
${}^\e\Ind_{\fp_\et}^{\fg_\da}(A')\in\cq^\e_\et(\fg_\da)$.}
\nl
Using \cite{\GRA, 7.5} we see that some shift of $A'$ is a direct summand of
$\ind_{\fq_\et}^{\fl_\et}(A)$ for some $\hfl,\fq$ as in 4.2 where $A$ is
a simple cuspidal perverse sheaf in $\cq(\hfl_\et)$. It follows that some shift
of ${}^\e\Ind_{\fp_\et}^{\fg_\da}(A')$ is a direct summand of 
$${}^\e\Ind_{\fp_\et}^{\fg_\da}(\ind_{\fq_\et}^{\fl_\et}(A)).\tag b
$$
It is then ehough to show that the complex (b) belongs to
$\cq^\e_\et(\fg_\da)$. This follows from the definitions using the 
transitivity property 4.2(a).

\mpb

The functor 
$${}^\e\Ind_{\fp_\et}^{\fg_\da}:\cq(\fl_\et)@>>>\cq^\e_\et(\fg_\da)$$
(where $(\fp_*,L,P_0,\fl,\fl_*,\fu_*)$ is as in 4.1(a)) is called
{\it spiral induction}. 

\mpb

Let $\ck^\e_\et(\fg_\da)$ be the abelian group 
generated by symbols $(A)$, one 
for each isomorphism class of objects of $\cq^\e_\et(\fg_\da)$, subject to 
the relations $(A)+(A')=(A\op A')$ (a Grothendieck group). Now 
$\ck^\e_\et(\fg_\da)$ is naturally an $\ca$-module by $v^n(A)=(A[n])$ for 
any $n\in\ZZ$. We shall write $A$ instead of $(A)$ (in $\cq^\e_\et(\fg_\da)$).
Clearly, $\ck^\e_\et(\fg_\da)$ is a free $\ca$-module with a finite
distinguished basis given by the various simple perverse sheaves in
$\cq^\e_\et(\fg_\da)$. Now $A,B\m(A:B)=\{A,D(B)\}\in\NN((v))$ (see 0.12) 
defines a pairing
$$(:):\ck^\e_\et(\fg_\da)\T\ck^\e_\et(\fg_\da)@>>>\ZZ((v)).\tag c$$
which is $\ca$-linear in the first argument, $\ca$-antilinear in 
the second argument (for $f\m\bar f$) and satisfies
$(b_1:b_2)=\ov{(b_2:b_1)}$ for all $b_1,b_2$ in $\ck^\e_\et(\fg_\da)$.

\subhead 4.5\endsubhead
In addition to the datum 4.1(a) we consider another datum
$$(\fp'_*,L',P'_0,\fl',\fl'_*,\fu'_*)\in\fP^\e\tag a$$ 
such that $\fp_N\sub\fp'_N$ for all $N\in\ZZ$ and $\fp_N=\fp'_N$ for 
$N\in\{\et,-\et\}$. We then have $\fu'_N\sub\fu_N$ for all $N\in\ZZ$ and 
$\fu_N=\fu'_N$ for $N\in\{\et,-\et\}$. We also have canonically $\fl_N=\fl'_N$ for
$N\in\{\et,-\et\}$ and $P_0\sub P'_0$. Let $\cp=P'_0/P_0$. Write $\r_{\cp!}\bbq
=\op_j\bbq[-2a_j]$ where 
$a_j$ are integers $\ge0$. Let $A\in\cq(\fl_\et)=\cq(\fl'_\et)$. We show:

(b) {\it Let $I={}^\e\Ind_{\fp_\et}^{\fg_\da}(A)$,
$I'={}^\e\Ind_{\fp'_\et}^{\fg_\da}(A)$. We have $I\cong\op_jI'[-2a_j]$.}
\nl
We consider the commutative diagram
$$\CD
\fl_\et@<a<<G_{\un0}\T\fp_\et@>b>>E'@>c>>\fg_\da\\
@V1VV        @V1VV                @VhVV     @V1VV  \\
\fl'_\et@<a'<<G_{\un0}\T\fp'_\et@>b'>>\tE'@>c'>>\fg_\da
\endCD$$
where the upper horizontal maps are as in 4.1(b), the lower horizontal are
the analogous maps when 4.1(a) is replaced by (a) and $h:E'@>>>\tE'$ is given 
by $(gP_0,z)\m(gP'_0,z)$. Note that $h$ is a $P'_0/P_0$-bundle. We can find a 
complex $A_1$ (resp. $A'_1$) on $E'$ (resp. $\tE'$) such that $I=c_!A_1$, 
$I'=c'_!A'_1$. We have $A_1=h^*A'_1$, hence
$$I=c_!A_1=c'_!h_!A_1=c'_!h_!h^*A'_1=c'_!(A'_1\ot h_!h^*\bbq)=
\op_jc'_!A'_1[-2a_j]=\op_jI'[-2a_j].$$
This proves (b).

\head 5. Study of a pair of spirals\endhead
This section serves as preparation for the next one, which aims to calculate the Hom space between two spiral inductions.

\subhead 5.1\endsubhead
In addition to $\et\in\ZZ-\{0\}$ which has been fixed in 2.9, in this section we fix $\e',\e''$ 
in $\{1,-1\}$. Let 
$$(\fp'_*,L',P'_0,\fl',\fl'_*,\fu'_*)\in\fP^{\e'},\qua
(\fp'',L'',P''_0,\fl'',\fl''_*,\fu''_*)\in\fP^{\e''}.$$
We show:

(a) {\it there exists a splitting $\tfl'_*$ of $\fp'_*$ and a splitting 
$\tfl''_*$ of $\fp''_*$ such that, if $\tL'_0=e^{\tfl'_0}\sub G$ and 
$\tL''_0=e^{\tfl''_0}\sub G$, then some maximal torus $\ct$ of $G_{\un0}$ is 
contained in both $\tL'_0$ and $\tL''_0$.} 
\nl
Let $\tfl'_*$ be any splitting of $\fp'_*$ and let $\tfl''_*$ be any splitting
of $\fp''_*$; let $\tL'_0=e^{\tfl'_0}\sub G$, $\tL''_0=e^{\tfl''_0}\sub G$. 
Recall that $P'_0$ (resp. $P''_0$) is a parabolic subgroup of $G_{\un0}$ with 
Levi subgroup $\tL'_0$ (resp. $\tL''_0$); hence there exists a maximal torus 
$\ct_0$ of $G_{\un0}$ contained in both $P'_0,P''_0$. Let ${}'\tL'_0$ (resp. 
${}'\tL''_0$) be the Levi subgroup of $P'_0$ (resp. $P''_0$) such that 
$\ct_0\sub{}'\tL'_0$ (resp. $\ct_0\sub{}'\tL''_0$). We can find 
$u'\in U_{P'_0}$, $u''\in U_{P''_0}$ such that $\Ad(u')\tL'_0={}'\tL'_0$, 
$\Ad(u'')\tL''_0={}'\tL''_0$. Now $\{\Ad(u')\tfl'_N;N\in\ZZ\}$ is a splitting 
of $\{\Ad(u')\fp'_N;N\in\ZZ\}=\fp'_*$ and $\{\Ad(u'')\tfl''_N;N\in\ZZ\}$ is a 
splitting of $\{\Ad(u'')\fp''_N;N\in\ZZ\}=\fp''_*$. Note that $\Ad(u')\tL'_0$,
$\Ad(u'')\tL''_0$ contain a maximal torus of $G_{\un0}$; (a) is proved.

\subhead 5.2\endsubhead
Let $(\fp'_*,L',P'_0,\fl',\fl'_*,\fu'_*)\in\fP^{\e'}$,
$(\fp_*,L,P_0,\fl,\fl_*,\fu_*)\in\fP^{\e''}$.
Let $A\in\cq(\fl_\et)$ be a simple cuspidal perverse sheaf. As in 4.3, we have 
$A=\cl^\sha[\dim\fl_\et]$ where $\cl$ is an irreducible local system on 
$\ovsc\fl_\et$. Let
$$B={}^{\e'}\Ind_{\fp_\et}^{\fg_\da}(\cl^\sha).$$
Let $\p':\fp'_\et@>>>\fl'_\et$ be the obvious map with kernel $\fu'_\et$. We want
to study the complex $K=\p'_!(B|_{\fp'_\et})\in\cd(\fl'_\et)$. As in 4.3, let
$$\dot\fg_\da=G_{\un0}\overset{P_{0}}\to{\times}\p\i(\ovsc\fl_\et),$$
where $\p:\fp_\et@>>>\fl_\et$ is the obvious map; let $\dot\cl$ be the local 
system on $\dot\fg_\da$ defined in terms of $\cl$ as in 4.3. As in 4.3,
we define $c':\dot\fg_\da@>>>\fg_\da$ by $c'(g,z)=z$. Let 
$$\dot\fp'_\et=\{(gP_0,z)\in G_{\un0}/P_0\T\fp'_\et;\Ad(g\i)z\in
\p\i(\ovsc\fl_\et)\}.$$ 
Note that $\dot\fp'_\et$ is the closed subvariety $c'{}\i\fp'_\et$ of 
$\dot\fg_{\da}$. 
The restriction of $\dot\cl$ from $\dot\fg_\da$ to 
$\dot\fp'_\et$ is denoted again by $\dot\cl$.
Now $c'$ restricts to a map $\dot\fp'_\et@>>>\fp'_\et$ whose composition with
$\p':\fp'_\et@>>>\fl'_\et$ is denoted by $\s:\dot\fp'_\et@>>>\fl'_\et$. 
We have $\s:(gP_0,z)\m\p'(z)$. Using 4.3(a) and proper base change, 
we see that $K=\s_!(\dot\cl)$.

We have a partition $\dot\fp'_\et=\cup_\Om\dot\fp'_{\et,\Om}$ into locally closed
subvarieties indexed by the various $(P'_0,P_0)$-double cosets $\Om$ in 
$G_{\un0}$ where 
$$\dot\fp'_{\et,\Om}=
\{(gP_0,z)\in\Om/P_0\T\fp'_\et;\Ad(g\i)z\in\p\i(\ovsc\fl_\et)\}.$$ 
Let $\s_\Om:\dot\fp'_{\et,\Om}@>>>\fl'_\et$ be the restriction of $\s$.
For any $\Om$ we set 
$$K_\Om=\s_{\Om!}(\dot\cl|_{\dot\fp'_{\et,\Om}})\in\cd(\fl'_\et).$$

\mpb

We say that $\Om$ is {\it good} if for some (or equivalently any) $g_0\in\Om$,
the following condition holds: setting $\fp''_N=\Ad(g_0)\fp_N$, 
$\fu''_N=\Ad(g_0)\fu_N$ (for $N\in\ZZ$), the obvious inclusion 
$$(\fp'_N\cap\Ad(g_0)\fp_N)/(\fp'_N\cap\Ad(g_0)\fu_N)@>>>
\Ad(g_0)\fp_N/\Ad(g_0)\fu_N$$
is an isomorphism for any $N\in\ZZ$ that is, 
$\Ad(g_0)\fp_N=(\fp'_N\cap\Ad(g_0)\fp_N)+\Ad(g_0)\fu_N$.
We say that $\Om$ is {\it bad} if it is not good.

\mpb

Until the end of 5.4 we fix an $\Om$ as above and we choose $g_0\in\Om$. Let 
$\fp''_N=\Ad(g_0)\fp_N$; then $\fp''_*$ is an $\e''$-spiral. It determines a 
datum $(\fp''_*,L'',P''_0,\fl'',\fl''_*,\fu''_*)\in\fP^{\e''}$.

By the change of variable $g=hg_0$ we may identify $\dot\fp'_{\et,\Om}$ with
$$\{(hP''_0,z)\in P'_0P''_0/P''_0\T\fp'_\et;
\Ad(h\i)z\in\Ad(g_0)\p\i(\ovsc\fl_\et))\}$$
which is the same as
$$\Xi=\{(h(P'_0\cap P''_0),z)\in P'_0/(P'_0\cap P''_0)\T\fp'_\et;
\Ad(h\i)z\in\p''{}\i(\ovsc\fl''_\et)\}$$
(in which $\p'':\fp''_\et@>>>\fl''_\et$ is the obvious map, with kernel 
$\fu''_\et$). In these coordinates, $\s_\Om:\dot\fp'_{\et,\Om}@>>>\fl'_\et$ 
becomes $(h(P'_0\cap P''_0),z)\m\p'(z)$.

We choose a splitting $\tfl'_*$ of $\fp'_*$ and a splitting $\tfl''_*$ of 
$\fp''_*$ as in 5.1(a); let $\tL'_0$, $\tL''_0$, $\ct$ be as in 5.1(a). 

Let $\mu',\mu''$ be elements of $Y_{G_{\un0},\QQ}$ such that
$\fp'_*={}^{\e'}\fp^{\mu'}_*$, $\tfl'_*={}^{\e'}\tfl^{\mu'}_*$,
$\fp''_*={}^{\e''}\fp^{\mu''}_*$, $\tfl''_*={}^{\e''}\tfl^{\mu''}_*$.
Let $r',r''$ in $\ZZ_{>0}$ be such that
$\l':=r'\mu'\in Y_{G_{\un0}}$, $\l'':=r''\mu''\in Y_{G_{\un0}}$.

As in 2.6(d) we have
$\l'(\kk^*)\sub\cz^0_{\tL'_0}$, $\l''(\kk^*)\sub\cz^0_{\tL''_0}$. Since $\ct$ 
is a maximal torus of $\tL'_0$, we must have $\cz^0_{\tL'_0}\sub\ct$ hence 
$\l'(\kk^*)\sub\ct$. Similarly, since $\ct$ is a maximal torus of $\tL''_0$, 
we have $\cz^0_{\tL''_0}\sub\ct$ hence $\l''(\kk^*)\sub\ct$. Since both
$\l'(\kk^*),\l''(\kk^*)$ are contained in the torus $\ct$, we must have
$\l'(t')\l''(t'')=\l''(t'')\l'(t')$ for any $t',t''$ in $\kk^*$. Hence, if for
any $k',k''$ in $\ZZ$ and $i\in\ZZ/m$ we set 
$${}_{k',k''}\fg_i=\{x\in\fg_i;\Ad(\l'(t'))x=t'{}^{k'}x,
\Ad(\l''(t''))x=t''{}^{k''}x, \quad \forall t',t''\in \kk^{*}\},$$
then we have $\fg=\op_{k',k'',i}({}_{k',k''}\fg_i)$.

\mpb

For any $N\in\ZZ$ we have a direct sum decomposition
$$\fp'_N\cap\fp''_N=(\tfl'_N\cap\tfl''_N)\op
(\fu'_N\cap\tfl''_N)\op(\tfl'_N\cap\fu''_N)\op(\fu'_N\cap\fu''_N).\tag a$$
This follows immediately from the decompositions
$$\fp'_N\cap\fp''_N=\op_{k',k'';k'\ge Nr'\e',k''\ge Nr''\e''}
({}_{k,k'}\fg_{\un N}),$$
$$\tfl'_N\cap\tfl''_N=\op_{k',k'';k'=Nr'\e',k''=Nr''\e''}
({}_{k,k'}\fg_{\un N}),$$
$$\fu'_N\cap\tfl''_N=\op_{k',k'';k'>Nr'\e',k''=Nr''\e''}
({}_{k,k'}\fg_{\un N}),$$
$$\tfl'_N\cap\fu''_N=\op_{k',k'';k'=Nr'\e',k''>Nr''\e''}
({}_{k,k'}\fg_{\un N}),$$
$$\fu'_N\cap\fu''_N=\op_{k',k'';k'>Nr'\e',k''>Nr''\e''}
({}_{k,k'}\fg_{\un N}).$$
For $N\in\ZZ$ let $\fq''_N$ be the image of $\fp'_N\cap\fp''_N$ under the 
obvious map $\fp''_N@>>>\fl''_N$; let $\fq''=\op_{N\in\ZZ}(\fq''_N)$, a Lie 
subalgebra of $\fl''$. We show:

(b) {\it $\fq''$ is a parabolic subalgebra of $\fl''$ compatible with the 
$\ZZ$-grading of $\fl''$.}
\nl
For $N\in\ZZ$ we set $\tfq''_N=\tfl''_N\cap\fp'_N$. Let 
$\tfq''=\op_{N\in\ZZ}\tfq''_N$, a Lie subalgebra of $\tfl''$. From (a) we see 
that the obvious isomorphism $\tfl''@>\si>>\fl''$ carries $\tfq''$ to $\fq''$.
Hence (b) follows from (c) below:

(c) {\it $\tfq''$ is a parabolic subalgebra of $\tfl''$ compatible with the 
$\ZZ$-grading of $\tfl''$.}
\nl
We have
$$\tfq''=\op_{k',N\in\ZZ;k'\ge Nr'\e'}({}_{k',Nr''\e''}\fg_{\un N}).$$
We define $\l_1\in Y_{\tL''}$ by $\l_1(t)=\l'(t^{r''})\l''(t^{-r'\e'\e''})$ 
for all $t\in\kk^*$. Then $\Ad(\l_1(t))$ acts on the subspace 
${}_{k',Nr''\e''}\fg_{\un N}$ of $\tfl''$ as $t^{k'r''-r'r''N\e'}$; the last 
exponent of $t$ is $\ge0$ if and only if $k'\ge r'N\e'$ which is just the 
condition that ${}_{k',Nr''\e''}\fg_{\un N}$ is one of the summands in the 
direct sum decomposition of $\tfq''$. This proves (c).

\mpb

For $N\in\ZZ$ let $\fq'_N$ be the image of $\fp'_N\cap\fp''_N$ under the 
obvious map $\fp'_N@>>>\fl'_N$; let $\fq'=\op_{N\in\ZZ}\fq'_N$, a Lie 
subalgebra of $\fl'$. 

For $N\in\ZZ$ we set $\tfq'_N=\tfl'_N\cap\fp''_N$. Let 
$\tfq'=\op_{N\in\ZZ}\tfq'_N$, a Lie subalgebra of $\tfl'$. The following 
result is proved in the same way as (b),(c).

(d) {\it $\fq'$ is a parabolic subalgebra of $\fl'$ compatible with the 
$\ZZ$-grading of $\fl'$; $\tfq'$ is a parabolic subalgebra of $\tfl'$ 
compatible with the $\ZZ$-grading of $\tfl'$.}

\mpb

We set ${}^!\tfq''=\op_N{}^!\tfq''_N$, ${}^!\tfq'=\op_N({}^!\tfq'_N)$, where
$${}^!\tfq''_N=\op_{k'\in\ZZ;k'>Nr'\e'}({}_{k',Nr''\e''}\fg_{\un N}),\qua
{}^!\tfq'_N=\op_{k''\in\ZZ;k''>Nr''\e''}({}_{Nr'\e',k''}\fg_{\un N}).$$ 
The proof of (c) shows also that ${}^!\tfq''$ is the nilradical of $\tfq''$ 
and that 
$$\op_{N\in\ZZ}({}_{Nr'\e',Nr''\e''}\fg_{\un N})$$ 
is a Levi subalgebra of $\tfq''$. Similarly, ${}^!\tfq'$ is the nilradical of 
$\tfq'$ and 
$$\op_{N\in\ZZ}({}_{Nr'\e',Nr''\e''}\fg_{\un N})$$ 
is a Levi subalgebra of $\tfq''$. Thus, 

(e) {\it $\tfq',\tfq''$ have a common Levi subalgebra, namely
$\op_{N\in\ZZ}({}_{Nr'\e',Nr''\e''}\fg_{\un N})$. }

\subhead 5.3\endsubhead
In this subsection we assume that $\Om$ is bad. Then for some $N$, 
$\tfl''_N\cap\fp'_N$ is strictly contained in $\tfl''_N$. Hence $\tfq''$ is a 
proper parabolic subalgebra of $\tfl''$ (see 5.2(c)). We will show that
$$K_\Om=\s_{\Om!}(\dot\cl|_{\dot\fp'_{\et,\Om}})=0\in\cd(\fl'_\et).\tag a$$
This is equivalent to the following statement:

(b) {\it for any $y\in\tfl'_\et$, the cohomology groups $H^j_c$ of the variety
$$\{(h(P'_0\cap P''_0),z)\in P'_0/(P'_0\cap P''_0)\T\fp'_\et; z-y\in\fu'_\et,
\Ad(h\i)z\in\p''{}\i(\ovsc\fl''_\et)\}$$
with coefficients in the local system defined by $\dot\cl$, are zero for all 
$j\in\ZZ$.}
\nl
(We have identified $\tfl'_\et$, $\fl'_\et$ via $\p'$.)
Considering the fibres of the first projection of the last variety to
$P'_0/(P'_0\cap P''_0)$, we see that it suffices to show that

(c) {\it for any $h\in P'_0$ and any $y\in\tfl'_\et$, the cohomology groups 
$H^j_c$ of the variety
$$\{z\in\fp'_\et;z-y\in\fu'_\et,\Ad(h\i)z\in\ovsc\tfl''_\et+\fu''_\et\}$$
with coefficients in the local system defined by $\dot\cl$, are zero for all 
$j\in\ZZ$.}
\nl
(We have used that $\p''{}\i(\ovsc\fl''_\et)=\ovsc\tfl''_\et+\fu''_\et$.)

If $z$ is as in (c), then we have automatically $\Ad(h\i)z\in\fp'_\et$; since
$\ovsc\tfl''_\et+\fu''_\et\sub\fp''_\et$, the condition that
$\Ad(h\i)z\in\ovsc\tfl''_\et+\fu''_\et$ implies 
$\Ad(h\i)z\in\fp'_\et\cap\fp''_\et$. By 5.2(a), we can then write uniquely 
$\Ad(h\i)z=\g+\nu'+\nu''+\mu$ where
$$\g\in\tfl'_\et\cap\tfl''_\et,\nu'\in\fu'_\et\cap\tfl''_\et,
\nu''\in\tfl'_\et\cap\fu''_\et, \mu\in\fu'_\et\cap\fu''_\et.\tag e$$
The condition that $\Ad(h\i)z\in\ovsc\tfl''_\et+\fu''_\et$ can be expressed as
$\g+\nu'\in\ovsc\tfl''_\et$. The condition that $z-y\in\fu'_\et$ is equivalent 
to $\Ad(h\i)z-\Ad(h\i)y\in\fu'_\et$ or (if we define $y'\in\tfl'_\et$ by
$\Ad(h\i)y-y'\in\fu'_\et$) to $\g+\nu''=y'$. Note that $y',\g,\nu''$ are 
uniquely determined by $h,y$. Hence the variety in (c) can be identified with
$$(\g+(\fu'_\et\cap\tfl''_\et))\cap\ovsc\tfl''_\et)\T(\fu'_\et\cap\fu''_\et).$$
Under this identification, the local system $\dot\cl$ is 
the pullback of $\cl$ (viewed as a local system on $\ovsc\tfl''_\et$) from the first factor. 
Now the desired vanishing of cohomology follows from the vanishing property
\cite{\GRA, 4.4(c)} of $\cl$, since in our case 
$\tfq''=\op_N(\tfl''_N\cap\fp'_N)$ 
is a proper parabolic subalgebra of $\tfl''$ with nilradical
$\op_N(\tfl''_N\cap\fu'_N)$.

\subhead 5.4\endsubhead
In this subsection we assume that $\Om$ is good. Then for any $N$ we have
$\tfl''_N\cap\fp'_N=\tfl''_N$ that is, $\tfl''_N\sub\fp'_N$. We also have
$\tfq''=\tfl''$. Thus $\tfq''$ is reductive so it is equal to its Levi 
subalgebra $\op_{N\in\ZZ}({}_{Nr'\e',Nr''\e''}\fg_{\un N})$ (see 5.2(e)) which
is then equal to $\tfl''$ and is also a Levi subalgebra of $\tfq'$ (see 
5.2(e)). Thus,

(a) {\it $\tfl''$ is a Levi subalgebra of $\tfq'$.}
\nl
Now $\Ad(g_0)$ defines an isomorphism $\fl@>\si>>\fl''$. Composing this with 
the inverse of the obvious isomorphism $\tfl''@>\si>>\fl''$ we obtain an 
isomorphism of $\ZZ$-graded Lie algebras $\fl@>\si>>\tfl''$. Using this, we 
can transport $\cl$ (a local system on $\ovsc\fl_\et$, see 5.1) to a local 
system $\cl''$ on $\ovsc\tfl''_\et$. Let $\cl''{}^\sha\in\cd(\tfl''_\et)$
be as in 0.11. Then 
$$\ind_{\tfq'_\et}^{\tfl'_\et}(\cl''{}^\sha)\in\cq(\tfl'_\et)$$
is defined as in 1.3 (we identify $\tfl''$ with the reductive quotient of
$\tfq'$, see (a)). We now state the following result.

(b) {\it We have 
$K_\Om=\ind_{\tfq'_\et}^{\tfl'_\et}(\cl''{}^\sha)[-2f]$ where}
$$f=\dim\fu'_0-\dim(\fu'_0\cap\fp''_0)+\dim(\fu'_\et\cap\fu''_\et).$$
Let $\tQ'_0=e^{\tfq'_0}\sub\tL'_0$, a parabolic subgroup of $\tL'_0$. Let
$$\Xi'=\{({}'h\tQ'_0,{}'z)\in\tL'_0/\tQ'_0\T\tfl'_\et;
\Ad({}'h\i){}'z\in\ovsc\tfl''_\et+{}^!\tfq_\et\}.$$
Define $c'':\Xi'@>>>\tfl'_\et$ by $c''({}'hQ'_0,{}'z)={}'z$.
By the argument in \cite{\GRA, 6.6} (for $\tL'$ instead of $G$) we have
$$\ind_{\tfq'_\et}^{\tfl'_\et}(\cl''{}^\sha)=c''_!\dot\cl''\tag c$$ 
where
$\dot\cl''$ is a certain local system on $\Xi'$ determined by $\cl''$ and such 
that $\D^*\dot\cl''=\dot\cl$ where $\D:\Xi@>>>\Xi'$ ($\Xi$ as in 5.2) is the 
map induced by the canonical maps $P'_0@>>>\tL'_0$ (with kernel $U_{P'_0}$) 
and $\fp'_\et@>>>\tfl'_\et$ (with kernel $\fu'_\et$); $\dot\cl$ is the local 
system on $\Xi$ considered in 5.2. We consider the following statement:

(d) {\it $\D$ is an affine space bundle with fibres of dimension $f$.}
\nl
Assuming that (d) holds, we have 
$$K_\Om=c''_!\D_!\dot\cl=c''_!\dot\cl''\ot\D_!\bbq=c''_!\dot\cl''[-2f]$$ 
and we see that (b) follows from (c). It remains to prove (d).

Let ${}'h\in\tL'_0$, ${}'z\in\tfl'_\et$ be such that $({}'hQ'_0,{}'z)\in\Xi'$.
Setting $h={}'hu$, $z={}'z+\tz$, we see that $\D\i({}'hQ'_0,{}'z)$ can be
identified with
$$\{(u(U_{P'_0}\cap P''_0),\tz)\in 
(U_{P'_0}/(U_{P'_0}\cap P''_0))\T\fu'_\et;\Ad(u\i)\Ad({}'h\i)({}'z+\tz)\in
\ovsc\tfl''_\et+\fu''_\et\}.$$
It suffices to show that
$$\{(u,\tz)\in U_{P'_0}\T\fu'_\et;\Ad(u\i)\Ad({}'h\i)({}'z+\tz)\in
\ovsc\tfl''_\et+\fu''_\et\}\tag e$$
is isomorphic to $U_{P'_0}\T(\fu'_\et\cap\fu''_\et)$. If $(u,\tz)$ are as in 
(e), we 
have automatically $\Ad(u\i)\Ad({}'h\i)({}'z+\tz)\in\fp'_\et$ (since 
${}'z+\tz\in\fp'_\et$ and ${}'hu\in P'_0$). Setting
$\Ad({}'h\i){}'z=a\in\ovsc\tfl''_\et+{}^!\tfq_\et$ (where $a$ is fixed) and
$\Ad(u\i)\Ad({}'h\i)\tz=\tz'\in\fu'_\et$, we see that the variety (e) may be 
identified with the variety
$$\{(u,\tz')\in U_{P'_0}\T\fu'_\et;\Ad(u\i)a+\tz'\in
\ovsc\tfl''_\et+(\fp'_\et\cap\fu''_\et)\}.\tag f$$
By 5.2(a) we can write uniquely
$$\Ad(u\i)a+\tz'=\g+\nu+\mu$$
where $\g\in\ovsc\tfl''_\et$, $\nu\in\tfl'_\et\cap\fu''_\et$, 
$\mu\in\fu'_\et\cap\fu''_\et$. 
Setting $\hz=\mu-\tz$ we see that (f) can be identified with the variety of all
quintuples $(u,\hz,\g,\nu,\nu')$ in 
$$U_{P'_0}\T\fu'_\et\T\ovsc\tfl''_\et\T(\tfl'_\et\cap\fu''_\et)\T(\fu'_\et\cap\fu''_\et)
$$  
such that 
$$\Ad(u\i)a=\g+\nu+\hz.\tag g$$
Since $a\in\tfl'_\et$, we have $\Ad(u\i)a-a\in\fu'_\et$ for $u\in U_{P'_0}$. 
Hence in (g) we have $\g+\nu=a$ and $\hz=\Ad(u\i)a-a$. In particular, $\g,\nu$
are uniquely determined. Thus, our variety may be identified with 
$U_{P'_0}\T(\fu'_N\cap\fu''_N)$. This completes the proof of (d), hence that 
of (b).

\subhead 5.5\endsubhead
From the results in 5.3 and 5.4 we can deduce, using the argument in
\cite{\GRA, 8.9} (based on \cite{\GRA, 1.4}), the following result.

\proclaim{Proposition 5.6} We have $K\in\cq(\fl'_\et)$; moreover, we have
(noncanonically) $K\cong\op_\Om K_\Om$, where $\Om$ runs over good $(P'_{0}, P_{0})$-double cosets in $G_{\un0}$. 
\endproclaim

\head 6. Spiral restriction\endhead
We introduce the spiral restriction functor which is adjoint to the spiral induction. The main result 
in this section is Proposition 6.4, which calculates the inner product $\{,\}$ (in the sense of 0.12) of two spiral inductions.

\subhead 6.1. Definition of spiral restriction \endsubhead
In addition to $\et\in\ZZ-\{0\}$ which has been fixed in 2.9, in this section we fix $\e',\e''$ 
in $\{1,-1\}$. Let 
$(\fp'_*,L',P'_0,\fl',\fl'_*,\fu'_*)\in\fP^{\e'}$. Let 
$\p':\fp'_\et@>>>\fl'_\et$ be the obvious map. For any $B\in\cd(\fg_{\da})$ 
we set
$${}^{\e'}\Res_{\fp'_\et}^{\fg_\da}(B)=\p'_!(B|_{\fp'_\et})
\in\cd(\fl'_\et).$$
We show:

(a) {\it If $B\in\cq^{\e''}_\et(\fg_\da)$ then
${}^{\e'}\Res_{\fp'_\et}^{\fg_\da}(B)\in\cq(\fl'_\et)$.}
\nl
To prove this we can assume that $B$ is in addition a simple perverse sheaf.
Then, using the definition of $\cq^{\e''}_\et(\fg_\da)$, we see that it is 
enough to prove (a) in the case where 
$B={}^{\e''}\Ind_{\fp_\et}^{\fg_{\da}}(\cl^\sha)$, 
with $(\fp_*,L,P_0,\fl,\fl_*,\fu_*)\in\fP^{\e''}$, $\cl^\sha$ as in 5.2. In this 
case, (a) follows from 5.6.

\mpb

We thus have a functor
${}^{\e'}\Res_{\fp'_\et}^{\fg_\da}:\cq^{\e''}_\et(\fg_\da)@>>>
\cq(\fl'_\et)$ called {\it spiral restriction}.

We have the following result.
\proclaim{Proposition 6.2. (Adjunction)}Let $C\in\cq(\fl'_\et)$, and let 
$B\in\cq^{\e''}_\et(\fg_\da)$. For any $j\in\ZZ$ we have
$$d_j(\fl'_\et;C,{}^{\e'}\Res_{\fp'_\et}^{\fg_\da}(B))=
d_{j'}(\fg_\da;{}^{\e'}\Ind_{\fp'_\et}^{\fg_\da}(C),B)\tag a$$
where $j'=j+2\dim\fu'_0$.
\endproclaim
The proof is almost a copy of that of \cite{\GRA, 9.2}. We omit it.

\mpb

For $B\in\cd(\fg_\da)$ we set
$${}^{\e'}\tRes_{\fp'_\et}^{\fg_\da}(B)=
{}^{\e'}\Res_{\fp'_\et}^{\fg_\da}(B)[\dim\fu'_\et-\dim\fu'_0].$$
With this notation, the equality (a) can be reformulated without a shift from
$j$ to $j'$ as follows.
$$d_j(\fl'_\et;C,{}^{\e'}\tRes_{\fp'_\et}^{\fg_\da}(B))=
d_j(\fg_\da;{}^{\e'}\tInd_{\fp'_\et}^{\fg_\da}(C),B).\tag b$$

\subhead 6.3\endsubhead
Let $(\fp'_*,L',P'_0,\fl',\fl'_*,\fu'_*)\in\fP^{\e'}$,
$(\fp_*,L,P_0,\fl,\fl_*,\fu_*)\in\fP^{\e''}$. Let $A\in\cq(\fl_\et)$, 
$A'\in\cq(\fl'_\et)$ be cuspidal perverse sheaves. As in 4.3 we have 
$A=\cl^\sha[\dim\fl_\et]$, $A'=\cl'{}^\sha[\dim\fl'_\et]$ where $\cl$ (resp. 
$\cl'$) is a local system on $\ovsc\fl_\et$ (resp. $\ovsc\fl'_\et$).

We denote by $X$ the set of all $g\in G_{\un0}$ such that the $\e''$-spiral
$\{\Ad(g)\fp_N;N\in\ZZ\}$ 
and the $\e'$-spiral $\fp'_*$ have a common 
splitting. If $g\in X$ there is a unique isomorphism of $\ZZ$-graded Lie 
algebras $\l_g:\fl@>>>\fl'$ such that the compositions

$\Ad(g)\fp_N\cap\fp'_N@>>>\fp'_N@>>>\fl'_N$,

$\Ad(g)\fp_N\cap \fp'_N@>\Ad(g\i)>>\fp_N@>>>\fl_N@>\l_g>>\fl'_N$  
\nl
coincide for any $N$ (the unnamed maps are the obvious imbeddings or 
projections). Moreover, $\l_g$ is induced by an isomorphism $L@>>>L'$. Let 
$X'$ be the set of all $g\in X$ such that $\l_g:\fl_\et@>\si>>\fl'_\et$ carries 
$\cl$ to the dual of $\cl'$. For any $g\in X'$ we set
$$\t(g)=-\dim\fra{\fu'_0+\Ad(g)\fu_0}{\fu'_0\cap\Ad(g)\fu_0}+
\dim\fra{\fu'_\et+\Ad(g)\fu_\et}{\fu'_\et\cap\Ad(g)\fu_\et}.$$
Note that both $X$ and $X'$ are unions of $(P'_0,P_0)$-double cosets in 
$G_{\un0}$ and that $\t(g)$ depends only on the double coset of $g$. We have 
the following result.

\proclaim{Proposition 6.4}Let 
$$\Pi=\sum_{j\in\ZZ}d_j(\fg_\da;{}^{\e'}\tInd_{\fp'_\et}^{\fg_\da}(A'),
{}^{\e''}\tInd_{\fp_\et}^{\fg_\da}(A))v^{-j}\in\NN((v)).$$
We have
$$\Pi=(1-v^2)^{-r}\sum_{g_0}v^{\t(g_0)}$$
where $r$ is the dimension of the centre of $\fl$ and the sum is taken over a 
set of representatives $g_0$ for the $(P'_0,P_0)$-double cosets in $G_{\un0}$ 
that are contained in $X'$. In particular, if $\Pi\ne0$ then $X'\ne\emp$.
\endproclaim
Using 6.2, we have
$$\Pi=\sum_{j\in\ZZ}d_j(\fl'_\et;A',{}^{\e'}\tRes_{\fp'_\et}^{\fg_\da}
({}^{\e''}\tInd_{\fp_\et}^{\fg_\da}(A)))v^{-j}
=\sum_{j\in\ZZ}d_{j+s}(\fl'_\et;A',K)v^{-j}$$ 
where
$s=\dim\fu_0+\dim\fu_\et+\dim\fu'_\et-\dim\fu'_0+\dim\fl_\et$ and
$$K={}^{\e'}\Res_{\fp'_\et}^{\fg_\da}
({}^{\e''}\Ind_{\fp_\et}^{\fg_\da}(\cl^\sha))$$ 
is as in 5.2. Using the description of $K$ in 5.3(a), 5.4(b), 5.6, we see that
$$\Pi=\sum_{j\in\ZZ}\sum_gQ_j(g)v^{-j+s-2f(g)}\tag a$$
where $g$ runs over a set of representatives for the $(P'_0,P_0)$-double 
cosets in $G_{\un0}$ which are good (see 5.2) and 
$$Q_j(g)=d_j(\tfl'_\et;A',\ind_{\tfq'_\et}^{\tfl'_\et}(\cl''{}^\sha)),$$
$$f(g)=\dim(\fu'_0/(\fu'_0\cap\Ad(g)\fp_0)+\dim(\fu'_\et\cap\Ad(g)\fu_\et);$$
the following notation is used:

$\tfl'_*$ is a certain splitting of $\fp'_*$, $\tfl''_*$ is a certain 
splitting of $\{\Ad(g)\fp_N;N\in\ZZ\}$, $\tfq'=\op_{N\in\ZZ}\tfq'_N$ (where 
$\tfq'_N=\tfl'_N\cap\Ad(g_0)\fp_N$)
is a parabolic subalgebra of $\tfl'=\op_N\tfl'_N$ whose with Levi subalgebra 
$\tfl''=\op_N\tfl''_N$; $A'$ is viewed as an object of $\cq(\tfl'_\et)$ via the
obvious isomorphism $\tfl'_\et@>>>\fl'_\et$ and $\cl''{}^\sha\in\cq(\tfl''_\et)$ 
corresponds to $\cl^\sha$ via the isomorphism
$\fl_\et@>\Ad(g)>>\Ad(g)\fp_\et/\Ad(g)\fu_\et=\tfl''_\et$. 

By the implication (a)$\imp$(c) in \cite{\GRA, 10.6}, we have $Q_j(g)=0$ 
unless $\tfq'=\tfl'$. In this case, since $\tfl''$ is a Levi subalgebra of 
$\tfq'$, we must have $\tfl'=\tfl''$ so that $g\in X$. Conversely, if 
$g\in X$, then the $(P'_0,P_0)$-double coset of $g$ is good. Indeed, let 
$\tfl'_*$ be a splitting of $\fp'_*$ which is also a splitting for 
$\{\Ad(g)\fp_N;N\in\ZZ\}$. We have 
$$\Ad(g)\fp_N=\tfl'_N\op\Ad(g)\fu_N\sub(\fp'_N\cap\Ad(g)\fp_N)
+\Ad(g)\fu_N\sub\Ad(g)\fp_N$$
and our claim follows. Thus the sum in (a) can be taken over a set of 
representatives $g$ for the $(P'_0,P_0)$-double cosets in $G_{\un0}$ that are 
contained in $X$ and for such $g$ we have
$Q_j(g)=d_j(\tfl'_\et;A',\cl''{}^\sha)$ where $\tfl'=\tfl''$, 
$\cl''{}^\sha\in\cq(\tfl''_\et)$ are as above. Using \cite{\GRA, 15.1}, we see 
that in the sum over $g$ in (a) we can take $g\in X'$ and that 
the contribution of such $g$ to the sum is $(1-v^2)^{-r}v^{s-2f(g)-d}$ where 
$d=\dim\fl_\et$. It remains to show that for $g$ as above we have 
$s-2f(g)-d=\t(g)$. It is enough to show that 

(b) $\fu'_0\cap\Ad(g)\fp_0=\fu'_0\cap\Ad(g)\fu_0$,

(c) $\dim(\Ad(g)\fu_0)=\dim\fu'_0$.
\nl
Now (b),(c) hold since $\Ad(g)\fp_0$, $\fp'_0$ are parabolic subalgebras of 
$\fg_{\un0}$ with nilradicals $\Ad(g)\fu_0$, $\fu'_0$ and with a common Levi 
subalgebra. This completes the proof of the proposition.

\subhead 6.5\endsubhead
In the special case where
$$(\fp'_*,L',P'_0,\fl',\fl'_*,\fu'_*)=(\fp_*,L,P_0,\fl,\fl_*,\fu_*)$$
and $A'\cong D(A)$, the sum $\sum_gv^{\t(g)}$ in Proposition 6.4 is over a 
nonempty set of $g$ (we have $1\in X'$) hence the sum is nonzero and $\Pi$ in 
6.4 is nonzero. In particular, we see that 
$${}^{\e''}\tInd_{\fp_\et}^{\fg_\da}(A)\ne0.\tag a$$

\subhead 6.6. The map $\ps$ from simple perverse sheaves to $\un\fT_{\et}$ \endsubhead
Let $B$ be a simple perverse sheaf in $\cq^{\e''}_\et(\fg_\da)$. We associate to
$B$ an element of $\un\fT_\et$ (see 3.5) as follows. We can find 
$(\fp_*,L,P_0,\fl,\fl_*,\fu_*)\in\fP^{\e''}$ and $A$ as in 6.3 such that
$${}^{\e''}\tInd_{\fp_\et}^{\fg_\da}(A)\cong B[d]\op C$$
where $d\in\ZZ$ and $C\in\cq^{\e''}_\et(\fg_\da)$. Let $\tfl_*$ be a splitting 
of $\fp_*$. Let $\tfl=\op_N\tfl_N$, $\tL=e^{\tfl}\sub G$, 
$\tL_0=e^{\tfl_0}\sub G$ and let $\tC$ be the simple perverse sheaf on 
$\tfl_\et$ corresponding to $A$ under the obvious isomorphism 
$\tfl_\et@>\si>>\fl_\et$. Then $(\tL,\tL_0,\tfl,\tfl_*,\tC)$ is an object of 
$\fT_\et$ and its $G_{\un0}$-orbit is independent of the choice of splitting, 
by 2.7(a). Now let $(\fp'_*,L',P'_0,\fl',\fl'_*,\fu'_*)\in\fP^{\e'}$, $A'$ be 
as in 6.3 (with $\e'=\e'')$ and assume that 
$${}^{\e'}\tInd_{\fp'_\et}^{\fg_\da}(A')\cong B[d']\op C'$$
where $d'\in\ZZ$ and $C'\in\cq^{\e''}_\et(\fg_\da)$. We choose a splitting 
$\tfl'_*$ of $\fp'_*$ and we associate to it a system 
$(\tL',\tL'_0,\tfl',\tfl'_*,\tC')$ just as $(\tL,\tL_0,\tfl,\tfl_*,\tC)$ 
was defined in terms of $\tfl$; here $\tC'$ corresponds to $A'$. Using 4.1(d),
we see that
$${}^{\e'}\tInd_{\fp'_\et}^{\fg_\da}(D(A'))\cong D(B)[-d']\op D(C').$$
Let $\Pi$ be as in 6.4 (with $A'$ replaced by $D(A')$ and $\e'=\e''$). From 
the definition of $\Pi$ in 6.4 we have also 
$$\Pi=\{B[d]\op C,D(B)[-d']\op D(C')\}=v^{d-d'}\text{ plus an element in }
\NN((v)).$$
(We use 0.12.) In particular we have $\Pi\ne0$ hence $X'$ in 6.4 is nonempty.
It follows that $(\tL',\tL'_0,\tfl',\tfl'_*,\tC')$ and 
$(\tL,\tL_0,\tfl,\tfl_*,\tC)$ are in the same $G_{\un0}$-orbit. This proves 
that $B\m(\tL,\tL_0,\tfl,\tfl_*,\tC)$ associates to $B$ a well-defined element
$\psi(B)\in\un\fT_\et$.

\subhead 6.7\endsubhead
For any $\x\in\un\fT_\et$ let ${}^\x\cq^{\e'}_\et(\fg_\da)$ be the full 
subcategory of $\cq^{\e'}_\et(\fg_\da)$ whose objects are direct sums of shifts
of simple perverse sheaves $B$ in $\cq^{\e'}_\et(\fg_\da)$ such that
$\psi(B)=\x$ (see 6.6); let ${}^\x\ck^{\e'}_\et(\fg_\da)$ be the (free) 
$\ca$-submodule of $\ck^{\e'}_\et(\fg_\da)$ with basis given by the simple 
perverse sheaves $B$ in ${}^\x\cq^{\e'}_\et(\fg_\da)$. Clearly, we have
$$\ck^{\e'}_\et(\fg_\da)=\op_{\x\in\un\fT_\et}{}^\x\ck^{\e'}_\et(\fg_\da).$$ 

\head 7. The categories $\cq(\fg_\da)$, $\cq'(\fg_\da)$\endhead
In this section we consider two categories of perverse sheaves $\cq(\fg_\da)$, $\cq'(\fg_\da)$ defined
in terms of spiral induction. 
The simple objects in $\cq(\fg_\da)$ are supported on $\fg^{nil}_{\da}$, while those in
$\cq'(\fg_\da)$ have Fourier-Deligne transforms  supported on $\fg^{nil}_{\da}$. 
We also complete the proof of the main theorem 0.6.

\subhead 7.1\endsubhead  
Let $(\co,\cl)\in\ci(\fg_\da)$. Let
$A_1$ be the simple perverse sheaf on $\fg_\da$ such that $\supp(A_1)$ is the 
closure $\bco$ of $\co$ in $\fg_\da$ and $A_1|_\co=\cl[\dim\co]$.

Choose $x\in\co$ and $\ph\in J_\da(x)$; define $\fp^x_*,\tfl^\ph_*,\tL^\ph,P_0$
as in 2.9. Then $\cq(\tfl^\ph_\et)$ is defined in terms of $\tfl^\ph_*,\tL^\ph$
and for any $A'\in\cq(\tfl^\ph_\et)$ we can consider
$$I(A'):={}^{\doet}\Ind_{\fp^x_\et}^{\fg_\da}(A')\in\cq^{\doet}_\et(\fg_\da),$$
see 4.1. We show:

(a) {\it If $A'\in\cq(\tfl^\ph_\et)$, then the support of $I(A')$ is contained 
in $\bco$.}
\nl
Let $y\in\fg_\da$ be in the support of $I(A')$. We must show that $y\in\bco$.
From the definition of $I(A')$, there exists $g\in G_{\un0}$ and 
$z\in\fp^x_\et$ such that $\Ad(g)(z)=y$. 
Since the support of $I(A')$ and $\bco$
are $G_{\un0}$-invariant we may replace $y$ by $\Ad(g\i)y$ hence we may assume
that $y\in\fp^x_\et$. Using 2.9(e), we see that $\fp^x_\et$ is equal to the
closure of the $P_0$-orbit of $x$ in $\fp^x_\et$, which is clearly contained
in $\bco$. This proves (a).  

\mpb

Recall that $x\in\ovsc\tfl^\ph_\et$ (see 2.9(b)) hence 
$\ovsc\tfl^\ph_\et\sub\co$. By 2.9(c), $\cl_1:=\cl|_{\ovsc\tfl^\ph_\et}$ is an 
irreducible $\tL^\ph_0$-equivariant local system on $\ovsc\tfl^\ph_\et$. Let 
$\cl_1^\sha\in\cd(\tfl^\ph_\et)$ be as in 0.11 and let 
$A=\cl_1^\sha[\dim\tfl^\ph_\et]$. We show:

(b) {\it $I(\cl_1^\sha)|_\co$ is $\cl$.}
\nl
Let $E'_\co$ be the inverse image of $\co$ under $c:E'@>>>\fg_\da$ (where 
$c,E'$ are as in 4.1 with $\fp_*=\fp_*^x$, $\e=\doet$). From the definitions we see that 
it is enough to check that the map $c_{\co}:E'_\co@>>>\co$ (restriction of $c$) is bijective on 
$\kk$-points. 
Since $G_{\un0}$ 
acts naturally on both $E'_\co$ and $\co$ compatibly with $c$ and the action 
on $\co$ is transitive, it suffices to check that $c\i(x)$ is a single point, 
namely $(P_0,x)$. Let $(gP_0,x)\in c\i(x)$. We have $g\in G_{\un0}$, 
$\Ad(g\i)x\in\fp^x_\et$ hence $x\in\Ad(g)\fp^x_\et$. From 2.9(d) we deduce that 
$g\in P_0$ hence $(gP_0,x)=(P_0,x)$. This proves (b).

We show:

(c) {\it $I(\cl_1^\sha)$ is isomorphic to $\op_{j=1}^rA_j[t_j]$ where 
$t_1=-\dim\co$ and for any $j\ge2$, $A_j$ is a simple $G_{\un0}$-equivariant 
perverse sheaf on $\fg_{\da}$ 
with support contained in $\bco-\co$ and $t_j\in\ZZ$.}
\nl
This follows from the fact that $I(\cl_1^\sha)$ is a semisimple 
$G_{\un0}$-equivariant perverse sheaf on $\fg_\da$ (the decomposition theorem), taking into account (a),(b). 

\mpb

By 1.5(a) we can find a parabolic subalgebra $\fq$ of $\tfl^\ph$, a Levi
subalgebra $\fm$ of $\fq$ (with $\fq,\fm$ compatible with the $\ZZ$-grading of
$\tfl^\ph$) and a cuspidal $M_{0}:=e^{\fm_0}$-equivariant perverse sheaf 
$C$ on $\fm_\et$ such that some shift of $A$ is a direct summand of 
$\ind_{\fq_\et}^{\tfl^\ph_\et}(C)$. From the definition we have 
$$\Ps(\co,\cl)=(M,M_{0},\fm,\fm_*,C)\in\un\fT_\et,\tag d$$
where $M=e^{\fm}$, see 3.5.  

For any $N\in\ZZ$ let $\hfp_N$ be the inverse image of $\fq_N$ under the 
obvious map $\fp_N@>>>\fl_N$. Then by 2.8(a), 
$\hfp_*$ is an $\doet$-spiral and $\fm_*$ is 
a splitting of it, so that, by 4.2(a), we have
$${}^{\doet}\Ind_{\hfp_\et}^{\fg_\da}(C)
={}^{\doet}\Ind_{\fp^\ph_\et}^{\fg_\da}(\ind_{\fq_\et}^{\tfl^\ph_\et}(C)).$$
It follows that some shift of ${}^{\doet}\Ind_{\fp^\ph_\et}^{\fg_\da}(A)$
is a direct summand of ${}^{\doet}\Ind_{\hfp_\et}^{\fg_\da}(C)$ hence, using 
(c), we see that some shift of $A_1$  is a direct summand of  
${}^{\doet}\Ind_{\hfp_\et}^{\fg_\da}(C)$. 
In particular we have $A_1\in\cq^{\doet}_\et(\fg_\da)$ and
$\psi(A_1)=(M,M_{0},\fm,\fm_*,C)\in\un\fT_\et$, see 6.6 (with $\e=\doet$). Comparing with (d) we see that

(e) $\psi(A_1)=\Ps(\co,\cl)$.

\subhead 7.2. Characterization of $\cq^{\doet}_\et(\fg_\da)$ as orbital sheaves\endsubhead
Let $A'$ be a semisimple $G_{\un0}$-equivariant complex on $\fg_\da$. We show: 

(a) {\it We have $A'\in\cq^{\doet}_\et(\fg_\da)$ if and only if
$\supp(A')\sub\fg_\da^{nil}$.}
\nl
We can assume that $A'$ is a simple perverse sheaf. If 
$\supp(A')\sub\fg_\da^{nil}$, then we have $A'\in\cq^{\doet}_\et(\fg_\da)$ by the 
arguments in 7.1. Conversely, assume that $A'\in\cq^{\doet}_\et(\fg_\da)$.
We can find $(\fp_*,L,P_0,\fl,\fl_*,\fu_*)\in\fP^{\doet}$ and 
$A\in\cq(\fl_\et)$ such that some shift of $A'$ is a direct summand of
$B:={}^{\doet}\Ind_{\fp_\et}^{\fg_\da}(A)$. 
To show that $\supp(A')\sub\fg_\da^{nil}$ it is enough to show that
$\supp(B)\sub\fg_\da^{nil}$ or (with $c,A_1$ as in 4.1 with $\e=\doet$) that
$\supp(c_!A_1)\sub\fg_\da^{nil}$. This would follow if we can show that the 
image of $c$ is contained in $\fg_\da^{nil}$. By the definition of $c$ it is 
enough to show that $\fp_\et\sub\fg_\da^{nil}$. This follows from 2.5(d) 
applied with $N=\et$.

We now restate 7.1(e) as follows.

(b) {\it Let $A'$ be a simple perverse sheaf in $\cq^{\doet}_\et(\fg_\da)$ and let 
$(\co,\cl)\in\ci(\fg_\da)$ be such that $\supp(A')=\bco$ and 
$A'|_{\co}=\cl[\dim\co]$. Then $\psi(A')=\Ps(\co,\cl)$. (Notation of 3.5 and 
6.6 with $\e=\doet$.)}

\subhead 7.3\endsubhead 
We now give another proof of the following statement (see also 3.8(f)):

(a) {\it The map $\Ps:\ci(\fg_\da)@>>>\un\fT_\et$ in 3.5 is surjective.}
\nl
Let $(M,M_0,\fm,\fm_*,C)$ be an element of $\fT_\et$. We can find an 
$\doet$-spiral $\fp_*$ such that $\fm_*$ is a splitting of $\fp_*$. By 6.5(a), we
have ${}^{\doet}\tInd_{\fp_\et}^{\fg_\da}(C)\ne0$, that is, there exists a 
simple perverse sheaf $A'$ in $\cq^{\doet}_\et(\fg_\da)$ such that some shift of $A'$
is a direct summand of ${}^{\doet}\tInd_{\fp_\et}^{\fg_\da}(C)$. It follows that
$\psi(A')=(M,M_0,\fm,\fm_*,C)$ hence, by 7.2(b), we have 
$\Ps(\co,\cl)=(M,M_0,\fm,\fm_*,C)$ where $(\co,\cl)$ corresponds to $A'$ as in
7.2(b). This proves (a).

\subhead 7.4\endsubhead
Until the end of 7.7 we assume that $p>0$.
If $E,E'$ are finite dimensional $\kk$-vector space with a given
perfect bilinear pairing $E\T E'@>>>\kk$ then we have the Fourier-Deligne
transform functor $\Ph:\cd(E)@>>>\cd(E')$ defined in terms of a fixed
nontrivial character $\FF_p@>>>\bbq^*$ as in \cite{\GRA, 1.9}. 

\subhead 7.5. Fourier transform and spiral restriction\endsubhead
Let $B\in\cd(\fg_\da)$; we denote by $\Ph_\fg(B)\in\cd(\fg_{-\da})$ the 
Fourier-Deligne transform of $B$ with respect to the perfect pairing 
$\fg_\da\T\fg_{-\da}@>>>\kk$ defined by $\la,\ra$.

Let $\e'\in\{1,-1\}$. Let $(\fp'_*,L',P'_0,\fl',\fl'_*,\fu'_*)\in\fP^{\e'}$ and let
$$R_\et={}^{\e'}\Res_{\fp'_\et}^{\fg_\da}(B)\in\cd(\fl'_\et),\qua
R_{-\et}={}^{\e'}\Res_{\fp'_{-\et}}^{\fg_{-\da}}(\Ph_\fg(B))\in\cd(\fl'_{-\et}).$$
Then 

(a) {\it $R_{-\et}$ is the Fourier-Deligne transform of $R_\et$ with respect to 
the perfect pairing $\fl_\et\T\fl_{-\et}@>>>\kk$ defined by $\la,\ra$.}
\nl
The proof is almost the same as that of \cite{\GRA, 10.2}. We omit it.

\subhead 7.6. Fourier transform and spiral induction \endsubhead
Let $\e'\in\{1,-1\}$. Let 

$(\fp'_*,L',P'_0,\fl',\fl'_*,\fu'_*)\in\fP^{\e'}$.
\nl
Let $A\in\cd(\fl'_\et)$ be a semisimple complex; we denote by 
$\Ph_{\fl'}(A)\in\cd(\fl'_{-\et})$ the Fourier-Deligne transform of $A$ with 
respect to the perfect pairing $\fl'_\et\T\fl'_{-\et}@>>>\kk$ defined by 
$\la,\ra$; note that $\Ph_{\fl'}(A)$ is a semisimple complex. Let
$$\align&I_\et={}^{\e'}\tInd_{\fp'_\et}^{\fg_\da}(A)\in\cd(\fg_\da),\\&
I_{-\et}={}^{\e'}\tInd_{\fp'_{-\et}}^{\fg_{-\da}}(\Ph_{\fl'}(A))
\in\cd(\fg_{-\da}).\endalign$$
Then 

(a) {\it $I_{-\et}$ is the Fourier-Deligne transform of $I_\et$ with respect to 
the perfect pairing $\fg_\da\T\fg_{-\da}@>>>\kk$ defined by $\la,\ra$.}
\nl
The proof is almost the same as that of \cite{\ANTIORB, A2}. We omit it.

\subhead 7.7. Characterization of $\cq^{-\doet}_\et(\fg_\da)$ as anti-orbital sheaves\endsubhead
Let $B\in\cd(\fg_\da)$ be a semisimple complex; let
$B'=\Ph_\fg(B)\in\cd(\fg_{-\da})$ be its Fourier-Deligne transform, as in 7.5.
Note that $B'$ is again a semisimple complex. We show: 

(a) {\it We have $B\in\cq^{-\doet}_\et(\fg_\da)$ if and only if
$\supp(B')\sub\fg_{-\da}^{nil}$.}
\nl
We can assume that $B$ (and hence also $B'$) is a simple perverse sheaf.

Assume first that $B\in\cq^{-\doet}_\et(\fg_\da)$. We can find
$(\fp'_*,L',P'_0,\fl',\fl'_*,\fu'_*)\in\fP^{-\doet}$ and a cuspidal
perverse sheaf $C$ in $\cq(\fl'_\et)$ such that some shift of $B$ is a direct
summand of ${}^{-\doet}\tInd_{\fp'_\et}^{\fg_\da}(C)$.
Using 7.6(a) we see that some shift of $B'$ is a direct summand of
${}^{-\doet}\tInd_{\fp'_{-\et}}^{\fg_{-\da}}(C')$
where $C'=\Ph_{\fl'}(C)\in\cd(\fl'_{-\et})$ (notation of 7.6). By 
\cite{\GRA, 10.6}, $C'$ is a cuspidal perverse sheaf in $\cq(\fl'_{-\et})$.
It follows that $B'\in\cq^{-\doet}_{-\et}(\fg_{-\da})$. Using 7.2(a) (with $\et,\da$
replaced by $-\et,-\da$) we deduce that $\supp(B')\sub\fg_{-\da}^{nil}$.

Conversely, assume that $B$ is such that $\supp(B')\sub\fg_{-\da}^{nil}$. Using 
7.2(a), we see that $B'\in\cq^{-\doet}_{-\et}(\fg_{-\da})$. We can find
$(\fp'_*,L',P'_0,\fl',\fl'_*,\fu'_*)\in\fP^{-\doet}$ and a cuspidal
perverse sheaf $C'_1$ in $\cq(\fl'_{-\et})$ such that some shift of $B'$ is a 
direct summand of ${}^{-\doet}\tInd_{\fp'_{-\et}}^{\fg_{-\da}}(C'_1)$.
We can find a cuspidal perverse sheaf $C_1$ in $\cq(\fl'_\et)$ such that 
$C'_1=\Ph_{\fl'}(C)$ (we use again \cite{\GRA, 10.6}). Using 7.6(a), we see 
that some shift of $\Ph_{\fg}(B)$ is a direct summand of
$\Ph_{\fg}({}^{-\doet}\tInd_{\fp'_\et}^{\fg_\da}(C_1))$
hence some shift of $B$ is a direct summand of
${}^{-\doet}\tInd_{\fp'_\et}^{\fg_\da}(C_1)$
so that $B\in\cq^{-\doet}_\et(\fg_\da)$. This completes the proof of (a).

\subhead 7.8\endsubhead
The assumption on $p$ in 7.4 is no longer in force.
From 7.2(a) we see that $\cq^{\doet}_\et(\fg_\da)$ (hence also $\ck^{\doet}_\et(\fg_\da)$)
is independent of $\et$ as long as $\un\et=\da$. 
We shall write $\cq(\fg_\da),\ck(\fg_\da)$ instead of 
$\cq^{\doet}_\et(\fg_\da)$, $\ck^{\doet}_\et(\fg_\da)$.
From 7.7(a) we see that $\cq^{-\doet}_\et(\fg_\da)$ (hence also 
$\ck^{-\doet}_\et(\fg_\da)$) is independent of $\et$ as long as $\un\et=\da$ (at least
when $p>0$, but then the same holds for $p=0$ by standard
arguments.) We shall write $\cq'(\fg_\da),\ck'(\fg_\da)$ instead of 
$\cq^{-\doet}_\et(\fg_\da)$, $\ck^{-\doet}_\et(\fg_\da)$. For $\x\in\un\fT_\et$ we write 
${}^\x\cq(\fg_\da),{}^\x\ck(\fg_\da)$ instead of ${}^\x\cq^{\doet}_\et(\fg_\da)$, 
${}^\x\ck^{\doet}_\et(\fg_\da)$ and we write ${}^\x\cq'(\fg_\da),{}^\x\ck'(\fg_\da)$ 
instead of ${}^\x\cq^{-\doet}_\et(\fg_\da)$, ${}^\x\ck^{-\doet}_\et(\fg_\da)$.

\subhead 7.9. Proof of Theorem 0.6 \endsubhead
Let $\x\in\un\fT_\et$. 
Let $K\in\cd_{G_{\un0}}(\fg_\da^{nil})$. We say that
$K\in\cd_{G_{\un0}}(\fg_\da^{nil})_\x$ if
any simple perverse sheaf $B$ 
which appears in a perverse cohomology sheaf of $K$ 
satisfies $\psi(B)=\x$; note that $B$ belongs to $\cq^{\doet}_\et(\fg_\da)$, see 
7.2(a); hence $\psi(B)$ is defined as in 6.6. 

Now let $\x,\x'$ in $\un\fT_\e$ be such that $\x\ne\x'$. Let 
$K\in\cd_{G_{\un0}}(\fg_\da^{nil})_\x$, 
$K'\in\cd_{G_{\un0}}(\fg_\da^{nil})_{\x'}$. We show:

(a) $\Hom_{\cd_{G_{\un0}}(\fg_\da^{nil})}(K,K')=0$.  
\nl
We can assume that $K=B[n],K'=B'[n']$ where
$B,B'$ are simple perverse sheaves in 
$\cq^{\doet}_\et(\fg_\da)$ such that $\psi(B)=\x,\psi(B')=\x'$ and $n,n'$ are 
integers. We see that it is enough to
prove (a) in the case where $K={}^{\doet}\tInd_{\fp'_\et}^{\fg_\da}(A')[n]$,
$K'={}^{\doet}\tInd_{\fp_\et}^{\fg_\da}(A)[n']$  
with $n,n'\in\ZZ$, $\fp_*,\fp'_*,A,A'$ as in 6.4, and $\e'=\e''=\doet$, 
since some shifts of $B$ and $B'$ appear as direct summands of such $K$ and $K'$.  By 0.12(a),  we have an isomorphism 
$$\Hom_{\cd_{G_{\un0}}(\fg_\da^{nil})}(K,K')=\DD_0(\fg_\da^{nil}, G_{\un0}; K, D(K'))^*.$$
Hence

(b) $\dim\Hom_{\cd_{G_{\un0}}(\fg_\da^{nil})}(K,K')=
d_{n-n'}(\fg_\da^{nil};{}^{\doet}\tInd_{\fp'_\et}^{\fg_\da}(A'),{}^{\doet}\tInd_{\fp_\et}^{\fg_\da}(D(A)))$.
\nl
Here we use 4.1(d). Since $\x\ne\x'$, the set $X'$ defined in 6.4 for the pair $(D(A),A')$ is empty. 
Therefore the right side of (b) is zero by 6.4. Then (a) follows from (b). We see that Theorem 0.6 holds.

\head 8. Monomial and quasi-monomial objects\endhead
The results in this section are parallel to those in 1.8-1.9. 
They serve as preparation for the next section.

\subhead 8.1\endsubhead
Let $\e=\doet$. We denote by $\fR^\e$ the set of all data of the form
$(\fp_*,L,P_0,\fl,\fl_*,\fu_*,A)$ where 
$(\fp_*,L,P_0,\fl,\fl_*,\fu_*)\in\fP^\e$
(see 4.1) and $A$ is a perverse sheaf in $\cq(\fl_\et)$ which is
$\et$-semicuspidal (as in 1.8 with $H$ replaced by $L$).

\subhead 8.2\endsubhead
An object $B\in\cq(\fg_\da)$ is said to be {\it $\et$-quasi-monomial} if 
$B\cong\tInd_{\fp_\et}^{\fg_\da}(A)$ for some
$(\fp_*,L,P_0,\fl,\fl_*,\fu_*,A)\in\fR^\e$; if in addition $A$ is taken to be 
cuspidal, then $B$ is said to be {\it $\et$-monomial}. Using 1.8(b) and the 
transitivity property 4.2, we see that 

(a) {\it If $B\in\cq(\fg_\da)$ is $\et$-quasi-monomial then there exists an 
$\et$-monomial object $B'\in\cq(\fg_\da)$ 
such that $B'\cong B[a_1]\op B[a_2]\op \cdots\op B[a_k]$ for some sequence $a_1,a_2,\do,a_k$ in
$\ZZ$, $k\ge1$. In particular, in $\ck(\fg_\da)$ we have
$(B')=(v^{a_1}+\do+v^{a_k})(B)$.}
\nl
An object of $\cq(\fg_\da)$ is said to be {\it $\et$-good} if it is a direct sum 
of shifts of $\e$-quasi-monomial objects.

\proclaim{Proposition 8.3} Let $B\in\cq(\fg_\da)$. There exists $\et$-good
objects $B_1,B_2$ in $\cq(\fg_\da)$ such that $B\op B_1\cong B_2$.
\endproclaim
We can assume that $B$ is a simple perverse sheaf. We define 
$(\co,\cl)\in\ci(\fg_\da)$ by the requirement that $\supp B$ is the closure 
$\bco$ of $\co$ in $\fg_\da$ and $B|_\co=\cl[\dim\co]$. We prove the 
proposition by induction on $\dim\co$. Let $x\in\co$. We associate to $x$ an 
$\e$-spiral $\fp_*=\fp_*^x$ as in 2.9; we complete it uniquely to a system
$(\fp_*,L,P_0,\fl,\fl_*,\fu_*)\in\fP^\e$. By 7.1(c), there exists 
$A_1\in\cq(\fl_\et)$ such that 
${}^\e\Ind_{\fp_\et}^{\fg_\da}(A_1)\cong B[d]\op B'$
where $d\in\ZZ$ and $B'\in\fq(\fg_\da)$ has support contained in $\bco-\co$. 
We now use 1.9(a) for $L,A_1$ instead of $H,A_1$; applying 
${}^\e\Ind_{\fp_\et}^{\fg_\da}$ to the equality in 1.9(a) we obtain
$${}^\e\Ind_{\fp_\et}^{\fg_\da}(A_{1})\op C'_1\op C'_2\op\do\op C'_t=
C'_{t+1}\op\do\op C'_{t+t'}$$  
where each $C'_j$ is an $\et$-quasi-monomial object with a shift (we have used 
the transitivity property 4.2). Thus we have
$$B[d]\op B'\op C'_1\op C'_2\op\do\op C'_t=C'_{t+1}\op\do\op C'_{t+t'}.$$ 
Now the induction hypothesis implies that $B'$ is $\et$-good. From this and 
the previous equality we see that $B$ is $\et$-good. The proposition is proved.

\proclaim{Corollary 8.4}(a) The $\ca$-module $\ck(\fg_\da)$ is generated by the
classes of \lb $\et$-quasi-monomial objects of $\cq(\fg_\da)$.

(b) The $\QQ(v)$-vector space $\QQ(v)\ot_\ca\ck(\fg_\da)$ is generated by
the classes of $\et$-monomial objects of $\cq(\fg_\da)$. 
\endproclaim
(a) follows immediately from 8.3; (b) follows from (a) using 8.2(a).

\subhead 8.5\endsubhead
We show:

(a) {\it If $B_1,B_2$ are elements of $\ck(\fg_\da)$ then 
$\{B_1,B_2\}\in\QQ(v)$ (notation of 4.4(c)).}
\nl
By 8.3, we can assume that $B_1,B_2$ are classes of $\et$-quasi-monomial 
objects. By 8.2(a) we have $f_1B_1=B'_1$, $f_2B_2=B'_2$ where $B'_1,B'_2$
represent $\e$-monomial objects and $f_1,f_2$ are nonzero elements of $\ca$.
Thus, we can assume that $B_1,B_2$ represent $\et$-monomial objects. In this 
case the result follows from 6.4.

\head 9. Examples\endhead
In this section we consider examples where $G=SL(V)$ or $Sp(V)$. We assume that $m\ge2$ and $\et=1$  hence $\da=\un1$. We write "spiral"
instead of "$1$-spiral". 
We explicitly describe the spirals and the set of blocks $\un\fT_{1}$ in both cases, and describe the map $\Ps$ in the case $G=SL(V)$. 

\subhead 9.1. Spirals for the cyclic quiver\endsubhead
We preserve the notation from 0.3. Thus we assume that $G=SL(V)$ where 
$V=\oplus_{i\in\ZZ/m}V_i$. We have an induced $\ZZ/m$-grading on 
$\fg=\fs\fl(V)$, so that $\fg_{\un1}$ is the space of all maps in 0.3(a).
In general, we have $\fg_i=\op_{j\in\ZZ/m}\Hom(V_j, V_{j+i})$.

The datum $\l\in Y_{G_{\un0},\QQ}$ is the same as a $\QQ$-grading on each 
$V_i$, i.e., $V_i=\oplus_{x\in\QQ}({}_xV_i)$ such that 
$\sum_i\sum_xx\dim({}_xV_i)=0$. Given such a $\QQ$-grading on each $V_i$, the 
corresponding spiral $\fp_*=\{\fp_N\sub\fg_{\un N}\}_{N\in\ZZ}$ takes the 
following form
$$\fp_N=\{\ph\in\fs\fl(V)|\ph({}_xV_j)\sub
\op_{x'\ge x+N}({}_{x'}V_{j+\un N}),\frl j\in\ZZ/m,x\in\QQ\}.$$
A splitting $\fm_*=\{\fm_N\sub\fg_{\un N}\}_{N\in\ZZ}$ of the spiral $\fp_*$ 
takes the form
$$\fm_N=\{\ph\in\fs\fl(V)|\ph({}_xV_j)\sub {}_{x+N}V_{j+\un N},\frl j\in\ZZ/m,
x\in\QQ\}.$$
For such a grading ${}_xV_i$ we may introduce a quiver $Q_\l$ as follows. Let 
$J_\l$ be the finite set of pairs $(i,x)\in\ZZ/m\T\QQ$ such that 
${}_xV_i\ne0$. Then $Q_\l$ has vertex set $J_\l$ and an edge 
$(i,x)\to(i+1,x+1)$ if both $(i,x)$ and $(i+1,x+1)$ are in $J_\l$.  Then 
$Q_\l$ is a disjoint union of directed chains (that is, quivers of type $A$
with exactly one source and exactly one sink). We may identify $\fm_1$ with 
the representation space of the quiver $Q_\l$ with vector space ${}_xV_i$ on 
the vertex $(i,x)\in J_\l$.

Let $B$ be the set of chains in $Q_\l$, and let 
$J_\l=\sqc_{\b\in B}({}_\b J_\l)$ be the corresponding decomposition of the 
vertex set. Let ${}_\b V:=\op_{(i,x)\in\b}({}_xV_i)$. Then we have 
$V=\op_{\b\in B}({}_\b V)$.
Let $M=e^{\fm}, M_0=e^{\fm_0}$ where $\fm=\op_N\fm_N$. Then
$M=S(\prod_{\b\in B}GL({}_\b V))$, 
$M_0=S(\prod_{(i,x)\in J_\l}GL({}_xV_i))$. The center $Z_M$ is the subgroup 
of $M$ where each factor in $GL({}_\b V)$ is a scalar matrix.

\subhead 9.2. Admissible systems for the cyclic quiver\endsubhead
Let $d$ be a divisor of $n=\dim V$. Suppose that the following hold:

(1) Each ${}_xV_i$ has dimension $\le1$.

(2) Each connected component of the quiver $Q_\l$ is a directed chain 
containing exactly $d$ vertices.
\nl
In this case, $M_0$ is a maximal torus of $G$ stabilizing each line ${}_xV_i$ 
for $(i,x)\in J_\l$. The open $M_0$-orbit $\ovsc{\fm}_1\sub\fm_1$ consists of
representations of $Q_\l$ where all arrows are nonzero (hence isomorphisms). 
The stabilizer of an element in $\ovsc{\fm}_1$ under $M_0$ is exactly $Z_M$, 
which acts by a scalar $z_\b$ on each chain $\b\in B$, such that 
$(\prod_{\b\in B}z_\b)^d=1$. We see that $\pi_0(Z_{M})\cong\mu_{d}$. For any 
primitive character $\c:\mu_d\to\bbq^*$, we have a rank $1$ $M_0$-equivariant
local system $C_\c$ on $\ovsc{\fm}_1$ on whose stalks $Z_M$ acts via $\c$. 
This is a cuspidal local system because it is the restriction of the cuspidal 
local system on the regular nilpotent orbit of $\fm$ with central character 
$\c$. Let $\tC_\c$ be the cuspidal perverse sheaf on $\fm_1$ corresponding to
$C_\c$. The system $(M,M_{0},\fm,\fm_{*},\tC_\c)$ is admissible. It is easy to 
see that any admissible system is of the form we just described.

Given such a grading $\l$, we define a function $f:B\to\ZZ/m$ such that 
$f(\b)=i$ where $(i,x)$ is the head (origin) of the chain $\b$. Each vertex 
$(i,x)\in J_\l$ lies in a unique chain $\b\in B$ whose head is of the form 
$(f(\b),x')$. Then $x-x'=y$ is an integer between $0$ and $d-1$ and 
$f(\b)+\un y=i$ in $\ZZ/m$. This implies that 
$\dim V_i=\sha\{x\in\QQ|(i,x)\in J_\l\}$ is the same as the number of 
pairs $(\b,y)\in B\T\{0,1,\cdots,d-1\}$ such that $f(\b)+\un y=i$. Choosing a 
bijection between $\{1,2,\cdots,n/d\}$ and $B$, the function $f$ may be viewed
as a function $\{1,2,\cdots,n/d\}\to\ZZ/m$ satisfying 0.7(b). 
Changing the 
bijection amounts to pre-composing $f$ with a permutation of 
$\{1,2,\cdots,n/d\}$. Summarizing the above discussion, we get a canonical 
bijection between $\un\fT_1$ and the set of equivalence classes of triples 
$(d,f,\c)$ as in 0.7(a). 

\subhead 9.3. The map $\Ps$ for the cyclic quiver\endsubhead
We preserve the notation from 9.1. 
Let $(\co,\cl)\in\ci(\fg_{\un1})$. 
For each 
element $e\in\co$, there exists a decomposition of $V$ into Jordan blocks 
$\{{}_\a W\}_{\a\in B_e}$ compatible with the $\ZZ/m$-grading in the following
sense. Each Jordan block ${}_\a W$ is a direct sum of finitely many 
$1$-dimensional subspaces indexed by $0,1,\cdots$, i.e., 
${}_\a W=({}_\a W_0)\op({}_\a W_1)\op\do$ such that 

(1)  ${}_\a W_N\subset V_{h(\a)+\un N}$ for some $h(\a)\in\ZZ/m$ (location of the 
head of the Jordan block $\a$);

(2) $e$ maps ${}_\a W_N$ isomorphically to ${}_\a W_{N+1}$ whenever $N\ge0$ and ${}_\a W_{N+1}\ne0$. 
\nl
The datum of $\{{}_\a W\}_{\a\in B_e}$ as above is the equivalent to the datum of an element $\phi\in J_{\un1}(e)$, see 2.3.  
From this we may define a quiver $Q_e$ whose 
vertex set $J_e$ consists of pairs $(\a,N)\in B_e\T\ZZ_{\ge0}$ such that 
${}_\a W_N\ne0$, and there is an edge $(\a,N)\to(\a,N+1)$. Each vertex 
$(\a,N)$ is labelled with the element $h(\a)+\un N\in\ZZ/m$. The isomorphism 
class of $Q_e$ together with the labelling by elements in $\ZZ/m$ is 
independent of the choice of $e$ in $\co$ and the choice of the Jordan block 
decomposition. Therefore we denote this labelled quiver by $Q_{\co}$, with 
vertex set $J_{\co}$ and set of chains $B_{\co}$.

Let $d'=\gcd\{|\a|\}_{\a\in B_{\co}}$ (here $|\a|$ is the number of vertices 
of the chain $\a$). Then for any $e\in\co$, there is a canonical isomorphism 
$\pi_{0}(G_{\un0}(e))\cong\mu_{d'}$. The local system $\cl$ on $\co$ 
corresponds to a character $\r$ of $\mu_{d'}$, which has order $d$
dividing $d'$  
and a unique factorization
$$\r: \mu_{d'}@>>>\mu_d@>\c>>\bbq^*$$
such that $\c$ is injective (here the first map $\mu_{d'}@>>>\mu_d$ is given 
by $z\m z^{d'/d}$). Now we define a new quiver $Q^{[d]}_{\co}$ by removing 
certain edges from each chain of $Q_{\co}$ such that each chain of 
$Q^{[d]}_{\co}$ has exactly $d$ vertices. Let $B$ be the set of chains of 
$Q^{[d]}_{\co}$; then $B$ can be identified with the set 
$\{1,2,\cdots,n/d\}$. Define $f:\{1,2,\cdots,n/d\}\cong B\to\ZZ/m$ to be the 
map assigning to each $\b\in B$ the label of its head. This way we get a 
triple $(d,f,\c)$ as in 0.7(b) whose  equivalence class is well-defined.

\proclaim{Proposition 9.4}In the case of cyclic quivers, the map 
$\Ps:\ci(\fg_{\un1})\to\un\fT_1$ 
sends $(\co,\cl)$ to the admissible system in 
$\un\fT_1$ which corresponds to the equivalence class of the triple $(d,f,\c)$
defined above under the bijection 0.7(a).  
\endproclaim
Let $e\in\co$, and let $V=\op_{\a\in B_e}({}_\a W)$, 
${}_\a W={}_\a W_0\op{}_{\a}W_1\op\cdots$ be a Jordan block decomposition, 
where ${}_\a W_N\sub V_{h(\a)+{\un N}}$ 
for $\a\in B_e,N\in\ZZ_{\ge0}$. Let 
$L$ be the Levi subgroup of a parabolic subgroup of $G$ such that $L$
stabilizes the decomposition $V=\op_{\a\in B_e}({}_\a W)$. Then $\fl=\fL L$ 
has a $\ZZ$-grading induced from the $\ZZ$-grading on each of ${}_\a W$. In 
particular, $\fl_1$ is the space of representations of the quiver $Q_e$. The 
system $(L, L_0, \fl,\fl_*)$ is the system 
$(\ti L^\ph,\ti L^\ph_0,\ti\fl^\ph,\ti\fl^\ph_{*})$ attached to some 
$\ph\in J_{\un1}(e)$ as in 2.9. 
Then $e$ is in the open $L_0$-orbit $\ovsc{\fl}_1$ of $\fl_1$, which is contained 
in the regular nilpotent orbit of $\fl$. 

Let ${}_\a L=SL({}_\a W)$ 
be the subgroup of $L$ which acts as 
identity on all blocks ${}_{\a'}W$ for $\a'\ne\a$. Then 
${}_\a\fl=\fL({}_\a L)$ carries a $\ZZ$-grading compatible with that on $\fl$.
For each interval $[a,b]\sub\ZZ_{\ge0}$, let ${}_\a W_{[a,b]}\sub{}_\a W$ be 
the direct sum of ${}_\a W_N$ for $a\le N\le b$. We decompose 
${}_\a W$ into $|\a|/d$ parts each of dimension $d$:
$$ {}_{\a}W=\op_{j=1}^{|\a|/d}({}_{\a}W_{[(j-1)d,jd-1]}).\tag a$$
Let ${}_\a M\subset{}_\a L$ be the subgroup stabilizing the decomposition (a).
Then the Lie algebra ${}_\a\fm$ of ${}_\a M$ inherits a $\ZZ$-grading from 
that of ${}_\a\fl$, and the open orbit ${}_\a\ovsc\fm_1$ carries a local 
system ${}_\a C_\c$ corresponding to the character $\c$ of 
$\mu_d\cong \pi_0(Z({}_{\a}M))$. 
Let ${}_\a\tC_\c$ be the cuspidal perverse sheaf on ${}_\a\fm_1$ 
corresponding to ${}_\a C_\c$.
Define a parabolic subalgebra 
${}_\a\fq\sub{}_\a\fl$ to be the stabilizer of the filtration 
${}_\a W_{[|\a|-d,|\a|-1]}\sub {}_\a W_{[|\a|-2d,|\a|-1]}\sub\cdots\sub
 {}_\a W={}_\a W_{[0,|\a|-1]}$. Then ${}_\a\fq$ is compatible with the 
$\ZZ$-grading on ${}_\a\fl$ and ${}_{\a}\fm$ is a Levi subalgebra of 
${}_\a\fq$. The induction 
$$\ind^{{}_\a\fl_1}_{{}_\a\fq_1}({}_\a \tC_\c)$$ 
restricted to ${}_\a\ovsc\fl_1$ is isomorphic to $\cl|_{{}_\a\ovsc\fl_1}$, because the map $c$ in 1.3 (applied to ${}_{\a}\fl, {}_\a\fq, {}_\a\fm$ in place 
of $\fh,\fp,\fl$) is an isomorphism when restricted to ${}_\a\ovsc\fl_1$. 
Therefore the middle extension of $\cl|_{{}_\a\ovsc\fl_1}$ to 
$\fl_1$ appears as a direct summand of $\ind^{{}_{\a}\fl_1}_{{}_\a\fq_1}({}_\a C_\c)$.
Therefore, under the map 
defined in 1.5(b), the image of $({}_\a\ovsc\fl_1,\cl|_{{}_\a\ovsc\fl_1})$ 
is 
$$({}_\a M,{}_\a M_0,{}_\a\fm,{}_\a\fm_1,{}_\a\tC_\c).$$

Let ${}_\a\ti M\sub GL({}_\a W)$ be the stabilizer of the decomposition (a).
Let $M=S(\prod_{\a\in B_e}({}_\a\ti M))\sub L$ with Lie algebra 
$\fm\sub\op({}_\a\ti\fm)$ and the induced $\ZZ$-grading from each 
${}_\a\ti\fm=\fL({}_\a\ti M)$. The open $M_0$-orbit on $\fm_1=\op({}_\a\fm_1)$
is $\ovsc\fm_1=\prod({}_\a\ovsc\fm_1)$. Let $C_\c=\boxtimes({}_\a C_\c)$ on 
$\ovsc\fm_1$. Let $\tC_\c$ be the cuspidal perverse sheaf on $\fm_1$ 
corresponding to $C_\c$.
 By the compatibility of the assignment in 1.5(b) with direct 
products, in the situation $H=L$, the pair 
$(\ovsc\fl_1,\cl|_{\ovsc\fl_1})$ maps to $(M,M_0,\fm,\fm_{*},\tC_\c)$. 
Therefore, $(M,M_0,\fm,\fm_*,\tC_\c)$ 
is the admissible system attached to 
$(\co,\cl)$ through the procedure in 2.9. By 9.2, the admissible system 
$(M,M_0,\fm,\fm_*,\tC_\c)$ 
corresponds to the triple $(d,f,\c)$ defined in 9.3 
before the statement of this proposition. This finishes the proof.

\subhead 9.5. The symplectic quiver\endsubhead
Let $V$ be a finite-dimensional vector space over $\kk$ with a
nondegenerate symplectic form $\o$. Assume that $m$ in 0.1 is even. 
Let $\ti\fS_m=\{j; j=k/2; k=\text{ an odd integer}\}$ and let $\fS_m$ be the 
set of equivalence classes for the relation $\si$ on $\ti\fS_m$ given by
$j\si j'$ if $j-j'\in m\ZZ$. Note that the involution $j\m-j$ of $\ti\fS_m$ 
induces an involution of $\fS_m$ denoted again by $j\m-j$. 

For any $N\in\ZZ$, the map $j\m N+j$ of $\ti\fS_m$ onto itself induces a 
map of $\fS_m$ onto itself which depends only on $\un N$ and is denoted by
$j\m\un N+j$. 

The set $\fS_m$ consists of $m$ elements represented by
$$\{\fra{1}{2},\fra{3}{2},\cdots,\fra{m-1}{2},\fra{m+1}{2},\cdots,
m-\fra{1}{2}\}.$$

Consider a grading on $V$ indexed by $\fS_m$:
$$V=\bigoplus_{j\in\fS_m}V_j,\tag a$$
such that $\o(V_j,V_j')=0$ unless $j'=-j$ (as elements of $\fS_m$). Using the 
symplectic form, for $j\in\fS_m$ we may identify $V_j$ with the dual of 
$V_{-j}$.

We assume that $G=Sp(V)$ 
and that the $\ZZ/m$-grading of 
$\fg=\fs\fp(V)$ is given by  
$$\fg_i=\{\ph\in\fs\fp(V)|\ph(V_j)\sub V_{i+j},\frl j\in\fS_{m}\},\qua\frl 
i\in\ZZ/m.\tag b$$
In particular, an element $\ph\in\fg_{\un1}$ 
is a collection of maps 
$\ph_i:V_{i-\fra{1}{2}}\to V_{i+\fra{1}{2}}$, $i\in\ZZ/m$, which can be 
represented by a cyclic quiver
$$\CD
V_{-\fra{1}{2}}@<\ph_{-1}<<V_{-\fra{3}{2}}@<<< \do@<\ph_{\fra{m}{2}+1}<<
V_{\fra{m+1}{2}}\\
@V\ph_0VV         @.                  @.                @A\ph_{\fra{m}{2}}AA\\
V_{\fra{1}{2}}@>\ph_1>>V_{\fra{3}{2}}@>>> \do @<\ph_{\fra{m}{2}-1}<<
V_{\fra{m-1}{2}}\\
\endCD$$    
The condition $\ph\in\fs\fp(V)$ becomes that
$$\ph_{-i}=-\ph^*_i, \quad \frl i\in\ZZ/m.\tag c$$
Here $\ph^*_i:V^*_{i+\fra{1}{2}}\to V^*_{i-\fra{1}{2}}$ is the adjoint of 
$\ph_i$, which can be viewed as a map $V_{-i-\fra{1}{2}}\to V_{-i+\fra{1}{2}}$
under the identifications $V^*_{i+\fra{1}{2}}\cong V_{-i-\fra{1}{2}}$, 
$V^*_{i-\fra{1}{2}}\cong V_{-i+\fra{1}{2}}$ using the symplectic pairing. In 
particular, for $i=0$, 
$\ph_0:V_{-\fra{1}{2}}=V^*_{\fra{1}{2}}\to V_{\fra{1}{2}}$ can be viewed as a 
vector $\ph_0\in V^{\ot2}_{\fra{1}{2}}$. The condition (b) for $i=0$ is 
equivalent to saying that $\ph_0\in\text{\rm Sym}^2(V_{\fra{1}{2}})$. 
Similarly, we may view $\ph_{\fra{m}{2}}$ as a vector in 
$V^{\ot2}_{\fra{m+1}{2}}$, and the condition (c) for $i=\fra{m}{2}$ is 
equivalent to saying that 
$\ph_{\fra{m}{2}}\in\text{\rm Sym}^2(V_{\fra{m+1}{2}})$.  

We call a representation of the quiver above in which $V_{-j}=V^*_j$, and (c) 
holds a {\it symplectic representation}. In other words, $\fg_{\un1}$ 
is the space of symplectic representations of the quiver above.

We have $G_{\un0}\cong\prod_{\fra{1}{2}\le j\le\fra{m-1}{2}}GL(V_j)$, where 
$GL(V_j)\cong GL(V_{-j})$ acts diagonally on both $V_j$ and 
$V_{-j}=V_{m-j}=V^*_j$.

\subhead 9.6. Spirals for the symplectic quiver\endsubhead
Each element $\l\in Y_{G_{\un0},\QQ}$ is 
the same datum as a $\QQ$-grading on each $V_j$, $j\in\fS_m$, i.e., 
$V_j=\oplus_{x\in\QQ}({}_xV_j)$ such that under the symplectic form $\o$, 
$\o({}_xV_j, {}_{x'}V_{-j})=0$ unless $x+x'=0$. Then ${}_{-x}V_{-j}$ can be 
identified with the dual of ${}_xV_j$ for all $(j,x)\in\fS_m\T\QQ$.  The 
spiral $\fp_{*}$ associated to this grading is
$$\fp_N=\{\ph\in\fs\fp(V)|\ph({}_xV_j)\sub\op_{x'\ge x+N}({}_{x'}V_{j+\un N }),
\frl j\in\fS_m,x\in\QQ\}.$$ 
A splitting $\fm_*$ of the spiral $\fp_*$ takes the form
$$\fm_N=\{\pi\in\fs\fp(V)|\ph({}_xV_j)\sub {}_{x+N}V_{j+\un N}, \frl 
j\in\fS_m, x\in\QQ\}.$$ 
To each such grading, we may attach a quiver $Q_\l$ as we did for the cyclic 
quiver (since the symplectic quiver is a special case of a cyclic quiver). 
There is an involution on $Q_\l$ sending $(j,x)\in J_\l$ to 
$(-j,-x)\in J_\l$. This involution stabilizes at most two chains $Q'_\l$ and $Q''_\l$
of $Q_\l$. The set of vertices of $Q'_\l$ (possibly empty) is 
$J'_\l:=\{(x,x)|{}_xV_x\ne0\}\subset J_\l$. The set of vertices of 
$Q''_\l$ (possibly empty) is $J''_\l:=\{(x-\fra{m}{2},x)| {}_xV_{x-\fra{m}{2}}\ne0\}\sub J_\l$.

\subhead 9.7. Admissible systems for the symplectic quiver\endsubhead
Suppose that the following hold:

(1) For each $(j,x)\in J-(J'_\l\sqcup J''_\l)$, we have $\dim {}_xV_j=1$. 

(2) The chains in $Q_\l$ other than $Q'_\l$ and $Q''_\l$ all consist of a 
single vertex.

(3) Let $\sha J'_\l=2a'$ for some $a'\in\ZZ_{\ge0}$. When $a'>0$, 
$(-a'+\fra{1}{2},-a'+\fra{1}{2})$ is the head of $J'_\l$ and 
$(a'-\fra{1}{2},a'-\fra{1}{2})$ is the tail. Then 
$\dim{}_xV_x=a'+\fra{1}{2}-|x|$ for all $(x,x)\in J'_\l$.

(4) Let $\sha J''_\l=2a''$ for some $a''\in\ZZ_{\ge0}$. When $a''>0$, 
$(-a''-\fra{m-1}{2},-a''+\fra{1}{2})$ is the head of $J''_\l$ and 
$(a''-\fra{m+1}{2},a''-\fra{1}{2})$ is the tail. Then 
$\dim {}_xV_{x-\fra{m}{2}}=a''+\fra{1}{2}-|x|$ for all 
$(x-\fra{m}{2},x)\in J''_\l$.
\nl  
Under these assumptions, $\fm_1=\fm'_1\op\fm''_1$, where $\fm'_1$ is the space
of representations of the quiver $Q'_\l$ with dimension vector 
$\dim {}_xV_x=a'+\fra{1}{2}-|x|$ and satisfying the duality condition 
$\psi_i=-\psi^*_{-i}$ (where $\psi_i:{}_{i-\fra{1}{2}}V_{i-\fra{1}{2}}\to 
{}_{i+\fra{1}{2}}V_{i+\fra{1}{2}}$)
for all $i\in\{-a'+1,\cdots,a'-1\}$. 
Similarly, $\fm''_1$ is the space of representations of the quiver $Q''_\l$ 
with dimension vector $\dim {}_xV_{x-\fra{m}{2}}=a''+\fra{1}{2}-|x|$ and 
satisfying 
the duality condition $\psi_i=-\psi^*_{-i}$. The open $M_0$-orbit 
$\ovsc{\fm}_1$ consists of those representations of $Q'_\l$ and $Q''_\l$ where
each arrow has maximal rank (either injective or surjective). 

Let $V'=\op {}_xV_x$ and $V''=\op {}_xV_{x-\fra{m}{2}}$. Let 
$V^{\dagger}=\op_{(j,x)\notin J'_\l\cup J''_\l}({}_xV_j)$. 
Then we have 
$V=V'\op V''\op V^{\dagger}$. This decomposition is preserved by $M$, and 
$M\cong Sp(V')\T Sp(V'')\T T^{\dagger}$, 
where $T^{\dagger}$ is the maximal torus in $Sp(V^{\dagger})$ 
stabilizing each line ${}_xV_j\sub V^{\dagger}$. 
The center $Z_M$ is isomorphic to 
$\{\pm1\}\T\{\pm1\}\T T^{\dagger}$ under this decomposition. The stabilizer of
a point in $\ovsc\fm_1$ under $M_0$ is exactly $Z_M$. Let $C$ be the rank one 
local system on $\ovsc\fm_1$ on whose stalks $\p_0(Z_{M})$ acts nontrivially 
on both factors of $\{\pm1\}$. Then $C$ is cuspidal because it is the 
restriction of the unique cuspidal local system on $\fm$. Let $\tC$ be the
cuspidal perverse sheaf on $\fm_1$ defined by $C$. The system 
$(M,M_0,\fm,\fm_*,\tC)$ is admissible. Moreover, any admissible system is of 
this form. Under $G_{\un0}$-conjugacy, the only invariant of an admissible 
system is given by the numbers $a'$ and $a''$. Since 
$\dim V'_j+\dim V''_j\le\dim V_j$, we have the following inequality for all 
$j\in\fS_m$  
$$\align&
\dim V_j\ge \sha\{-a'+\fra{1}{2}\le x\le a'-\fra{1}{2}|x\equiv j\mod m\ZZ\}
\\&
+\sha\{-a''+\fra{1}{2}\le x\le a''-\fra{1}{2}|x\equiv j+\fra{m}{2}\mod m\ZZ\}.
\tag a\endalign$$
To summarize, we have a natural bijection
$$\un\fT_1\lra\{(a',a'')\in\ZZ_{\ge0}\T\ZZ_{\ge0}\text{ satisfying (a) for 
all }j\in\fS_m\}.\tag b$$
The map $\Psi:\ci(\fg_{\un1})\to\un\fT_1$ 
for the symplectic quiver as well as 
other graded Lie algebras of classical type  will be described in a sequel 
to this paper using the combinatorics of symbols.


\widestnumber\key{ABC}
\Refs
\ref\key{\KL}\by D.Kazhdan and G.Lusztig\paper Proof of the Deligne-Langlands
conjecture for Hecke algebras\jour Invent.Math.\vol87\yr1987\pages153-215
\endref
\ref\key{\KO}\by B.Kostant\paper The principal three dimensional subgroup and 
the Betti numbers of a complex simple Lie group\jour Amer.J.Math.\vol81\yr1959
\pages973-1032\endref
\ref\key{\INTERS}\by G.Lusztig\paper Intersection cohomology complexes on a 
reductive group\jour Inv. Math.\vol75\yr1984\pages205-272\endref
\ref\key{\CSI}\by G.Lusztig\paper Character Sheaves, I\jour Adv. Math.\vol56
\yr1985\pages193-237\endref  
\ref\key{\CSII}\by G.Lusztig\paper Character Sheaves, II\jour Adv. Math.\vol57
\yr1985\pages226-315\endref  
\ref\key{\GRA}\by G.Lusztig\paper Study of perverse sheaves arising from graded
Lie algebras\jour Adv.Math.\vol112\yr1995\pages147-217\endref
\ref\key{\ANTIORB}\by G.Lusztig\paper Study of antiorbital complexes\jour
Contemp.Math.\vol557\yr2011\pages259-287\endref
\ref\key{\RIR}\by L.Rider and A.Russell\paper Perverse sheaves on the nilpotent
cone and Lusztig's generalized Springer correspondence\jour arxiv:1409.7132
\endref
\ref\key{\ST}\by R.Steinberg\book Endomorphisms of linear algebraic groups
\bookinfo Memoirs of Amer. Math.Soc. 80\publ Amer.Math.Soc.\yr1968\endref
\ref\key{\VI}\by E.B.Vinberg\paper The Weyl group of a graded Lie algebra\jour
Izvestiya Akad. Nauk S.S.S.R.\vol40\yr1976\pages488-526\endref
\endRefs
\enddocument